%% file: guideblob3.tex.vcleaned.tex
\documentclass[12pt]{article} 
\usepackage{latexsym} \usepackage{amssymb}
\usepackage{epic} \usepackage{eepic} 
\usepackage[all]{xy}
\usepackage{amscd} \usepackage{graphicx}
\input amssym.def
\input amssym
\begin{document}
\newtheorem{theo}{Theorem}      \newtheorem{cor}{Corollary}[theo]   
\newtheorem{de}{Definition}     \newtheorem{pr}{Proposition} 
\newtheorem{co}{Corollary}[pr]  \newtheorem{rem}{Remark} 
\newtheorem{lem}{Lemma} 
\newenvironment{nullenvir}{}{}
\newenvironment{proof}[1][Proof:]{\begin{nullenvir} \noindent {\em #1}  \footnotesize}%
                      {\end{nullenvir} \hfill \Qed}
\newenvironment{smalleq}{\begin{nullenvir} \noindent  \footnotesize}%
                      {\end{nullenvir} \hfill \Qed}
\newcounter{minidef}[section]
\renewcommand{\theminidef}{\thesection.\arabic{minidef}}
\newcommand{\mdef}{\refstepcounter{minidef} 
\medskip \noindent ({\bf \theminidef}) }

\newcommand{\mpr}[1]{\refstepcounter{minidef} 
\medskip \noindent ({\bf \theminidef}) {\bf Proposition}. {\em #1}}

\newcommand{\UU}{\underline{\sqcup}}  \newcommand{\UUU}{\sqcup}  
\newcommand{\e}{\epsilon}        
\newcommand{\lam}{\lambda}  
\newcommand{\la}{\langle}        \newcommand{\ra}{\rangle}
\newcommand{\ha}{\#}             \newcommand{\rmap}{\rightarrow}
\newcommand{\aaa}{\alpha}        \newcommand{\ab}{\alpha,\beta}
\newcommand{\aab}{a(\ab )}       \newcommand{\Weyl}{Weyl }
\newcommand{\HH}{H \!\!\! I}              
\newcommand{\C}{\mathbb C }
\renewcommand{\Re}{\mathbb R }
\newcommand{\Z}{\mathbb Z }
\newcommand{\N}{\mathbb N }
\newcommand{\Q}{\mathbb Q }
\def\H{\Bbb H }
\def\sixrm{} \def\sevrm{} \def\twlrm{} \def\tenrm{}
\def\nset(#1){ \{ #1 \}_{ \underline{n} }}
\def\Sym(#1){\Sigma(#1)}                  
\def\Sy(#1){\Sigma_{#1}}                  
\def\sym(#1){\mbox{\LARGE s}(#1)}       
\def\sy(#1){\mbox{\LARGE s}_{#1}}       
\def\Sm(#1){\Sigma_{#1}}                  
\def\Ee(#1){{\bf E}_{#1}}                 
\def\Eee(#1){{\bf E}_{\{ #1 \}_{\underline{n}}}}  
\def\Ss(#1){{\bf S}_{#1}}                 
\def\Sss(#1){{\bf S}_{\{ #1 \}_{\underline{n}}}}  
\def\ka(#1){\kappa_{#1}}                  
\def\SS(#1){{\cal S}_{#1}}                
\def\ul(#1){_{\underline{#1}}}            
\def\Tm{{\cal T}}                         
\def\ul(#1){_{\underline{#1}}}            
\newcommand{\cs}{\C \, \sym(n)}           
\newcommand{\Gg}{{\cal G}}                
\newcommand{\com}{\bullet}                
\def\Ai(#1){ A^{ #1 \cdot } }             
\def\Aij(#1,#2){ A^{ #1  #2 } }           
\def\One{\mbox{\bf $1 \!\!\!\! 1$}}                  
\def\Gg{{\cal G}}                         
\def\Fg{{\cal F}}                         
\def\Vq{V_Q^{\otimes n}}                  
\def\Wm(#1){{\cal S}_{#1}}                
\def\wm(#1,#2){{}_{#1}{\cal S}_{#2}}      
\def\Bp{B_p}                              
\def\Bb(#1){B_p[#1]}                      
\def\Pp(#1){P_n[#1]}                      
\def\Ps(#1){P_n[#1] \! /}                 
\def\Ph{\hat{P}}                          
\def\Is(#1){\sim^{#1}}                    
\def\Ff(#1,#2){F_{#1}^{(#2)}}             
\def\Ef(#1,#2){E_{#1}^{(#2)}}             
\def\bbc(#1){((\beta_1)(\beta_2)...(\beta_{#1}))}     
\def\Ln{{\cal L}_{n}}                     
\def\smap{s}                              
\def\tmap{t}                              
\def\pmap{\psi}                           
\def\Amap(#1){{\cal A}_{#1}}              
\def\Rr(#1){R_{#1}}                       
\def\Cr(#1){C_{#1}}                       
\def\Zz{\zeta}                            
\def\Ww(#1){\mbox{\LARGE $w$}_{#1}}       
\def\choo(#1,#2){ \left( \begin{array}{c} #1 \\ #2 \end{array} \right) }
\def\Qed{$\Box$}                          
\def\staq(#1){\stackrel{#1}{=}}           
\def\ra{\rightarrow}
\def\ses(#1,#2,#3){0\ra #1 \ra #2 \ra #3 \ra 0}  
\newcommand{\beq}{\begin{equation} }
\def\eql(#1){ \begin{equation} \label{#1} 
}
\newcommand{\eq}{\end{equation} }
\def\lab(#1){\label{#1}
}
\def\prl(#1){ \begin{pr} \label{#1} 
}
\def\mat{ \left( \begin{array} }    \def\tam{ \end{array}  \right) }

\def\del(#1){ \begin{de} \label{#1} 
}
\def\starr(#1){ \stackrel{ #1 }{\longrightarrow} }
\def\Hnq{H_n(q)}
\def\Hn{H_n}    
\def\A{{\cal A}}
\def\mod{\mbox{mod} }
\def\Res(#1,#2){\mbox{Res}^{#1}_{#2}}
\newcommand{\WRes}{\mbox{Res}}
\newcommand{\WInd}{\mbox{Ind}}
\def\blob{ the blob algebra }
\def\TL{ Temperley-Lieb }
\def\MS{ \cite{MartinSaleur}}
\def\DD(#1,#2,#3){\Delta_{#1}({#3} #2 )}
\def\RR(#1,#2,#3){R_{#1}({#3} #2 )}
\def\ignore#1{}
     \newcommand{\NB}{NB}
\newcommand{\xfig}[1]{\input{./xfig/#1.eepic}}
\newcommand{\xfigla}[1]{\input{./xfig/#1.latex}}
\newcommand{\xfigeps}[1]{\includegraphics{./xfig/#1.eps}}
\newcommand{\mprop}[1]{\noindent{\bf Proposition.}{\em #1}}
\newcommand{\Aring}{{\mathcal A}}  
\newcommand{\Aquot}{{\mathcal K}}  

\def\TL{T}
\def\braid#1{A_{#1}\mbox{--braid}}
\def\abraid#1{\hat{A}_{#1}\mbox{--braid}}
\def\bbraid#1{{B}_{#1}\mbox{--braid}}
\def\Nset#1{\underline{#1}}

\[ \]
\begin{center}
{\Large {\bf
Generalized Blob Algebras and Alcove Geometry}}
\\
\vspace{.4in}
   Paul P Martin 
   \footnote{\footnotesize
                Mathematics Department,
                City University, 
                Northampton Square, London EC1V 0HB, UK.}
and
   David Woodcock$^{1}$
\end{center}
\newcommand{\MartinWoodcockLevy}{MartinWoodcockLevy}
\newcommand{\pre}{}
\newcommand{\MartinWoodcock}{MartinWoodcock}
\section{Introduction}\label{s1}

\newcommand{\aH}[1]{H(#1)} 
\newcommand{\AH}{\mbox{H}}        
\newcommand{\LAM}{\Lambda}
\renewcommand{\DD}{\Delta}
\newcommand{\BB}{{\sf B}}
\def\PPP(#1,#2){P^{#1}_{#2}}%
\def\PP(#1){P_{#1}}%
\def\EEE(#1,#2,#3){e^{#1 #2}_{#3}}%
\def\EEEu(#1,#2,#3){\mbox{\bf e}^{#1 #2}_{#3}}%
\def\EE(#1){e^{#1 }}%
\newcommand{\Alg}{{\sf A}}
\newcommand{\Ide}{e}     
\newcommand{\weights}{\Lambda}
\newcommand{\AKweights}{\Lambda}
\newcommand{\ym}{m}
\newcommand{\End}{\mbox{End}}
\newcommand{\Bratteli}{Bratteli}
\newcommand{\Pglob}{{\mathfrak P}}
\newcommand{\Alc}{{\mathfrak A}}
\newcommand{\Bbas}{{\mathfrak B}}

Soergel has given a beautiful procedure \cite{Soergel97A,Soergel97B} 
for analysing tilting modules 
for quantum groups at roots of unity 
through parabolic Kazhdan--Lusztig polynomials. 
The procedure itself may be applied formally to
an alcove geometry, {\em without} reference to representation theory. 
Hence it may be applied, in principle, in cases which are beyond the scope of  
Soergel's proof of representation theoretic interpretation. 
It is interesting then to try to find
algebras for which the resultant combinatorial data 
{\em has}
a representation theoretic interpretation,
even though Soergel's proof is not applicable. 
The output of the usual procedure
in type--$A$ may be mapped by Ringel duality \cite{Donkin98agc} 
to the content of projective modules for certain quotients
of ordinary Hecke algebras. (There it may be understood in terms of
idempotent decompositions of 1 \cite{\MartinWoodcock98b}.) 
This leads to a determination of decomposition numbers for standard
modules of the Hecke algebras themselves. 
Here we consider generalising the implementation of the 
procedure on this Ringel dual side. 
We do this by constructing generalized Hecke algebra quotients which
(mildly) generalize the usual role of alcove geometry. 

 
One example where the formal procedure gives the correct decomposition
numbers is
the blob algebra $b_n$ \cite{MartinSaleur94a,\MartinWoodcock2000\pre} 
(a certain two parameter {\em affine} Hecke algebra quotient). 
We demonstrate the procedure for this example in \S\ref{decomp nos b} below. 
There is a set of key properties of $b_n$ (see \S\ref{blob def}), which 
it has in common with the ordinary Hecke algebra quotients 
(see \S\ref{affine hecke}), which 
may serve to explain the phenomenon.  
In this paper we discuss generalisations of $b_n$ which also
possess these properties.  

To generalise $b_n$ suitably  we first place it 
in the context of affine/cyclotomic Hecke and Ariki--Koike--Levy algebras 
\cite{BroueMalle93,ArikiKoike94,MartinLevy94}
(although these are {\em not} themselves the generalisations we require). 
The study of these algebras is interesting both abstractly
and also since they are useful in studying solutions to the reflection
equation in integrable statistical mechanics (see \cite{MartinLevy94}
for references). 
This parallels the role of ordinary Hecke algebras in solving the
Yang--Baxter equation. 
In both cases the `physical' representation theory focuses
attention on specific quotients, and implies that decomposition number
data should be organized in a certain specific way. In the ordinary
Hecke case this is the `Soergel' rather than the `LLT way'
\cite{LascouxLeclercThibon96} 
(as complete data sets these are equivalent 
but computationally they are not \cite{MartinWoodcock98b}). 
Thus while the algebras we shall
construct have representation theory which is accessible in principle
by LLT methods \cite{Ariki96}, this does not remove the need for a
generalized (dual) Soergel approach. 
In \S\ref{affine hecke} we discuss quotients of affine Hecke algebras 
generalizing $b_n$ which, like $b_n$, realize 
certain key ingredients of this approach 
--- in particular they possess  a weight space. 
In \S\ref{linkage} we imbue this space with an alcove structure and
verify a linkage principle\cite{Jantzen87}. 
These generalized algebras 
are quotients of cyclotomic Hecke algebras by certain primitive
and central idempotents. 
In \S\ref{s3} we show how the simplest non--trivial such quotient may
effectively be {\em identified} with $b_n$. 
The rest of the paper discusses outstanding technical issues
in showing  
the validity of the generalized Soergel procedure for 
the generalized algebras 
(the primitive and central idempotents of the cyclotomic Hecke
algebras are computed in \S\ref{s gen d}). 


One property of the ordinary Hecke quotients {\em not} possessed by these
generalisations is the defining representation on
`tensor space', realizing Ringel
duality 
\cite{Donkin93} with a
quantum group quotient (cf. 
\cite{DipperJamesMathas99,ArikiTerasomaYamada95,SakamotoShoji99}). 
Such a faithful tensor space representation 
is not manifestly necessary for our purpose, but would be very useful. 
In \S\ref{s4} we address this problem,  culminating
in the construction of some 
intriguing new concrete representations of $b_n$ which are candidates. 
(The serendipitous constructions of a number of other interesting new
representations 
of $b_n$ and its generalisations 
are outlined in \S\S\ref{other reps}--\ref{last s}.) 


We will argue that the representation theory of these algebras, while
{\em containing} that of ordinary Lie theoretic objects, is in a sense
more simply described. Given that the representation theory of $b_n$
itself is known for $q$ root of unity in arbitrary charactersitic
\cite{CoxGrahamMartin01}, the possibility that open questions in
ordinary Lie representation theory may be accessible by this route
makes these algebras particularly interesting objects for study. 

\subsection{Alcove geometry and decomposition numbers for $b_n$}
\label{decomp nos b}\label{decnos}
We will recall the definition \cite[\S2]{MartinWoodcock2000} of the
blob algebra  $b_n=b_n(q,\ym)$ in \S\ref{blob def}. 
The decomposition numbers 
in the `doubly critical' case 
($q$ a primitive $l^{th}$ root of 1; 
$\ym$ an integer, $|\ym| < l$ --- see \S\ref{s3}) in
characteristic zero are determined in \cite[\S9]{MartinWoodcock2000}
by algebraic methods. 
A formal application of Soergel's procedure to this case works as
follows. 

First recall, quite generally, that a
Euclidean space with reflection hyperplanes
removed has a set  (${\Alc}$, say) of connected components, 
called {\em alcoves} \cite[Ch.6]{Jantzen87}\cite{Bourbaki81}.  
For $s$ a reflection hyperplane and $B$ an alcove 
we denote by $Bs$ the image of $B$ in $s$. 
Each nonempty intersection of the closure of an alcove with a 
hyperplane is called a {\em wall} of the alcove 
(here we will confuse each
such wall with the hyperplane containing it). 
We make no assumption about the relation of hyperplanes to the origin
$0$, except that the origin lies in the interior of an alcove,
called the fundamental alcove and denoted $A^0$. 
Define $|B|$ as the number of hyperplanes between $B$ and 0
(so $|A^0|=0$). 

In algebraic Lie theory one starts with a set of ordinary
(non--affine) reflections generating the ordinary Weyl group. 
The ($q$)--group weight space is the 
underlying Euclidean space with its origin $\rho$--shifted \cite{Jantzen87}. 
In particular, even when an affine reflection
is added the origin is at a fixed position at the base of the dominant
region. In our case there is effectively no {\em ordinary} Weyl group and no 
dominant region, that is to say, 
the placement of {\em all} hyperplanes is
controlled by parameters of the algebra. Thus 
the weight space for 
the blob algebra,
just as for $sl_2$, 
is the Euclidean space 
associated to the $A_1$ Coxeter system, i.e. it is effectively
$\Re$ \cite[\S6]{MartinWoodcock2000}. 
However now, cf. $sl_2$, {\em all} integral weights are dominant, 
that is to say, simple modules may be indexed by $\Z$ 
(we shall explain this, 
in the context of our generalised construction, shortly). 
In the $b_n$ case, the reflection hyperplanes  are just points, 
and those generating the affine Weyl group lie at $-\ym$ and $l-\ym$. 
(And an alcove $B \in {\Alc}$ is a connected component of $\Re$
with the reflection points removed.) 
A reflection is `upward' if $|B|<|Bs|$ 
(cf. the usual $sl_2$ situation). 

For each alcove $A$ one defines a map 
$$n_A : {\Alc} \rightarrow \Z[v]$$
where $v$ is a formal parameter, as follows. 
(For simplicity we ignore features of Soergel's procedure which do not
arise in our case.)
Firstly $n_A(A)=1$ and $n_A(B) \neq 0$ implies $|B| \leq |A|$. 
Note that $n_{A^0}$ is determined immediately by this, 
and proceed inductively on $>$. 
For each alcove $A$ there will be a wall $s$ of $A$ such that $|As|>|A|$. Each
alcove $B$ has one wall which is in the affine Weyl orbit of $s$, and we
will write $B.s$ for the image of $B$ in that wall (thus $As=A.s$). 
Then, with $n_A$ known, set
$$n'_{As} (B.s) = 
\left\{ \begin{array}{ll} 
n_A(B) + v^{-1} n_{A}(B.s)  \;\;\;\;\; & |B.s| > |B| \\
v^{-1} n_A(B) + n_{A}(B.s) & |B.s| < |B|
\end{array}\right.$$
and define $n_{As}$ by 
$$n_{As}(C) = n'_{As}(C) 
           - \sum_{B \; : \; |B|<|As|} n'_{As}(B) |_{v=0} \; n_{B}(C) .$$ 
(That this procedure is well defined is not trivial \cite{Soergel97A}.)
 
\newcommand{\Deltax}[2]{\Delta(#1)}
\newcommand{\Lx}[2]{L(#1)}         
Evaluating $n_A(B)|_{v=1}$,  
this construction is (formally) computing the standard module
content of tilting modules in a Ringel dual algebra. 
(Any two weights which are in different affine Weyl orbits \cite{Jantzen87}
are in different blocks \cite[Ch.1]{Benson95}. 
Thus each block intersects each alcove in at most 1 weight
and, fixing a block, it is the modules with these weights which
$n_A(B)|_{v=1}$ describes.)
The corresponding data for $b_n$ is, 
in effect, the standard module content of projective modules. 
By reading by column
instead of row, as it were, 
we convert this to the simple module content of standard
modules \cite{Donkin93}
(truncation to a finite column interval, such as that
pictured in the example which follows, 
corresponds to localisation to some finite $n$ --- see
ingredient~I2 below). 
As noted, this construct is entirely formal, however in fact

\mprop{ For $\lambda \in A$ and $\mu
\in B$ the $b_n$ standard composition multiplicity
$$ [ \Deltax{\mu}{n} : \Lx{\lambda}{n} ] 
        = \left\{ \begin{array}{ll} 
             n_A(B)|_{v=1} & \mbox{$\mu$ in the affine Weyl orbit of
               $\lambda$; $|\lambda| \leq n$} \\
                         0 & \mbox{otherwise} \end{array} \right. . $$  
Further, the power of the formal parameter $v$ determines the Loewy layer.}

{\em Proof:} 
The composition multiplicity data is determined in \cite[\S9]{MartinWoodcock2000}. 
A table illustrating the 
computation of the polynomials $n_A(B)$ is as follows.
\[
\input{./xfig/pKL.eepic}
\]
The row position in the table gives $A$ and the column position $B$; 
and $n_{A^0}(A^0)$ is shaded. 
\footnote{The rows and columns, and hence the alcoves, are labelled using $\Z$. 
These labels should not be confused with points in the underlying
space $\Z$ (each alcove contains $l-1$ such points) or with weights.} 
The table is complete for the rows shown, except the top and bottom rows. 
For the top row, 
the arrows within the table illustrate the contributions to the
$n'_{As}$ from a particular $n_A(B)$ (the shaded lines are the relevant
walls for reflection in each case). The arrow outside the table
illustrates a required subtraction to obtain $n_A$ in the bottom row 
(which is complete except for this subtraction). 
The pattern is clear, and 
one sees immediately that the formal procedure reproduces the
multiplicity and layer data. \Qed


The case $b_n(q,1)$ contains the ordinary Ringel dual, 
$\End_{U_qsl_2}(V^{\otimes n}_2)$,
as a quotient. Note that the ordinary Soergel procedure is 
embedded in this version
(in the `dominant' region of case $m=1$) accordingly. 
 

The `idempotent splitting' analysis described in
\cite{MartinWoodcock98b} 
applies in principle in this situation, giving a heuristic 
explanation of why Soergel's procedure is relevant here. 
Following this paradigm 
there are a set of natural generalisations 
for which an analogous method should work. 



The ingredients
are (in precis, see also 
\cite{\MartinWoodcock98b,MartinWoodcockLevy2000}):
\begin{enumerate}
\item[I1] A tower of unital algebras $\Alg_n \subset \Alg_{n+1}$ over a
  ring with indeterminates, and a multiplicity-free
  \cite{VershikOkunkov96} semisimple
  specialisation \cite{ClineParshallScott99} (split, and we will only 
  consider characteristic 0 here).

For a tower as in I1, let $\BB_{\Alg_-}$  
denote the \Bratteli\ diagram of the semisimple case. Regarded as
a set,  $\BB_{\Alg_-}$   will here mean the vertex set of this graph. 

\item[I2]
(i) A quasi-hereditary global limit via an idempotent 
$\Ide \in \Alg_{m}$ (some $m$) and isomorphisms \cite{Green80} 
\eql(id sub)
\Ide \Alg_{n+m} \Ide \cong \Alg_n
\eq
(and hence a tower of {\em recollement} \cite{CPS88B}).
(Exclude consideration of specialisations in which $e$ is not well  defined.) 
\newline 
(ii) A map $\Pglob_{\Alg}$ from $\BB_{\Alg_-}$ to a global 
index set  $\weights$, which localises at each $n$ to an index set 
$\weights(n)$ for standard modules $\Delta(\mu)$
  \footnote{Standard in the quasi--hereditary sense, but we might also
    hope these modules are `nice' in some Kazhdan--Lusztig sense
    \cite{KazhdanLusztig79,GarsiaMclarnan88}.}  
  of $\Alg_{n}$ 
(i.e., such that 
$\weights(n)\hookrightarrow \weights(n+m)$ 
via the full embedding of $\Alg_{n}$--mod in $\Alg_{n+m}$--mod 
consequent on eqn.(\ref{id sub}), 
while $\dot{\cup}_{n} \weights(n) \cong \BB_{\Alg_-}$). 

Let $\WRes(\mu) \subset \weights$ denote the set of weights of standard
factors of the restriction $\Res({\Alg_{n+1}},{\Alg_{n}})\Delta(\mu)$, 
and $\WInd(\mu)$ of the corresponding induction. 
\item[I3] A space $V$ (for definiteness we
  will assume this is a real Euclidean space) and map $\weights
  \hookrightarrow V$ with the following properties. 
  The convex hull
  of $\WRes(\mu)$ intersects $\weights$ in a subset of $\{\mu\} \cup \WInd(\mu)
  \cup \WRes(\mu)$ ({\em locality} of induction and restriction); 
  the set $\rho_{\mu}$ of reflections in $V$ which fix $\{\mu\}$ and
  $\WRes(\mu)$ fixes $\weights$; 
  the group $W$ generated by
  $\cup_{\mu}\rho_{\mu}$ is an affine Weyl group
  \cite[\S4.2]{Humphreys90}; and $\weights$ is a subset of the set of point
  facets in the alcove geometry induced by $W$ on $V$. 
\item[I4] Control of bases for the algebras and standard modules ---
  including the means in principle to compute Gram matrices in the
  case of indeterminate parameters.
\item[I5] Explicit forms for the simplest primitive idempotents (in
  particular any primitive {\em and} central idempotents).
\item[I6] A linkage principle \cite{Jantzen87} --- $\Delta(\mu), \;
  \Delta(\nu)$ are in different blocks if there does not exist any $w
  \in W_l$ ($W_l$ a suitably rescaled version of $W$, depending on the
  specialisation) such that $w\mu=\nu$.
\end{enumerate}
\subsection{Role and realization of ingredients 1, 2(ii) and 3}
\label{affine hecke}\label{role of 123}

Recall \cite{Lusztig89C} (and cf. \cite{Lambropoulou94})
the affine Hecke algebra $\aH{n}$ defined by generators 
$\{1,X,g_1,\cdots,g_{n-1} \}$
and relations 
\eql(e01)
g_i g_{i \pm1} g_i=g_{i \pm1} g_ig_{i \pm1} \hspace{.8in} 
g_i g_j=g_j g_i  \hspace{.2in} \mbox{ $i \neq j \pm 1 $ }
\eq
\eql(e02)
g_1 X g_1 X = X g_1 X g_1  \hspace{1.06in} 
g_j X=X g_j  \hspace{.2in} \mbox{ $j>1$ } . 
\eq
\eql(local)
(g_i-q)(g_i+q^{-1})=0    
\eq
The {\em cyclotomic} Hecke algebra 
\cite{BroueMalle93} $\AH=\AH(n,d)$ is the quotient
$\Psi_d$ of $\aH{n}$ by 
\eql(e03)
\prod_{i=1}^{d} (X- \lambda_i) \; = 0.
\eq
Here $q,\lambda_1,\cdots,\lambda_d,\cdots$ are parameters, which we may begin
by regarding as indeterminates.
Write $\Aring$ for $\Z[q,q^{-1},\lambda_1,\ldots,\lambda_d]$ and
$\Aquot$ for the quotient field. Write $\AH^{\Aring}(n,d)$ for $\AH(n,d)$ over
$\Aring$. This is a free module over $\Aring$ \cite{ArikiKoike94} (see
\S\ref{prelims}). 
The case $d=2$,
$\lambda_1=-\lambda_2^{-1}$ is essentially the $B$--type Hecke algebra
(cf. \cite{Hoefsmit74,Martin89}).
For any $d' > d$ let $\Psi_d : \AH(n,d') \rightarrow \AH(n,d)$ also denote
the quotient by equation(\ref{e03}). 
Denote by $\AH(-,d)$ the sequence of inclusions
\eql(e03x) \AH(n,d) \subset \AH(n+1,d) . \eq


Usually we will fix an $\Aring$--algebra $k$ which is a field, 
and which {\em as a field} is $\C$, and consider 
$\AH(n,d)=\AH^{\Aring}(n,d) \otimes_{\Aring} k$.
The semi-simple generic structure of $\AH(n,d)$ over $\C$ is well
known, through that of the specialisation to the group algebra of the group 
$Z_d \wr S_n$ 
(confer \cite{JamesKerber81,Hoefsmit74,ClineParshallScott99} as in
\cite{ArikiKoike94}).
We recall it briefly.  
An {\em integer partition} $\mu$ of degree $n$ is a
list $(\mu_1,\mu_2, \cdots)$ of non-negative integers such that $\mu_i
\geq \mu_{i+1}$ and $\sum_i \mu_i =n$. 
There is a natural
correspondence with {\em Young diagrams} of degree $n$.  
Denote by 
$\AKweights_n=\AKweights_n^d$ the set of
ordered lists of $d$ integer partitions, of summed degree $n$ (called
$d$-partitions of degree $n$).  
For example
$$\AKweights^2_2 = \{((2),0),((1^2),0),((1),(1)),(0,(2)),(0,(1^2)) \} . $$ 
The conjugacy classes of $Z_d \wr S_n$ are readily seen to be 
indexed by $\AKweights_n^d$ \cite{MacDonald79,MartinWoodcockLevy2000}, 
thus $\AH(n,d)$ has 
\newcommand{\DDAK}[1]{\Delta_{#1}}
simple modules $\DDAK{\mu}$ indexed by $\mu \in \AKweights_n^d$. 
Similarly the Bratteli diagram $\BB=\BB_{\AH(-,d)}$ 
of the natural tower of semi-simple algebras
$\AH(-,d)$ is
determined by
\eql(restrictionrule01)
{\Res({},{})}_{\AH(n,d)}^{\AH(n+1,d)} \DDAK{\mu}
= \bigoplus_{i=1}^d \bigoplus_{j} \DDAK{\mu-e^i_j }
\eq
where the sum over $j$ is over possible subtractions of one box from the $i^{th}$
Young diagram of $\mu$.

For each $n>1$, there are $2d$ one dimensional irreducible 
representations $R_{\pm l }$ ($l=1,2,\cdots,d$) of $\AH(n,d)$, given by 
$R_{\pm l }(X)=\lambda_l$ and $R_{\pm l }(g_i)=\pm q^{\pm 1}$. 
The representation $R_{ l }$ corresponds to the 
module $\DDAK{ \mu}$ with multipartition $\mu=(,,\mu_l,,)$
in which 
all component integer partitions are empty except the $l^{th}$ partition, 
which is either $(n)$ (case $l>0$) or $(1^n)$ (case $l<0$).  
For each $n$ 
we may associate an unique primitive (and central) idempotent to each of these 
representations, in the algebra over generic $k$ \cite{Hamermesh62}. 
We write these idempotents $\EEE(\pm,l,n)$. 
For $l \in \{ 1,2, ...,d \}$ 
the element $\EEE(\pm,l,n)$ of $\AH(n,d)$ uniquely obeys 
\eql(ideqs)
(g_i \mp q^{\pm 1}) \EEE(\pm,l,n)=0 \hspace{.1in} (i=1,2,\cdots,n-1),
\hspace{.35in} (X - \lambda_l ) \EEE(\pm,l,n)= 0, 
\hspace{.35in} (\EEE(\pm,l,n) )^2=  \EEE(\pm,l,n) .
\eq
The inclusion (\ref{e03x}) allows us to regard $\EEE(\pm,l,n)$ as an
idempotent in $\AH(n+1,d)$, albeit neither primitive nor central in
general. Indeed the idempotent will be expressible as a sum of
primitive idempotents in reciprocity with the rule
(\ref{restrictionrule01}).%
\footnote{Let $1=\sum_{\mu} e^{\mu}$ be the unique decomposition of
  $1$ into primitive {\em central} idempotents of $\AH(n+1,d)$. Then
  $\EEE(\pm,l,n) =\sum_{\mu} \EEE(\pm,l,n) e^{\mu}$ 
  is this decomposition; i.e. it is also unique, even though the
  decomposition of $1$ into primitive idempotents is not.} 


In \S\ref{s gen d} of this paper 
we give explicit formulae for all $\EEE(\pm,l,n)$ for all $d$. 
For now
we will be concerned particularly with
$\EE(\pm l) := \EEE(\pm,l,2)$. 
The reason for this is the desire for a small but significant generalisation 
of the set of dominant weights and the weight spaces underlying 
Soergel's procedure for case $A_{m-1}$ (i.e. $U_q sl_{m}$). 
Although the induction and restriction rules are straightforward, 
and satisfy ingredient 1, the `weight space' of $\AH(-,d)$ 
(in which distance $d(\mu,\lambda)$ is the minimum number of steps 
on $\BB$ from $\mu$ to $\lambda$) 
is somewhat unmanageable geometrically, cf. our desired ingredients~2 and~3. 
What is wanted is {\em something like} an analogue 
for $\AH(-,d)$ of the quotients 
\eql(ringel)
H^m_n \cong \mbox{End}_{U_q sl_{m}}(V_m^{\otimes n})
\eq
of the ordinary Hecke algebra $H_n$. 
The \Bratteli\ diagram $\BB_{H_{-}}$ 
of the ordinary Hecke algebra is the Young tree, but  
via (\ref{ringel}) the quotients $H^m_n$ are the natural incarnations of the 
Ringel duals of $A_{m-1}$ quantum groups, 
and hence may be associated to the same weight spaces, and satisfy
ingredients 1--6. Let us briefly review this. 
Let $e^{\pm}_{m+1}$ denote the two primitive and central idempotents
of $H_{m+1}$
(for simplicity assume $[m+1]!\neq 0$ for now). 
The $H^m_n$ are such that 
\eql(q001)
0 \longrightarrow H_n e^-_{m+1} H_n  
       \longrightarrow H_n \longrightarrow H^m_n \longrightarrow 0
\eq 
is exact. For $[m]!\neq 0$ they 
are quasi--hereditary and satisfy ingredient 2, for example, through 
\eql(eHe=H)
e^-_{m} H^m_{n+m} e^-_{m} \cong H^m_{n} . 
\eq
Let $\weights^{1,m}$ denote the set of Young diagrams of $< m$ rows, 
regarded as a subset of $\Z^{m-1}$. 
Let $v=(1,1,\ldots,1)\in\Z^{m-1}$ and let $\Z^{m-1}/v$ denote the 
corresponding quotient set. 
Note that the injective map from 
$\Z^{m-2}$ into $\Z^{m-1}$ given by 
$(\mu_1,\mu_2,\ldots,\mu_{m-2}) \mapsto (\mu_1,\mu_2,\ldots,\mu_{m-2},0)$
has image a set of representative elements of 
$\Z^{m-1}$ in $\Z^{m-1}/v$. 
Denote by $\Pglob^{m-1}$ the corresponding surjective map from 
$\Z^{m-1}$ to $\Z^{m-2}$ (and also its restriction to 
$\weights^{1,m}$, whose image is $\weights^{1,m-1}$). 
Note that  $\BB_{H^m_-}$, the set of weights for $H^m_n$ for all $n$, 
is $\weights^{1,m+1}$. 
In the sense of ingredient 2(ii) 
the set of isomorphisms (\ref{eHe=H}) collapses $\BB_{H^m_-}$
into $\weights^{1,m}$, which is the set of dominant
weights of $sl_m$,
via $\Pglob^{m+1}$ \cite{MartinWoodcock98b}.
That is 
$\Pglob_{H^m} = \Pglob^{m}$; 
and $\weights^{1,m}(n)$, the index set for $H^m_n$, is the subset
of $\weights^{1,m}$ of diagrams of degree $\leq n$ and congruent to
$n$ modulo $m$. 

\newcommand{\HD}{H^{{\cal D}}}%
\newcommand{\DSI}{{\cal D}}%
\newcommand{\bigplus}{+}%

Delightfully, we find that only the $n=2$ idempotents are needed 
for an analogue of equation(\ref{q001}) for $\AH(n,d)$ 
(see also \cite{\MartinWoodcockLevy2000\pre}). 
Denote the sum of the ideals generated by 
$\{ \EEE(-,l,2) \, | \, \mbox{all $l$} \}$ by 
$$
\DSI_d =\DSI_d(n) := \bigplus_{l=1}^d \AH \EEE(-,l,2) \AH  
$$
(we will modify this definition very slightly later). 
Define algebra $\HD=\HD(n,d)$ by $\HD(1,d)=\AH(1,d)$ and for $n>1$ by
exactness of the sequence 
\eql(q002)
0 \longrightarrow   \DSI_d
       \longrightarrow \AH \longrightarrow \HD \longrightarrow 0   . 
\eq
The idea is to restrict consideration to the subset 
$\AKweights^{\DSI}$ of $\AH(-,d)$--weights in 
which each integer partition is the trivial partition of that degree. 
Such an $\AH(-,d)$--weight is characterised by a sequence of $d$ 
non-negative integers, i.e. the degrees of the component integer
partitions (in this way we have an action of $\Pglob^d$ on $\AKweights^{\DSI}$). 
This sequence need not be ordered as an integer partition, and hence the set of 
such weights maps onto the set of all weights of $A_{d-1}$
(i.e. not just the usual dominant weights). 
For example with $d=3$, the weight $((2),(4),(3))$ becomes $(2,4,3)$,
and $\Pglob^3((2,4,3))=(-1,1)$:
\[
\begin{array}{ccccccc}
( \input{./xfig/young2.latex}, \input{./xfig/young4.latex}, \input{./xfig/young3.latex} )
& \hspace{.1in}  \mapsto  \hspace{.2in} &
\input{./xfig/young243.latex} 
& \hspace{.2in}  \mapsto  \hspace{.2in} &
(-1,1) 
\hspace{.1in} \mapsto \hspace{.2in} 
\includegraphics{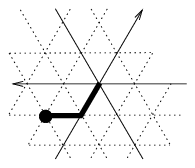}
\\ \in \weights^3_9 && && \in \mbox{$A_2$--weights}
\\
&&& \mapsto & \!\!\!\!\!\!\!\!-2\omega_1+\omega_2=(-2.(1,0)  \hspace{.4in} \input{./xfig/A2eg_alt.eepic}
\\ &&&      &   +1.(1,1)) 
\end{array}
\]
(the second row just shows the same weight in terms of fundamental
weights); 
while with $d=2$
\[
( \input{./xfig/young2.latex}, \input{./xfig/young4.latex} )
 \hspace{.1in}  \mapsto  \hspace{.2in} 
\input{./xfig/young24.latex} 
 \hspace{.2in}  \mapsto  \hspace{.1in} 
-2 \in \Z 
\]
this last, note, being the weight set used in \S\ref{decnos}. 
Note that $\DDAK{\mu} \in \AH\!-\!\!\mod$ is also in 
$\HD\!-\!\!\mod$ if and only if $\mu \in \AKweights^{\DSI}$ 
(in which case, as an $\HD$--module we write it $\DD (\Pglob^d(\mu))$). 
To see this note, from (\ref{restrictionrule01}), that the restriction of
$\DDAK{\mu}$ to $n=2$ contains a copy of one of the excluded one
dimensional modules (not necessarily a {\em direct} summand in
general) if and only if $\mu$ contains at least one 
integer partition with a second part (i.e. a Young diagram with
a second row).  


By construction then, the tower $\HD(-,d)$ has ingredients 1,
3 (restriction is local via (\ref{restrictionrule01}) and induction via Frobenius
reciprocity) and, at least formally, the final part of 2. 
Of course the construction is most interesting if it can be made to
include non--semisimple specialisations (else the {\em fundamental
  alcove} of \S\ref{decnos} is the whole weight space and the Soergel
procedure is trivial). 
In the remainder of this paper we address ingredients 2, 4, 5 and 6 
from this point of view. 
In particular we identify the $d=2$ case with the blob algebra. 
This has useful implications for all $d>2$. 
Note that, fixing $k$, there are a number of distinct ways to impose a
quotient relation of the form of $\Psi_d$ 
on $\AH(n,d')$,
corresponding to the choice
of factors to remove in the strengthening of the relation~(\ref{e03}). 
For $\lambda'=(\lambda_1,\lambda_2,...,\lambda_{d'})$ and $\lambda$ any subsequence
of this of length $d'-d$ let $\Psi_d^{\lambda}$ denote the
strengthening by omission of the factors $(X-\lambda_i)$ with
$\lambda_i$ in $\lambda$. 
The quotient $\Psi_d^{(\lambda_i)}$ commutes with
the quotient to $\HD$ so 

\mprop{
\label{hmmm}
Fixing $k$, 
there are $d+1$ ways ($\Psi_d^{(\lambda_i)}$, $i=1,2,..,d+1$) to quotient to pass
from $\HD(n,d+1)$ to $\HD(n,d)$ ($n>0$). 
}

Returning to the example of $d=3$, we see that as $n$ varies the image
of $\Psi_2^{(\lambda_i)}$ sweeps out a $2\pi/3$ radian arc of weight
space 
(there will be a picture illustrating this in section~\ref{s4.3}), with 
the union of these arcs over $i=1,2,3$ giving the complete space. 
If $k$ gives a non--semisimple specialisation then, 
as we will see, there is at least one $i$ such that 
the corresponding tower of blob (i.e. $d=2$) algebras is critical 
(in the sense of \cite{MartinSaleur94a,MartinWoodcock2000}, or
\S\ref{decnos}), i.e. it has one or more reflection points. 
As $n$ varies a given $d=2$ reflection point 
sweeps out a straight line in this arc in $d=3$ weight space
(see \S\ref{s4.3} for details). 
This is then a reflection
line of the $d=3$ alcove geometry, to which the
Soergel procedure may be applied.  

Our approach to ingredient 2 is through representation theory, and the
last part of the paper addresses this. It includes a `walk though'
review of some earlier work in statistical mechanics 
which explains our approach to this problem, and concludes by defining
certain representations of $b_n$ 
(and $\HD(n,d)$) which are candidates for `tensor
space' representations (i.e. they would establish ingredient 2 if
faithful, see \S\ref{Sreps} 
--- the question of faithfulness is {\em not} resolved here). 


\section{Preliminaries}\label{prelims}

Let involution $t:\Z [q,q^{-1}] \rightarrow \Z [q,q^{-1}]$ be given by
$t(q)=-q^{-1}$; 
and let $s$ act on $ \Z [\lambda_1,\cdots,\lambda_d ] $
by permuting the indices on the  $\{ \lambda_i \}$ cyclically.
Let $[n]=q^{n-1}+q^{n-3}+\cdots + q^{1-n}$.

For $A$ an algebra, let $Z(A)$ denote the centre of $A$; and for $B
\subseteq A$ let $Z_B(A)$ denote the centralizer of $B$ in $A$ (so
$Z_B(A) \supseteq Z(A)$ and 
$Z_A(A)=Z(A)$). 
 
\mdef
For $a=(a_1,a_2,\cdots,a_n)$ an $n$--tuple of natural numbers, 
the {\em Weyl orbit} of $a$ in $\N^n$ is the orbit
of the $S_n$ action permuting indices. 
Let  $Y=(Y_1,...,Y_n)$ be an $n$--tuple of variables in some
commutative ring $R$. 
Define 
\[ Y^{a} = \prod_{i=1}^n Y_i^{a_i}   \]
and the {\em monomial symmetric polynomial} (see
e.g. \cite{FultonHarris91}) 
\[ (Y^{a})^{\Sigma} = \sum_{a'} Y^{a'}  \]
where the sum is over all elements in the Weyl orbit of $a$. 

\newcommand{\Abasis}{{\bf A}}

\mdef
Let $R$ be a unique factorisation domain, $A$ 
a free $R$--module with basis $\Abasis$, 
and $K$ the field of fractions of $R$. 
Consider any $e \in A \otimes_{R} K$
and let $a_e \in A$ be such that 
${\bf e} = a_e e \in A \subset A \otimes_{R} K$. 
If there exists an $a_e$ such that 
the coefficient of some  ${\goth a} \in \Abasis$ in $\EEEu(,,)$  is 1,
then $a_e$ and $\EEEu(,,)$ are unique with this property up to a unit. 
If $A$ is an $R$--algebra and $e$ and idempotent,  
call such an $\EEEu(,,)$ a
{\em preidempotent}, 
and $a_e$ the corresponding {\em normalisation} of $e$.

\newcommand{\usually}{usually}
Now let $R$ above be $\Aring$, and 
consider $A$
as a collection of $\C$--algebras by specialisation.
We say that a property holds {\em generically} if it holds on an
(Zariski) open subset of parameter space, and that it holds {\em
  \usually} if 
the condition for failure to hold may be expressed as a single finite
polynomial in a single variable 
(e.g. $[2]$ is generically and \usually\ invertible, 
$(q^2\lam_1 -
\lam_2)$ is generically but not \usually\ invertible). 
 
Let $e=\frac{\EEEu(,,)}{a'_e}$ be {\em any} explicit expression 
for an element in $A$ over $\Aquot$ as above. 
For any point $x$ in parameter space
(i.e. any $k$) there is an open region with $x$ in its closure in
which the polynomial $a'_e$ has no root, so $e$ may be evaluated as a
limit at $k$. This process does not guarantee a unique finite limit. 
However, if $e$ is a primitive and central idempotent then 
any two finite limits must be the same, since they will have the
properties of a unique primitive and central idempotent in $A$ over
$k$. (I.e., they will induce the same simple projective module, 
$A \EEEu(,,) \otimes_{} k$.)


\mdef
Let $\sigma_i$ denote the elementary transposition $\sigma_i
= (i \; i+1) \in S_n$, so $\sigma_i(i)=i+1$ and so on 
\cite{Hamermesh62}.
Let $\Bbas_n$ be a maximal set of inequivalent reduced words 
in the generators $\{g_1,...,g_{n-1}\}$ of $H_n$. For each $w \in
\Bbas_n$ note that there is a natural (reduced expression for an) 
element of $S_n$ associated to it by substituting $g_i \leadsto \sigma_i$.



\mdef
Let  
\eql(e010)
X_1 := X \hspace{1in} X_j := g_{j-1}X_{j-1}g_{j-1}
\eq
so in $\aH{n}$  
\eql(e011)
 [X_j,X_k]=0 
\eq
\eql(e011')
[X_j,g_k]=0 \hspace{.2in} \mbox{ $j \neq k,k+1$ }
\eq
\eql(e012)
g_{k} X_{k+1}= X_kg_k +(q-q^{-1}) X_{k+1}
\hspace{1in}
g_{k} X_{k}= X_{k+1}g_k -(q-q^{-1}) X_{k+1}
\eq
and so 
\eql(e012')
[X_{k} +X_{k+1},g_k]=[ X_{k} X_{k+1} ,g_k]=0 . 
\eq
Now 
\[
X_{k}^{j} +X_{k+1}^{j} 
= ( X_{k}^{j-1} +X_{k+1}^{j-1} )( X_{k} +X_{k+1} )
  - (X_{k}^{j-2}   +   X_{k+1}^{j-2}) X_{k} X_{k+1}
\]
so
\[ [g_i , (X^{a})^{\Sigma} ] =0 \]
for all $i$ and any $a$. 


\newcommand{\gh}{\hat{g}}
\newcommand{\asp}{X^{\Sigma}}
\mdef
Let $\asp$ denote the algebra of symmetric polynomials in the
$X_i$s. Evidently $\asp \subseteq Z(\aH{n})$, and in fact 
Bernstein has noted that this is an equality (see appendix).
 


\mdef \label{AKbasis}
It follows 
from equations (\ref{e011} -- \ref{e012}) 
that any product of generators of $\AH(n,d)$ can be expressed as a 
$\Z[q,q^{-1},\lam_1,\ldots,\lam_d]$-linear 
combination of words from the set
\newcommand{\AKbas}{{\mathfrak C}}
\[
\AKbas_n^d = 
\{ X^{a} w \; | \; a \in \{0,1,\cdots,d-1 \}^n, \; w \in \Bbas_n
\}
\]
The dimension of this spanning set is clearly $d^nn!$, which is also the 
dimension of $\AH(n,d)$, thus 

\mprop{ \mbox{{\rm \cite{ArikiKoike94}}}
The set $\AKbas_n^d$ is a basis for \mbox{{\rm $\AH(n,d)$}}.}

Linearly extend $s$ and $t$ to act on $\AH(n,d)$, fixing $\AKbas_n^d$
pointwise. 


\mdef \label{testit}
\mprop{ 
\mbox{{\rm \cite{MartinWoodcockLevy2000\pre}}}
Let $B$ be a basis for $\AH(n-1,d)$, $D$ a basis for $<X>$ and 
$G=g_1g_2...g_{n-1}$. Then 
\[
\{ aGb, g_1aGb, g_2g_1aGb,\cdots,g_{n-1}\cdots g_2g_1aGb \; 
| \; a \in D, \; b \in B \}
\]
is a basis for \mbox{{\rm $\AH(n,d)$}}. 
}


\newcommand{\Epre}{z}
\newcommand{\polo}{{p}}
\mdef
Fix $d$ and set 
$\polo^l_n = \prod_{i\neq 1} (q^{2n-2} \lambda_l -\lambda_i)$. 
Note that  
\[
\Epre^l_n = \prod_{k=1}^n \left( \prod_{i\neq l} (X_k - \lambda_i) \right)
\]
lies in $Z(\AH(n))$ and obeys $\Psi_d((X-\lambda_l)\Epre^1_n) = 0$. 
Comparing with (\ref{ideqs}) we thus have 
\[
\EEE(\pm,l,n)=\frac{\Epre^l_n}{\prod_{k=1}^n \polo^l_k} e^{\pm}_{n} . 
\]
It follows that 
$\EEEu(\pm,l,n) = (\prod_{k=1}^n ([k] \polo^l_k)) \EEE(\pm,l,n)$ 
is a preidempotent for $ \EEE(\pm,l,n)$ (compare with the basis $\AKbas^d_n$).
Now define
$$
\DSI_d =\DSI_d(n) := \bigplus_{l=1}^d \AH \EEEu(-,l,2) \AH  
$$
(a modification of the definition in section~\ref{role of 123}) 
and $\HD$ accordingly. 
We will give another expression for $\EEE(\pm,l,n)$ shortly. 


\section{Standard modules and linkage}\label{linkage}

The generators and relations in 
equations(\ref{e01}) (and inverses) define the ordinary braid group $\braid{n}$.
Denote by $\bbraid{n}$ the extension by $g_0 = X$ (and inverse) obeying 
equations(\ref{e02}) (cf. \cite{Hoefsmit74,Martin91} and references therein). 
Thus $\aH{n}$ is a quotient of $\C\bbraid{n}$ by 
the quadratic relation in equation(\ref{local}). 

One realization of $\bbraid{n}$ is
as the group of braids on the cylinder, 
with $g_0$ becoming the pure braid in which the first string passes 
over
all the other strings and then around the cylinder. 
There is a natural `Young' embedding $\bbraid{n}\times\bbraid{m}
\hookrightarrow \bbraid{n+m}$. One places the second cylinder
concentrically inside the first, then allows the two cylinders to
converge in such a way that the nodes of $\bbraid{n}$, 
respectively $\bbraid{m}$,  remain
consecutive on the edges of the cylinder (while of course preserving
over/under information). 
There is a corresponding embedding $\aH{n}\times\aH{m} \hookrightarrow
\aH{n+m}$. The construction of `standard' modules (in the sense of
\cite{Cherednik91,Rogawski85}) follows from this. 
The quotient $\Psi_d$ complicates
this, in that the spectrum of $X_{n+1}$ 
(the image of $(1,X_1)$) is not that of
$X_1$ (the image of $(X_1,1)$). 
\newcommand{\lXr}{\langle X_i \rangle}%
\newcommand{\lXrFull}{\langle X_i \; | \; i=1,\ldots ,n \rangle}%
\newcommand{\cure}{{\epsilon}}%
\newcommand{\pip}{\pi}%
\newcommand{\ucure}{{\bf \cure^.}}
\newcommand{\ue}{{\bf e}}
\newcommand{\un}{\underline{n}}%
\newcommand{\nd}{d}%
\newcommand{\und}{\underline{\nd}}%
\newcommand{\stuff}{T}%
\newcommand{\rush}{\#}%
\newcommand{\Rest}{{\mbox{Res}}}%
\newcommand{\rG}{{\cal G}}%
\newcommand{\Ffset}{{\cal F}}%

\mdef \label{para3.1}
Let $\lXr$ denote the commutative subalgebra 
$\lXrFull \subseteq \AH(n,d)$.
Generically, as we will see, 
we may determine a {\em unique} basis of primitive (and of course
central) idempotents $\cure_{x}$ of $\lXr$ with $\sum_x \cure_x =
1$. There will be certain specialisations where this basis will not
make sense (certain idempotents will have preidempotents with
vanishing normalisation). 
In any case,
any primitive idempotent decomposition of 1 in $\AH(n,d)$ will be
different, but the unique primitive central idempotent decomposition 
of 1 in $\AH(n,d)$ will be expressible as a crudification of the above
(albeit depending on $k$). The generic case will be the least crude
(one idempotent per block/isomorphism class of simples); and it will be
necessary, formally, to combine certain of these generic idempotents (into
non--simple blocks) to make idempotents which make sense over $k$ in
non--semisimple cases.  
Each $\cure_x$ must obey 
$X_i \cure_{x} = x_i \cure_{x}$ with $x_i$ some scalar.
\newcommand{\Per}{\Pi}%
(Thus each induces a left $\AH$--module  $\Per_x := \AH \cure_{x}$. Since the
$H_n$ subalgebra of $\aH{n}$ maps isomorphically to its image in the
quotient we have $\AH \cure_{x} = H_n \cure_{x}$, of
rank $n!$.) 

Evidently $x_1 \in \{\lambda_i \}$, and with 
$\pip_j := \prod_{i \neq j} (X-\lambda_i)$ we have 
$X \pip_j = \lambda_j \pip_j$ and $\pip_j \cure_x \propto \cure_x$. 
For each such $\pip_j$ there exists a
minimal polynomial $\pip_{j,-}(X_2) = \prod_{k}(X_2-\lambda_{j,k})$ such that 
$\pip_j \pip_{j,-} =0$. 
Set $\pip_{j,k} = \prod_{l\neq k}(X_2-\lambda_{j,l})$ and 
$\ucure_{(j,k)}=\pip_j \pip_{j,k}$. 
Then $X_2 \ucure_{(j,k)} = \lambda_{j,k}  \ucure_{(j,k)}$. 
For each $\ucure_{(j,k)}$ there exists a polynomial
$\pip_{j,k,-}(X_3)$ such that 
$\pip_j \pip_{j,k} \pip_{j,k,-} =0$,
and so on. 
That is, the roots of such a polynomial are certain of the eigenvalues of the
$X_i$s. 


We can work out these {\em eigenvalues} of the
$X_i$s by looking at the properties of generically irreducible
representations as given in \S\ref{role of 123}. 
Let $\mu \in \AKweights_n^d$. 
A `standard' insertion of $\un = \{1,2,..,n \}$ 
into the boxes of $\mu$ is one such that 
deletion of the boxes containing $\{l,..,n\}$ produces a legitimate
$p$-partition for every $l$. 
Let $\stuff_{\mu}$ denote the set of all such standard insertions. 
For $i \in \un$ and $w \in \stuff_{\mu}$ 
there is a $k \in \{ 1,2,..,d \}$
such that $i$ appears in a box in the $k^{th}$ partition in $w$. 
Define $w(i) = k$.   
It will be evident from the restriction rule (\ref{restrictionrule01}) that
$\stuff_{\mu}$ may be used as a basis for $\DD_{\mu}$, once equipped
with a suitable action. 
 

\mdef \label{Xevals}
\mprop{
The action of the $X_i$s on this basis may be taken to be lower
triangular, with diagonal element of $X_i$ on $w \in \stuff_{\mu}$
given by $\lambda_{w(i)}q^{2w^i}$ where 
$w^i$ is the distance of $i$ off the main
diagonal in the $w(i)^{th}$ partition in $\mu$ 
(with distances {\em below} the diagonal
being negative).
}
\newline{\em Proof:} In case $n=1$ the claim holds since
$R_{i}(X)=\lambda_i$. Suppose the claim is true at level $n-1$. Then
the eigenvalues for $X_1,..,X_{n-1} $ at level $n$ are given by
restriction using the rule (\ref{restrictionrule01}). 
The eigenvalues for $X_n$ may be 
determined using the centrality of $\sum_i X_i$ etc.. \Qed


\mdef
Let $v_i$ denote the $i^{th}$ elementary vector in $\Z^{\nd}$. 
Describe a walk on $\Z^{\nd}$ (or $\Z^{\nd}/(1,1,..,1)$) 
by word $w=w_1w_2...$ in $\und$ such
that the vector between the $i^{th}$ and $i+1^{st}$ site visited is
$v_{w_i}$. 
For given word $w$ define 
$\rush_l(i)=\rush^w_l(i)=\sum_{j=1}^{l} \delta_{w_j,i}$. 
A reflection hyperplane $(i,j;x)$ is characterised by a pair 
$v_i,v_j$ ($i \neq
j \in \und$) not parallel to it, and the {\em signed} distance $x$ 
in the direction of $v_i$ of 
this hyperplane from 0 (i.e. the $x$ such that $0+xv_{i}$ lies on
it). A walk touches this hyperplane at $l$ if
\eql(touch)
\rush_l(i)-\rush_l(j)=x .
\eq 
Let $w$ be a walk which touches hyperplane
$(i,j;x)$ at $l$. 
The walk $w'$ got from $w$ by applying permutation $(ij)$
to every $w_t$ $t>l$ is called the {\em (affine) reflection} of  $w$ in
$(i,j;x)$ at $l$ 
(NB, the touching point $l$ is not in general uniquely defined by $w$
and $(i,j;x)$; or indeed existent). Every point of the reflection
after $l$ is the reflection of this point in $(i,j;x)$ in the usual
alcove geometry sense (see \cite[Ch.7]{Martin91}). 

If $w'$ meets another 
(not necessarily distinct)
hyperplane $(k,m;y)$ at $l'>l$ then of course
$w$ meets the image of  $(k,m;y)$ in $(i,j;x)$ at the same moment. 
If $w''$ is the reflection of $w'$ in  $(k,m;y)$ at $l'$ then we say
$w,w',w''$ in the same {\em walk orbit} of 
the reflection group generated by these hyperplanes 
(the $i^{th}$ points of these walks are in the same orbit
in the usual sense, for each $i$). We may think of folding up the
space along the 
set of hyperplanes in the group --- the walk orbits are the sets of walks which
are mapped into each other by this folding. 
For $G$ a reflection group generated by hyperplanes 
we write $w \sim^G w'$ if walks $w,w'$ in the same walk orbit of $G$. 
Each hyperplane partitions
space into  two parts (not counting the hyperplane itself). The
`outside' of the hyperplane  $(i,j;x)$, $x\neq 0$, is the part not
containing 0. 

\begin{figure}
\includegraphics{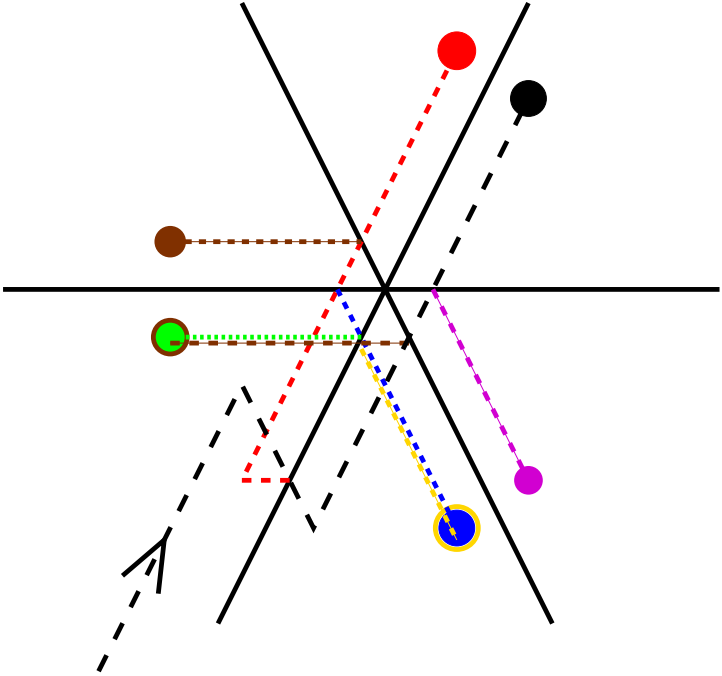}
 \caption{\label{walkies} Eight walks in a walk orbit for $d=3$
   (origin in the bottom left hand corner, endpoints circled).}
\end{figure}
\begin{figure}
\includegraphics{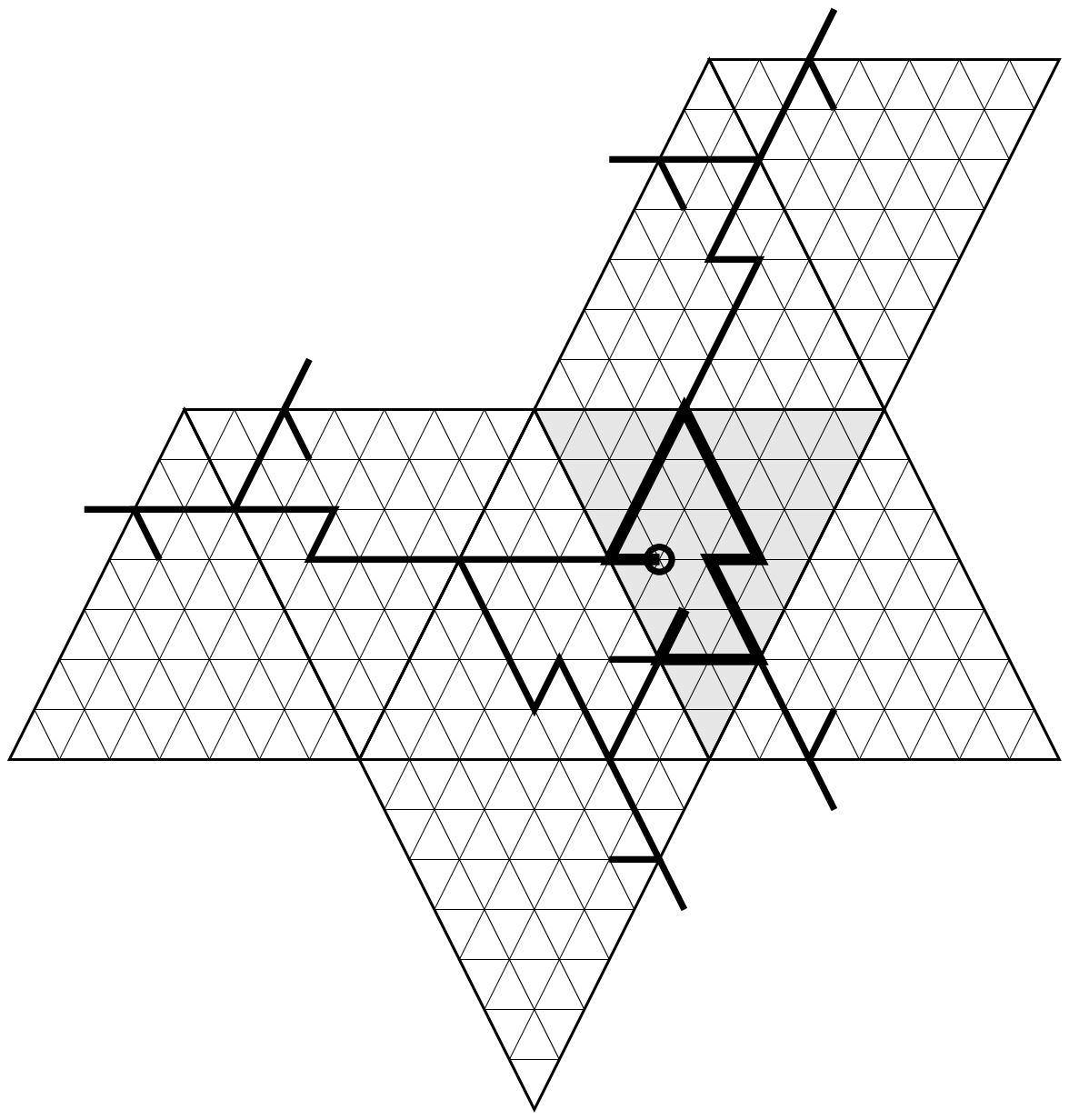}
 \caption{\label{A2p7tX} A walk,
   starting at 0 (circled), in the fundamental alcove (shaded) of
   a certain $d=3$ affine reflection group; and all walks in its walk
   orbit. Note that the two endpoints closest
   to 0 are each reached by two different walks.}
\end{figure}

We will restrict attention to hyperplanes {\em not} touching 0. 
Note that for each walk $w$ from 0 which finishes at a point $\mu$ outside some
hyperplane there is not in general a unique walk 
in its orbit which finishes at the image point of $\mu$, inside it. 
The orbit structure of walks is more complicated than that of points. 
For example, the endpoints of elements of the same walk orbit are
necessarily in
the same point orbit; but the converse does not follow. Further, 
a walk which stays in the interior of the fundamental alcove is in a
singleton orbit; and more generally a walk which touches hyperplanes a
total of $t$ times has $2^t$ elements in its orbit 
(see figure~\ref{walkies} and figure~\ref{A2p7tX}). 

For $G$ a reflection group generated by a set of reflection
hyperplanes $S$ let $\overline{S} \subset G$ denote the set of simple
reflections. 
It will be convenient to confuse these reflections with the
corresponding hyperplanes, noting that $\overline{S} \supset S$ in
general. In particular we may now compose hyperplanes by conjugation
in $G$, i.e. 
$(i,j;x) \circ (j,k;y)=(i,k;x+y)$  and so on. 
Let $\Aring' \subset \Aring$ consist of those elements of the form 
$(\lambda_i - q^{-2x} \lambda_j)$, all $i,j,x$. 
Define the {\em factor set} of $G$, $\Ffset (G) \subset \Aring'$ such that 
$(i,j;x) \in \overline{S} $ iff 
$(\lambda_i - q^{-2x} \lambda_j) \in \Ffset(G)$. 
The elements of $\Ffset(G)$ are called {\em factors}.  
Note that there is an equality of ideals 
\eql(ideal=)
\sum_{(i,j;x)\in S} \Aring (\lambda_i - q^{-2x} \lambda_j) = 
\sum_{(i,j;x)\in \overline{S}} \Aring (\lambda_i - q^{-2x} \lambda_j) 
\eq
since the composition of hyperplanes above corresponds to the addition
(up to a unit) of factors.


For a walk $w$ of length $n$ let $\mu(w)=(\rush_n(1),\rush_n(2),..)$, 
an {\em ordered} partition of $n$. 
To each walk  $w$ of length $n$ associate an element $\lambda^w \in
\Aring^n$ as follows:
\[ (\lambda^w)_i = \lambda_{w_i} q^{2(\rush_i(w_i) -1)}  . \]
For example $\lambda^{3331312} = (\lambda_3 , \lambda_3 q^2 ,
\lambda_3 q^4 , \lambda_1 , \lambda_3 q^6 , \lambda_1 q^2 , \lambda_2 )$. 


NB, This $\lambda^-$ gives an injective map from 
$\cup_{\mu \in \AKweights^{\DSI}(n)} \stuff_{\mu}$ into
  $\Aring^n$. 

\mdef \label{hyperp2k}
\mprop{ (i) Let $w'$ be the reflection of $w$ in $(i,j;x)$ at $l$, then
  every non--zero element of $\lambda^w -\lambda^{w'}$ is 
of the form
  $\pm q^{2\alpha}(\lambda_i - q^{-2x} \lambda_j)$ (some $\alpha\in\Z$). 
(NB, no mention of $l$ in the implication.) 

(ii) Let $w,w'$ be two walks. 
 If every element of $\lambda^w -\lambda^{w'}$
is divisible by
  $(\lambda_i - q^{-2x} \lambda_j)$ then $w,w'$ are related by
  one (or more) reflection in $(i,j;x)$ 
for some $l$ (resp. $l<l'<l''...$). 

(iii) Let $G$ be the group generated by a set of simple reflections 
$\{ (i,j;x), (k,l;y),.. \}$. 
Walks $w,w'$ are in the same walk orbit of $G$ 
(i.e. related by a sequence of reflections in hyperplanes in $G$) 
if and only if each non--zero $(\lambda^w -\lambda^{w'})_i $ 
is divisible by an 
element of $\Ffset (G)$. 
}
\newline{\em Proof:} (i) Consider the first point at which $w$ and $w'$
differ, which we may take to be $w'_{l+1} = (ij)w_{l+1} =j$. 
We have 
$(\lambda^w - \lambda^{w'})_{l+1} 
= \lambda_{i} q^{2(\rush^{w}_{l+1}(i) -1)} - 
  \lambda_{j} q^{2(\rush^{w'}_{l+1}(j) -1)}$.  
Now $\rush^w_{l+1}(i)=\rush^w_{l}(i)+1$
and $\rush^{w'}_{l+1}(j)=\rush^{w'}_{l}(j)+1$, so, 
cf. equation(\ref{touch}), the exponents differ by $2x$ as required. 
Since subsequent points $k$ with $w_k \neq w'_k$ are {\em all} those at which 
$\{w_k,w'_k \} \cap \{i,j\} \neq \emptyset$ 
the {\em difference} in exponents is preserved. 


(ii) is a special case of (iii). \ignore{{ 
In the first case we have 
$(\lambda^w - \lambda^{w'})_{m} 
= \lambda_{w_m} q^{2(\rush^{w}_{m}(w_m) -1)} - 
  \lambda_{w'_m} q^{2(\rush^{w'}_{m}(w'_m) -1)}$. 
This is not divisible by $(\lambda_i - q^{-2x} \lambda_j)$. 
The third case is similar (but fails on the $x$ rather 
than the $i$ or $j$). 
In the second case we have 
$(\lambda^w - \lambda^{w'})_{m} 
= \lambda_{w_m} q^{2(\rush^{w}_{m}(i) -1)} - 
  \lambda_{w_m} q^{2(\rush^{w'}_{m}(i) -1)}$. 
This can only be divisible by $(\lambda_i - q^{-2x} \lambda_j)$ if it
is zero, but this would require 
$\rush^{w}_{m}(i) = \rush^{w'}_{m}(i)$. 
}} 


(iii)
(Only if part) Note that, by construction, 
$(\lambda^w - \lambda^{w'})_i$ is, up to a unit, an element of $\Aring'$. 
Each pair $w\sim^G w'$ may be related by a series of simple reflections:
$w \sim^G w'' \sim^G .. \sim^G w'$, and each intermediate 
$(\lambda^{w^i} - \lambda^{w^{i+1}})_i$ is of the required form by (i). 
Evidently 
$(\lambda^w - \lambda^{w'})_i = 
\sum_i (\lambda^{w^i} - \lambda^{w^{i+1}})_i$. 
Now confer equation(\ref{ideal=}). 

(if part) Let $w|_{m}$ denote the first $m$ steps in walk $w$. 
Fix $w$ and consider $w'$ such that $w \not\sim^G w'$. 
Suppose that point $m$ is the first at which $w|_m \not\sim^G w'|_m$. 
Let $\tilde{w}|_{m-1} = w'|_{m-1}$ and $\tilde{w}_m$ be
such that $w|_m \sim^G \tilde{w}|_m$. 
Thus $\tilde{w}|_m$ is a reflection of $w'|_m$ in some hyperplane $h$ at
$m-1$, with $h \not\in G$. 
Thus $\lambda^{w'}-\lambda^{\tilde{w}}
=(0,0,..,0,(q^{2\alpha}(\lambda_i - q^{-2x} \lambda_j)))$ 
for some $\alpha,i,j,x$, where $i \neq j$ since $w'$ and $\tilde{w}$
agree at $m-1$. The non--zero factor is not a factor of $\Ffset(G)$
since $h \not\in G$. 
Now consider 
$\lambda^{w}-\lambda^{{w'}}=
(\lambda^{w}-\lambda^{\tilde{w}}) + (\lambda^{\tilde{w}}-\lambda^{w'})$. 
The $m^{th}$ element of the first summand on the right 
{\em is} divisible by some $f \in \Ffset(G)$ by the {\em only if} part, 
while that of the second summand is not, by the above argument. Thus 
the LHS is not. 
 \Qed


For example, consider the sequences 333, 331 and 321 
with $\lambda^{333}=(\lambda_3,q^2\lambda_3,q^4\lambda_3)$,
$\lambda^{331}=(\lambda_3,q^2\lambda_3,\lambda_1)$ and 
$\lambda^{321}=(\lambda_3,\lambda_2,\lambda_1)$. 
We have $(3,1;2) 333 = 331$ and $(3,2;1) 331 = 321$ and 
$\lambda^{333} - \lambda^{321} = 
(0,q^2(\lambda_3 - q^{-2} \lambda_2),q^4(\lambda_3 - q^{-4} \lambda_1))$.


\mdef \label{para3.5}
\mprop{ For each $w \in \cup_{\mu \in \AKweights^{\DSI}(n)} \stuff_{\mu}$ 
there is an element $\ucure_w$ of $\lXr$ as described above obeying 
\[ X_i \ucure_w = (\lambda^w)_i \ucure_w . \] }
\newline{\em Proof:} For a $d$-partition $\mu$ of form
$((\mu_1),(\mu_2),..)$ (i.e. $\mu \in \AKweights^{\DSI}$)
and $w \in \stuff_{\mu}$ identified with the
corresponding word $w$ we have $w_i = w(i)$ and $\rush_i(w_i) = w^i $. 
Now cf. proposition~\ref{Xevals}. \Qed


Let $R_x$ denote the representation of $\lXr$ corresponding to
any $\ucure_x$ as above. 
It follows that the $\lXr$--simple character of
$\Rest^{\AH}_{\lXr}(\DD_{\nu})$ is 
\eql(char formula01) 
[ \Rest^{\AH}_{\lXr}(\DD_{\nu}) : R_w ] = \delta_{\mu(w),\nu} . \eq
NB, this is not a {\em unique} characterisation unless:
\newline {\bf Cor.} {\em If $k$ is generic, i.e. $\lambda^-$ remains
  injective on passing from $\Aring^n$ to $(\Aring\otimes k)^n$, then
  $R_w \cong R_{w'}$ implies $w=w'$ and $\lXr^{\DSI} \subset \HD$ has
  enough simples to be semisimple.}


\mdef \label{ref gp}
We define $\rG(k)$, the {\em reflection group induced by $k$} as follows. 
For each triple $(i,j;x)$ such that 
$(\lambda_i -q^{2x}\lambda_j)\otimes k = 0$
include the hyperplane  $(i,j;x)$ as a  generator. 
Note that this is sensible in as much as 
$(i,j;x)\circ (j,k;y)=(i,k;x+y)$ while 
$\lambda_i -q^{2x}\lambda_j =0$ and
$\lambda_j -q^{2y}\lambda_k =0$ imply  
$\lambda_i -q^{2(x+y)}\lambda_k =0$. 
If $d'$ of the parameters $\{ \lambda_i \}$ are related in this way then the
group generated is $S_{d'}$, unless $q$ is a root of unity, in which
case it is the affine extension. We will usually refer to the group as affine
Weyl group regardless, to emphasize the fact that even in the finite
case the hyperplanes do not pass through 0. 
That is, we exclude from consideration any $k$ in which $\lambda_i
=\lambda_j$ does not imply $i=j$.

\mprop{ Over $k$, $R_w \cong R_{w'}$ iff $w,w'$ are in the 
same walk orbit of the affine Weyl
group induced by $k$ (i.e. related by a some series of reflections). 
}
\newline {\em Proof:} There is an isomorphism iff $\lambda^{w}-\lambda^{w'}$
vanishes over $k$. Every term is of the form 
$q^{\alpha}(\lambda_i -q^{2x}\lambda_j )$. 
Now apply proposition~(\ref{hyperp2k}). 


\mdef
\mprop{ [Linkage] If there exists a nontrivial homomorphism $\DD_{\mu}
  \rightarrow \DD_{\nu}$ over $k$ then $\mu,\nu$ lie in the same orbit 
of the affine Weyl group induced by $k$. 
}
\newline {\em Proof:} 
\[ \xymatrix{\ar@{}[dr]|{\Downarrow}
\DD_{\mu}   \ar@{->}[d]^{\Rest} \ar@{->}[r]  &  \DD_{\nu}  \ar@{->}[d]^{\Rest} \\
\sum_{w \in \stuff_{\mu}} R_w \ar@{->}[r]
                                         &  \sum_{w \in \stuff_{\nu}} R_w
}
\]
Now apply ((\ref{ref gp}) Proposition). 
\Qed

Note that this strengthens immediately to exclude interorbit maps from
any submodule of 
$\DD_{\mu}$ to any quotient of $\DD_{\nu}$, 
i.e. $\mu \not\sim^{\rG(k)} \nu$ implies $\DD_{\mu},\DD_{\nu}$ have no
composition factors in common. 
Thus, {\em under the assumption that every simple occurs in} 
$\{ \mbox{Head}(\DD_{\mu}) \, | \, \mu\in\AKweights^{\DSI} \}$, 
as for a quasi--hereditary algebra, 
we have linkage in the form of I6. 
 

\mdef
\mprop{
Let $w,w' \in \cup_{\mu \in \AKweights^{\DSI}(n)} \stuff_{\mu}$ 
be reflections of each other 
in $(i,j;x)$ at any $l$. Then over any $k$ in which 
$(\lambda_i -q^{2x} \lambda_j)$ vanishes we have an isomorphism of
left {\rm $\AH$}--modules
{\rm 
\[ \AH \cure_w \cong \AH \cure_{w'} . \]
}
}
\newline {\em Proof:} By ((\ref{ref gp}) Proposition) $\ucure_{w}$ and
$\ucure_{w'}$ induce isomorphic $\lXr$--modules. \Qed 


\mdef
Given the injectivity of $\lambda^-$ it will be convenient to be able
to refer to a walk $w$ either directly or via its image $\lambda^w$. 
This unifies the labelling schemes for $\cure_w$ in 
(\ref{para3.1}) and (\ref{para3.5}). 

Note that $R_x( \pip_{j} ) =0$ for all $R_x$ except those with
$x_1=\lambda_j$. 
Let $S_j$ be the set of possible values of $R_{x}(X_2)$
when $x_1=\lambda_j$
(i.e. $S_j = \{ \lambda_k \neq \lambda_j, q^2 \lambda_j,  q^{-2}
\lambda_j \}$). Then 
$R_{(\lambda_j,x_2,..)}( \prod_{s \in S_j} (X_2 - s)) =0$
for any such $x$, and $\pip_j \prod_{s \in S_j} (X_2 - s)$ lies in the
radical of $\lXr$, and hence some power of it vanishes. 
This tells us, up to multiplicity, the roots of $\pip_{j,-}$, and
hence of the preidempotent $\pip_{j}\pip_{j,k}$ 
(obtained by omitting a factor $(X_2-\lambda_{j,k})$, $\lambda_{j,k} \in S_j$,
from $\pip_{j}\pip_{j,-}$).
If $\lambda^-$ is injective then the radical is $\{ 0 \}$ and all the
roots can be distinguished. Iterating this argument we have

\mprop{ \label{eexplicit}
For $w$ a walk
\[ \ucure_w = \prod_{i=1}^{|w|} 
\prod_{\begin{array}{l} x\in \cup_{\nu, |\nu|=i}\stuff_{\nu} \\ 
                        \mbox{such that} \\
                        x|_{i-1} = w|_{i-1} \\
                        x_i \neq w_i  \end{array} }
(X_i - R_x(X_i))  
\] 
where the union is over {\em all} multipartitions satisfying the constraints.}

For example, 
$\ucure_{111} = 
 \ucure_{11} (X_3 - q^{-2} \lambda_1)\prod_{i=2}^d(X_3 - \lambda_i)$; 
and 
$$\ucure_{11122}=
\ucure_{111} \left( (X_4 - q^{-2} \lambda_1)(X_4 -q^6\lambda_1)
\prod_{2<i\leq d}(X_4 - \lambda_i) \right) \hspace{.75in}$$ 
$$\hspace{1.2975in} . \left(
(X_5 - q^{-2} \lambda_1)(X_5 - q^6\lambda_1)(X_5 - q^{-2} \lambda_2)
\prod_{2<i\leq d}(X_5 - \lambda_i) \right).$$ 

Note that unless the {\em omitted} factor
in $\pip_{j}$ or $\pip_{j,k}$ coincides (under $X_1 \leftrightarrow
X_2$) with a factor in the other, then the product $\pip_{j}\pip_{j,k}$ is
automatically in $Z(\AH(2,d))$ (simply multiply in all the factors
apparently required
for symmetry, then replace these with scalars using the eigenvector
property --- the excluded cases are those where one or more such scalar
vanishes).


For example, in case $n=2$, $d=2$ 
$\langle X_i \rangle$ has rank 6. 
We have 
$(X-\lambda_2)(X_2 - q^{-2}\lambda_1)(X_2 - q^2\lambda_1)(X_2 - \lambda_2)=0$ 
so, for example, 
$\ucure_{11}=\ucure_{\lambda_1, q^2\lambda_1}
=(X-\lambda_2)(X_2-q^{-2}\lambda_1)(X_2-\lambda_2)$ 
is a preidempotent with 
$X_2\ucure_{\lambda_1,q^2\lambda_1}
     =q^2\lambda_1\ucure_{\lambda_1,q^2\lambda_1}$;
while 
$\ucure_{\lambda_1\lambda_2}
     =(X-\lambda_2)(X_2 -q^{-2}\lambda_1)(X_2-q^2\lambda_1)$ 
is a preidempotent with 
$X_i\ucure_{\lambda_1\lambda_2}=\lambda_i\ucure_{\lambda_1\lambda_2}$. 
The preidempotent 
$\cure_{\lambda_1,q^2\lambda_1}$ lies in $Z(H)$ (provided $q^2 \neq 1$)
since it can be symmetrized:  
$(X_1-q^{-2}\lambda_1)\ucure_{\lambda_1,q^2\lambda_1}=
(\lambda_1-q^{-2}\lambda_1)\ucure_{\lambda_1,q^2\lambda_1}$.
On the other hand, $\ucure_{\lambda_1\lambda_2}$ 
cannot be rescaled to its symmetrized form 
because the symmetrizing factor
$X_2-\lambda_2$ kills it.

For $\AH(2,d)$ we have 
$\prod_{j \neq 1}(X-\lambda_j) 
 (X_2 -q^2\lambda_1)(X_2 -q^{-2}\lambda_1)
 \prod_{j \neq 1}(X_2 -\lambda_j) =0$, so similar considerations apply.
Indeed
they do 
for all $n$. 
In particular $\ucure_{11..1} \in Z(\AH)$ unless $q^2=1$. 
More generally, we may proceed as follows.


\mdef \label{eecommute}
For $\mu$ an ordered partition of $n$ let $H_{\mu}$
denote the corresponding Young subalgebra of $H_n  \subset \AH(n,d)$. 
For $w$ a walk let $e_w \in H_n$ denote the $q$--Young symmetrizer 
\cite{Cohn82II} \cite[\S9.3]{GoodmanWallach98} associated to $H_{\mu(w)}$. 
A walk $w$ is said to be {\em sorted} if it takes the form
111..22.. (more precisely, if $w_i \leq w_{i+1}$ for all consecutive
pairs of steps in $w$). 
There is a unique sorted walk in each $\stuff_{\mu}$, denoted
$w(\mu)$. 
A walk $w$ is said to be {\em direct} if all steps in a given
direction are taken consecutively (thus a sorted walk is direct). 

\mprop{
For $w$ sorted, and $(1+q^2)$ and each $\lambda_i$ invertible, 
$\ucure_w$ commutes with $H_{\mu(w)}$. 
}
\newline {\em Proof:} 
Let $\mu^t$ denote the $t^{th}$ interval of $\{1,2,..,n \}$ in the
ordered paritition $\mu(w)$, i.e. the $t^{th}$ set of integers fixed
under the action of the Young subgroup $S_{\mu(w)}$ of $S_n$ on $\{1,2,..,n \}$. 
We require to show, for each $t$, 
that the factors in $\ucure_w$ involving 
$X_{\mu^t} := \{X_i \, | \, i \in \mu^t \}$ 
constitute  a symmetric polynomial in these
variables, and hence commute with the  $t^{th}$ factor algebra in 
$H_{\mu(w)}$ (the remaining factor algebras commute with these
variables by equation(\ref{e011'})). 
Our walk is of the form $w=111..22..tt..$, and the factors in question are
(by ((\ref{eexplicit}) proposition)) those written out explicitly in:
\[
\ucure_{w}=\ucure_{111..22..} \prod_{i=t+1}^d(X_a -\lambda_i)
\left( 
\prod_{i=1}^{t-1}(X_a -q^{-2}\lambda_i)(X_a - q^{2\mu_i(w)}\lambda_i)\right)
 \hspace{.5in} \] \[ \hspace{.5in} 
\prod_{b\in {\mu^t}\setminus\{a\}} \left( (X_b - q^{-2}\lambda_t)
 \prod_{i=t+1}^d(X_b -\lambda_i)
\left( \prod_{i=1}^{t-1}(X_b -q^{-2}\lambda_i)(X_b - q^{2\mu_i(w)}\lambda_i)\right)
\right) 
\ldots
\]
(where $X_a$ is the first $X_i$ in $X_{\mu^t}$). It will be apparent
that this is rendered symmetric by multiplying by 
$ (X_a - q^{-2}\lambda_t)  $, but 
$ (X_a - q^{-2}\lambda_t) \ucure_w 
=  (\lambda_t - q^{-2}\lambda_t) \ucure_w $ so it is 
{\em already} symmetric provided that $(1-q^{-2})\lambda_t$ is invertible. \Qed

A similar property holds for direct walks. 

\mdef
\mprop{
If $w$ takes the form 111..22.. then $\cure_w e_w = e_w \cure_w$ (and
similarly for the preidempotent forms $\ucure$ and $\ue$). 
The modules {\rm $\DD_{\mu(w)} = \AH \ue_w \ucure_w$ }
are the left standard modules of {\rm $\AH$} with these weights. 
}
\newline {\em Proof:}
Note that $e_w \in H_{\mu(w)}$
and apply ((\ref{eecommute}) proposition). 


\mdef \label{standard injection}
Write $\mu \geq \nu$ if every change of direction in $w(\mu)$
occurs at the same step as 
one in $w(\nu)$ (NB, $w(\nu)$ may have
changes at other points as well). 
Note that $\mu \geq \nu$ implies 
$e_{w(\mu)} e_{w(\nu)}=e_{w(\mu)}$, and hence 
$\ue_{w(\mu)} \ue_{w(\nu)}= \kappa_{\nu} \ue_{w(\mu)}$ for some 
$\kappa_{\nu} \in \Aring$. 

\mprop{ Suppose that $w(\mu) \sim^{\rG(k)} w(\nu)$, and $\mu \geq \nu$,
  and $\kappa_{\nu} \otimes k \neq 0$. 
Then $\DD_{\mu} \hookrightarrow \DD_{\nu}$. }
\newline
{\em Proof:} 
$$ \DD_{\mu}  = 
         \AH \ue_{w(\mu)} \ucure_{w(\mu)} 
      =   \AH  \ue_{w(\mu)} \ue_{w(\nu)} \ucure_{w(\mu)} 
\hookrightarrow   
         \AH \ue_{w(\nu)} \ucure_{w(\mu)}   
         \cong 
         \AH \ue_{w(\nu)} \ucure_{w(\nu)} = \DD_{\nu} . $$



\section{The case $d=2$ and the blob algebra}\label{s3}
\subsection{Idempotents in $\AH(n,2)$}

The primitive and central idempotents  $\EEE(\pm,l,n)$
corresponding to the four one--dimensional  
representations of $\AH(n,2)$ over $\Aquot$
may be constructed as follows. 

{\em Fixing} $d=2$ define 
\[
\PP(n)=\PPP(+2,n)= q^{n-1}[n]\left( q^{2n-2} \lambda_2 - \lambda_1 \right) 
\]
and $\PPP(-2,n)= t \PP(n)$, $\PPP(+1,n)= s \PP(n)$ and 
$\PPP(-1,n)= st \PP(n)$. 

\mdef
\mprop{ Set
\eql(e05)
\alpha_j = \frac{-\lambda_1}{\PP(j)} \hspace{.3in} 
\beta_j = \frac{q \PP(j-1)}{\PP(j)} \hspace{.3in} 
\gamma_j = \frac{q^{2j-2}}{\PP(j)}  . 
\eq
Then
 $\EEE(+,2,0)=1$ and 
\eql(e04)
\EEE(+,2,j+1)=
\EEE(+,2,j) \left( \alpha_{j+1} + \beta_{j+1} g_j +\gamma_{j+1} X_{j+1} \right) 
\EEE(+,2,j)
\eq 
and $\EEE(-,2,j)=t \EEE(+,2,j)$, 
$\EEE(+,1,j)=s \EEE(+,2,j)$,  and 
 $\EEE(-,1,j)=ts \EEE(+,2,j)$. } 


\noindent {\em Proof:} see \S\ref{gen d} 
(or simply consider $g_{j-1} \EEE(+,2,j)$). 

Examples: 
\[   \EEE(-,2,1) =  \EEE(+,2,1) 
   = \frac{ X_1 - \lambda_1}{\lambda_2-\lambda_1}    \]
\eql(e07)
\EEE(-,2,2)\;\; = \;\; \left( \frac{X-\lambda_1}{\lambda_2- \lambda_1} \right) 
\frac{-\lambda_1+q^{-2} g_1 X g_1  -\; q^{-1} (\lambda_2- \lambda_1) g_1}
{(1+q^{-2}) (q^{-2} \lambda_2 - \lambda_1 )}
 \left( \frac{X-\lambda_1}{\lambda_2- \lambda_1}  \right) . 
\eq

\newcommand{\ZSigma}{\Sigma}
\newcommand{\ZPi}{\Pi}
\newcommand{\Zlam}{\lam_1+\lam_2}
\newcommand{\Zpi}{\lam_1 \lam_2}
\newcommand{\Dbas}{{\mathfrak D}}

Define $$
\ZSigma = (X_1+X_2-(\Zlam)) 
\hspace{1in} 
\ZPi = (X_1X_2 -\Zpi) .
$$
NB, $\Dbas^2_2 = \{ 1, \ZSigma, \ZPi, \ZSigma g_1 ,  \ZPi g_1 \} $ 
is a basis of 
$Z(\AH(2,2))$
(see appendix, eqn.(\ref{basisII})).


In terms of $\AKbas^2_2$ and $\Dbas^2_2$ we have
\eql(e07')
\EEE(-,2,2) = 
\frac{ \lam_1^2 - \lam_1 (X+X_2) +XX_2 
      +q^{-1}\lam_1 (X+X_2-(\lam_1+\lam_2))g_1
      -q^{-1} (XX_2 - \lam_1 \lam_2)g_1}%
{(\lam_1 - \lam_2)(1+q^{-2})(q^{-2}\lam_2-\lam_1)} 
\eq
$$
= \frac{(-\lam_1 \ZSigma + \ZPi)(1-q^{-1} g_1)}%
{(\lam_2-\lam_1)(q^{-2}\lam_2-\lam_1)(1+q^{-2})}
$$
which form manifests the centrality of this idempotent. 
Note that the preidempotent
$$
\EEEu(-,2,2) = (-\lam_1 \ZSigma + \ZPi)(1-q^{-1} g_1)
$$
coincides with its $s$--image (up to a unit in $\Aring$) 
if and only if $\lam_1 = \lam_2$ in $k$. 
However, 
$$
\EEE(-,2,2) + \EEE(-,1,2) = \frac{q^{-2}\ZSigma +(1+q^{-2})\ZPi}%
{(q^{-2}\lam_2 - \lam_1)(q^{-2}\lam_1 - \lam_2)}
\frac{( q-g_1 )}{[2]} . 
$$

\newcommand{\rad}{{\mbox{rad}}}

A remark is in order on denominators and idempotent decompositions of 1. 
The idempotent decomposition of 
$1=\EEE(+,2,1)+\EEE(+,1,1) \in \AH(1,2)$ 
is not defined in case $\lam_1=\lam_2$ over $k$, and the radical 
$\rad(\AH(1,2))=k\EEEu(+,1,1)$. 
Obviously $R_{\pm 1}=R_{\pm 2}$ for any $n$ in this case. 
The decomposition of 
$1= ( \EEE(+,2,1)+\EEE(+,1,1) )
              +\EEE(,,{((1),(1))})
              + ( \EEE(-,2,1)+\EEE(-,1,1) ) \in \AH(2,2)$ thus has
the same limitation --- the bracketed sums do not split. 
Similarly when $q=-q^{-1}$ we have
$R_{+i}=R_{-i}$ (any $n$). 
More interesting is the case $\lam_2 = q^2 \lam_1$. 
Here both $\EEE(-,2,2)$ and $\EEE(+,1,2)$ are undefined. Clearly 
$R_{+1} \neq R_{-2}$ (unless $q^2=-1$, $\lambda_i=0$) so 
$\EEE(-,2,2) + \EEE(+,1,2)$ is also undefined, but 
$\EEE(+,2,2)$ and $\EEE(-,1,2)$ are well--defined so 
$\EEE(-,2,2) + \EEE(+,1,2) + \EEE(,,{((1),(1))})$ 
does not split and $R_{+1}$,$R_{-2}$ 
(i.e. $\Delta_{((2),)}$ and $\Delta_{(,(1^2))}$)
both must be composition factors of $\Delta_{((1),(1))}$. 
NB, at first site this seems problematic for our proposed $\HD$ weight
space, however 
over this $k$ we may identify  $R_{-2}$ with the simple head of 
$\Delta_{((1),(1))}$ and label it accordingly 
(rather than by its $\AH(n,d)$ label, which is $(,(1^2))$). 
This is a good paradigm for the subtleties with labels in realising
ingredient 2(ii) (cf. \cite{Mathas96}). 
In terms of the $A_1$ integral weight set $\Z \subset \Re$ 
(recall ${\Pglob^2}((,(2)))=-2$, $\Pglob^2 (((1),(1)))=0$, $\Pglob^2(((2),))=2$)
we depict the standard modules and 
morphism by:
\[
\input{./xfig/reflect1.eepic}
\]


\mdef
It will be convenient to note the equality 
\eql(useful01)
(X - \lambda_2) ( -\lambda_1 \Sigma + \Pi ) = 0  
\eq
and its $s$ image; and hence that 
\eql(useful02)
( -\lambda_1 \Sigma + \Pi )  \EEE(-,2,2)
= ( -\lambda_1 \Sigma + \Pi ) 
       \frac{g_1-q}{-q^{-1}-q}  
\eq 
and $s$ image. 


\mdef
Now consider the algebras $\HD(n,2)$  
obtained by quotienting by  $\EEEu(-,1,2)=0$, $\EEEu(-,2,2)=0$. 

First consider $n=2$. 
A basis of the ideal $\DSI_2(2)$ is simply 
$\{ \EEEu(-,1,2), \EEEu(-,2,2)   \}$ ($\lam_1 \neq \lam_2$). 
The basis $\AKbas_2^2$ of $\AH(2,2)$ is 
\[ \{ 1,X_1,X_2, X_1X_2,g_1,X_1g_1,X_2g_1,X_1X_2g_1 \} . \]
We can think of using the relation $\EEE(-,2,2)=0$ 
to eliminate $X_1X_2g_1$, and then using 
its $s$-partner to eliminate $g_1$, 
to obtain a basis for $\HD(2,2)$  
(see also \S\ref{conjectures}). 
Specifically, it is convenient to represent the quotient relations as
\eql(e013)
(X_1+X_2-(\Zlam))(g_1-q)=0 \hspace{1in} (X_1X_2 -\Zpi)(g_1-q)=0 . 
\eq
N.B. These two relations generate the {\em same} ideal. 


\newcommand{\U}{u}
Alternatively, 
rewriting $g_i-q =: \U_i$ 
(so $\U_i^2=-(q+q^{-1})\U_i$) and $X-\lambda_1 =: v$ we have 
\eql(e08)
v\U_1v\U_1= ( \lambda_1 q^{-1} - \lambda_2 q  ) v\U_1 \hspace{1.17in}
\U_1v\U_1v= ( \lambda_1 q^{-1} - \lambda_2 q  ) \U_1v
\eq
and similarly for the image under $s$. 
Combining (or by applying $(g_1-q)$ to (\ref{e013}i) from the left) we obtain 
\eql(e09)
\U_1v\U_1 =( \lambda_1 q^{-1} - \lambda_2 q  ) \U_1
\eq
(and $s$ image). 
\newcommand{\blobe}{e_-}
\newcommand{\ThL}{T--L}
\subsection{$\HD(n,2)$ and the blob algebra}\label{blob def}

\newcommand{\mm}{m}

The blob algebra 
$b_n$ \cite{MartinSaleur94a} may be defined by generators 
$\{ \blobe ,U_1,U_2,...,U_{{n-1}} \}$ 
and relations 
$\blobe\blobe=\blobe$, 
$U_i U_i =-[2]U_i$,
\eql(blob4) U_i U_{i\pm 1} U_i = U_i , \eq 
\eql(blob3) U_1 \blobe U_1 = y_e U_1 , \eq 
generators commute otherwise, and $q,y_e$ are parameters.%
\footnote{The blob algebra is usually defined in terms of a basis of decorated
  Temperley--Lieb diagrams (hence its name). That the presentation
  here is equivalent may be verified by a straightforward (if tedious)
  generalisation of the corresponding exercise for the ordinary
  Temperley--Lieb algebra as in \cite{Martin91}.}
(The ordinary Temperley--Lieb algebra $T_n(q)$ is the subalgebra
generated by the $U_i$s.) 
The properties I1--6 are established for the blob algebra in
\cite{MartinWoodcock2000}.  

There is a `good' parameterisation of $b_n$ by $q$ and $\mm$, where 
$y_e = -\frac{[\mm -1]}{[\mm]}$. In this  parameterisation it is
convenient to replace $\blobe$ by the rescaled generator 
$U_0 = -[\mm]\blobe$, or rather to replace $\blobe$ with $U_0$, with
relations  $U_0^2 = -[\mm]U_0$ and 
\eql(blob03) U_1U_0U_1=[\mm-1] U_1 . \eq 
The variant with these relations is isomorphic to the original
except (obviously) when $[\mm ]= 0$.


\mdef \label{->blob}
\mprop{
Let $\mm,q,\lambda_1,\lambda_2$ be such that 
\[
{[\mm-1]}{(\lambda_1-\lambda_2)}  = {(\lambda_1 q^{-1} -\lambda_2 q)}{[\mm]}
\]
(NB, $\lambda_1 \neq \lambda_2$ unless $q^2=1$ and 
$\mm \rightarrow\infty$, or $[\mm]=0$; 
otherwise, putting 
$
\frac{\lambda_1}{\lambda_2} = q^{-2r}
$
then $\mm=r$ --- we may take $\lam_1=\frac{q^{m}}{q-q^{-1}}$, then
$\lam_1-\lam_2 = [\mm]$). 
Then there is an isomorphism 
$$\phi :
        \HD(n,2) \rightarrow b_n$$ 
given by 
$\phi(\U_i)=U_i$ 
and 
$\phi(v)=U_0$. 
There is a homomorphism $\psi : b_n \rightarrow \HD(n,2)$ given by 
$\psi(U_i)=\U_i$, $\psi(U_0)=v$. 
}


\noindent
{\em Proof:} 
It is straightforward to verify the cyclotomic Hecke relations under
$\phi$. A direct calculation also shows that the image under $\phi$ of the 
{\em numerator} of $\EEE(-,2,2)$, as in equation~(\ref{e07}), 
vanishes identically if we use the form of $\lam_1$ given in the
proposition. 
The image of  $\EEE(-,1,2)$ vanishes similarly. 
Thus $\phi$ is
a surjective homomorphism (except possibly when $q^2=1$). 
A dimension count reveals that this surjective map is 
generically an isomorphism. 

For $\psi$, 
equation(\ref{e09}) verifies relation(\ref{blob03}). 
It is now enough to check $\U_1\U_2\U_1=\U_1$ 
(i.e., in $\AH(3,2)$, $\U_1\U_2\U_1-\U_1=\U_1(\U_2\U_1-1)
= [3]! e^-_{3} \in \DSI_2$).
Write
$f=\U_1\U_2\U_1-\U_1$ and consider $ f \DSI_2 f $ 
(if $f \in \DSI_2$ then $f \in f \DSI_2 f $, at least
if $[3]!$ is invertible, and this
is a much smaller and simpler object to work with). 
This (pre)idempotent subalgebra includes 
$$
f (X_1 + X_2 -(\lam_1 + \lam_2)) f
= f ((1+q^{-2})X -(\lam_1 + \lam_2)) f
$$ 
$$
f (X_1 X_2 - \lam_1 \lam_2 ) f 
= f (-q^{-1} X g_1 X  - \lam_1 \lam_2 ) f 
$$
and 
$$
f Xg_1g_2(X_1 + X_2 -(\lam_1 + \lam_2)) f
= f (Xg_1g_2X +Xg_1g_2g_1 X g_1 -Xg_1g_2(\lam_1 + \lam_2)) f
$$ $$
=f (-q^{-1}Xg_1X +Xg_2g_1g_2 X g_1 -q^{-2}(\lam_1 + \lam_2)X) f
=f (-(q^{-1}+q^{-3}) Xg_1X  -q^{-2}(\lam_1 + \lam_2)X) f
$$
A linear combination of these is 
$$
[3]! ((1+q^{-4})\lam_1 \lam_2 -q^{-2}(\lam_1^2 + \lam_2^2)) f
$$
so $f \in f\DSI_2 f$ at least in an open subset of parameter space. 
The coefficient may be rewritten 
$[3]!q^{-2}\lam_1\lam_2(q-q^{-1})^2 [\mm+1][\mm-1]$
so this covers most interesting cases. The remainder require a more
tedious calculation. 
\Qed


\noindent {\bf Corollary. }{\em 
The finite characteristic $b_n$ decomposition matrices in
  \cite{CoxGrahamMartin01} are a subset of type--$B$ Hecke
  decomposition matrices via the correspondence in \S\ref{role of 123}.
}


\mdef
Recall that $\mm= \pm 1$ implies the existence of a 
quotient algebra isomorphic to the 
Temperley--Lieb algebra (in addition to the noted 
\ThL\
subalgebra). The quotient  identifies $X=1$ (or $q^2$). 

\section{On general $d$}\label{s gen d}

\subsection{Idempotents in $\AH(n,d)$} \label{gen d}

\def\lamb(#1){\lambda^{( #1 )}}

The primitive and central idempotents corresponding to the $2d$ one-dimensional  
representations of $\AH(n,d)$ over $\Aquot$ may be constructed as follows. 
Fixing $d$, define 
\eql(Pn)
\PP(n)=\PPP(+1,n)= 
q^{n-1}[n]\prod_{i \neq 1}\left( q^{2n-2} \lambda_1 - \lambda_i \right) 
\eq
$\PPP(-1,n)= t \PP(n)$, $\PPP(+2,n)= s \PP(n)$,  
$\PPP(+3,n)= s^2 \PP(n)$ and so on. 
Define $\lamb(0)=1$ and 
\[
\lamb(i) = \sum_{d \geq j_1 > j_2 > \cdots > j_i >1} 
\left( \prod_{l=1}^{i} -\lambda_{j_l} \right) 
\]
(sum over descending positive integer sequences $(j_1,j_2,...,j_i)$,
with $j_1 \leq d$). 

\mdef \mprop{
Set
\eql(e05--d)
\beta_{j+1} = \frac{q \PP(j)}{\PP(j+1)} 
\eq
\eql(e05-d)
\alpha^0_j = \frac{\prod_{i \neq 1}(-\lambda_i)}{\PP(j)} 
\hspace{.3in}
\alpha^{d-i}_{j+1}= \frac{ q^{2j} \lamb(i-1) + 
\sum_{l=2}^{i} q^{(i-1)(2j-2)}(q^{2j} -1) \lambda_1^{l-1} 
 \lamb(i-l) }{ P_{j+1}}
 \hspace{.3in} (1 \leq i < d) . 
\eq
Then
$\EEE(+,1,0)=1$  
\eql(e04-d)
\EEE(+,1,j+1)=
\EEE(+,1,j) 
\left( \beta_{j+1} g_j + \sum_{i=0}^{d-1} \alpha^i_{j+1} X_{j+1}^i  \right) 
\EEE(+,1,j)
\eq 
$\EEE(-,1,j)= t \EEE(+,1,j)$, and so on.  
}

For example  $\beta_1=0$, $\alpha^i_1 = \lamb(d-1-i)$, and 
\[   \EEE(+,1,1) 
= \frac{ \prod_{i \neq 1} (X_1 - \lambda_i)}
{\PP(1)} .    
\]
(There are more examples in \S\ref{s4.3}.) 

{\em Proof:}


We work by induction on $j$, with the example above as base. 
Firstly note 
that $\EEE(+,1,j+1) \EEE(+,1,j)= \EEE(+,1,j)\EEE(+,1,j+1)  = \EEE(+,1,j+1)$
so, from proposition \ref{testit}, equation(\ref{e04-d}) gives a 
correct {\em form} for the $R_{+1}$ idempotent up to coefficients. 
Then note that in this form it is sufficient,
cf. equation~(\ref{ideqs}), to check the  identity 
$(g_j -q) \EEE(+,1,j+1)=0 $, and normalisation (i.e. idempotency). 
The former is a direct calculation, and the later may be checked by 
evaluating in $R_{+1}$, i.e. 
substituting $q$ for $g_i$ and $\lambda_1$ for $X$ (in which case 
$\EEE(+,1,j+1) $ must evaluate to 1). 
To begin we rewrite the claimed expression for 
$(g_j - q) \EEE(+,1,j+1) $ as 
\[
(g_j - q) \left(
\EEE(+,1,j) \left( \beta_{j+1} g_{j} \right) \EEE(+,1,j)
+ \left(  \sum_{i=0}^{d-1} \alpha^i_{j+1} X_{j+1}^i \right)
        \EEE(+,1,j) \right)
\]
using commutation properties. By the inductive assumption this is  
\begin{smalleq}
\[
=(g_j - q) \left(
\EEE(+,1,j-1) \left( \beta_j g_{j-1} +   \sum_{i=0}^{d-1} \alpha^i_{j} X_{j}^i 
\right) 
\left( \beta_{j+1} g_{j} \right) \EEE(+,1,j)
+ \left(  \sum_{i=0}^{d-1} \alpha^i_{j+1} X_{j+1}^i \right) \EEE(+,1,j) \right) 
\]
\[
=\EEE(+,1,j-1) \beta_{j+1}  
 \left( \beta_j (g_j - q) g_{j-1} g_{j} 
+   (g_j - q) \left( \sum_{i=0}^{d-1} \alpha^i_{j} X_{j}^i \right) g_{j} 
\right) 
\EEE(+,1,j)
+ (g_j - q) \left(  \sum_{i=0}^{d-1} \alpha^i_{j+1} X_{j+1}^i \right) \EEE(+,1,j)
\] 
\[
=\EEE(+,1,j-1) \beta_{j+1}  
 \left( 
   (g_j - q) \left( \sum_{i=0}^{d-1} \alpha^i_{j} X_{j}^i \right) g_{j} 
\right) 
\EEE(+,1,j)
+ (g_j - q) \left(  \sum_{i=0}^{d-1} \alpha^i_{j+1} X_{j+1}^i \right) \EEE(+,1,j)
\]
Now use the commutation rules on the first summand 
to bring the factor $g_j$ forward through the $X_{j}^i $s:
{\footnotesize 
\[ 
=  
(g_j - q)  
 \left( \left( g_{j} \beta_{j+1}\alpha^0_j +\alpha^0_{j+1} \right)+
    \sum_{i=1}^{d-1} \left( 
 \beta_{j+1} \left( -q \alpha^i_{j} -(q-q^{-1})
 \left( \sum_{l=i+1}^{d-1} \alpha^l_j \left( X_{j} \right)^{l-i} \right) \right)
+\alpha^i_{j+1} \right)
X_{j+1}^i  
\right) 
\EEE(+,1,j)
 \]}
but $X_j \EEE(+,1,j) = q^{2j-2} \lambda_1 \EEE(+,1,j) $ so we have 
\[
=  
(g_j - q)  
 \left( \left( -q^{-1} \beta_{j+1}\alpha^0_j +\alpha^0_{j+1} \right)
+
\right. \hspace{3.6in} 
\] \[ \left. \hspace{.8in} 
    \sum_{i=1}^{d-1} \left( 
 \beta_{j+1} \left( -q \alpha^i_{j} -(q-q^{-1})
 \left( \sum_{l=i+1}^{d-1} \alpha^l_j \left( q^{2j-2} \lambda_1  \right)^{l-i} 
\right) \right)
+\alpha^i_{j+1} \right)
X_{j+1}^i  
\right) 
\EEE(+,1,j)
\]
and equating coefficients to zero we get
\[
\alpha^{0}_{j+1} = q^{-1} \alpha^{0}_{j} \beta_{j+1}
\]
\[
\alpha^{d-1}_{j+1} = q \alpha^{d-1}_{j} \beta_{j+1}
\]
\[
\alpha^{d-2}_{j+1} = 
\left( q \alpha^{d-2}_{j} +
(q-q^{-1}) \alpha^{d-1}_{j} (q^{2j-2}\lambda_1)  \right)\beta_{j+1}
\]
and 
\[
\alpha^{i}_{j+1} = 
\left( q \alpha^{i}_{j} +
(q-q^{-1}) \sum_{l=i+1}^{d-1} \alpha^{l}_{j} (q^{2j-2}\lambda_1)^{l-i} 
 \right)\beta_{j+1}
\hspace{1in} (0<i<d).
\]
Without regard for normalization, any one coefficient may be 
chosen arbitrarily, so without loss of generality try  
$\beta_{j+1}=\frac{qP_j}{P_{j+1}}$.
Then $\alpha^{d-1}_{j+1}=\frac{\rho^{d-1}_{j+1}}{P_{j+1}}$ 
where $\rho^{d-1}_{j+1}=q^{2j}$,  
$\alpha^{0}_{j+1}=\frac{\lambda^{(d-1)}}{P_{j+1}}$, and 
\[
\alpha^{d-2}_{j+1}=
\frac{ \rho^{d-2}_{j+1} \lambda_1 
+ \rho^{d-1}_{j+1} \lamb(1) }{ P_{j+1}}
\]
for some $\rho^{d-2}_{j+1}$. 
Then 
\[
\rho^{d-2}_{j+1} 
= q^2\rho^{d-2}_{j} +(q^2-1) q^{2j-2} q^{2j-2} = q^{2j-2}(q^{2j}-1) . 
\]
Similarly
\[
\alpha^{d-3}_{j+1}= \frac{ \rho^{d-3}_{j+1} \lambda_1^2 
+ \rho^{d-2}_{j+1} \lamb(1)\lambda_1 +\rho^{d-1}_{j+1} \lamb(2) }{ P_{j+1}}
\]
and for $1 \leq i < d$
\[
\alpha^{d-i}_{j+1}= \frac{ \sum_{l=1}^{i}\rho^{d-l}_{j+1} \lambda_1^{l-1} 
 \lamb(i-l) }{ P_{j+1}}
\]
for some $\rho^{-}_{j+1}$. 
Then    
\[
\rho^{d-3}_{j+1} = q^2 \rho^{d-3}_{j} 
+(q^2-1) (q^{2j-2}  q^{2j-4}(q^{2j-2}-1)  + (q^{2j-2})^3 )   = q^{4j-4}(q^{2j}-1)
\]
\[
\rho^{d-4}_{j+1} = q^2 \rho^{d-4}_{j} 
+(q^2-1) ((q^{2j-2}) q^{4j-8}(q^{2j-2}-1)
+ (q^{2j-2})^2  q^{2j-4}(q^{2j-2}-1)  + (q^{2j-2})^4 )   
= q^{6j-6}(q^{2j}-1)
\]
and 
\[
\rho^{d-i}_{j+1} = q^{(i-1)(2j-2)}(q^{2j}-1) \hspace{1in} (1<i<d).
\]
Finally, the normalization condition, 
\[
\lamb(d-1) + \sum_{i=1}^{d-1} \left( \left( q^{2j} \lamb(i-1) 
+ \sum_{l=2}^{i} q^{(l-1)(2j-2)} (q^{2j} -1) \lambda_1^{l-1} \lamb(i-l) \right) 
\left( q^{2j} \lambda_1 \right)^{d-i} \right) 
=P_{j+1}-q^2P_j
\] 
is verified by direct computation. 
\end{smalleq}



By inspection of these idempotents we see that the algebras will not be 
generic in specialisations in which 
$\frac{\lambda_i}{\lambda_j}=q^{2r}$ for some $i,j,r \in \N$; 
and in certain specialisations in which $q$ is a root of unity 
(this also follows immediately from \cite{Ariki96,Mathas96}). 
Let us disregard, for the moment, the cases in which $X$ is 
degenerate or non-invertible. Then 
noting that rescaling all the $\lambda_i$s by the same factor produces 
an isomorphic algebra we can fix $\lambda_1=1$ (say) and adopt as parameters
$q,$ $\frac{\lambda_i}{\lambda_1}$ \cite{Mathas96}. 
It is illuminating to proceed by example 
(and cf. \cite[\S 5]{MartinWoodcockLevy2000}).


\subsection{The case $d=3$}\label{s4.3}

The primitive and central idempotents corresponding to the six one-dimensional  
representations of $\AH(n,3)$ may be constructed as follows. 
Fixing $d=3$, define 
$
\PP(n)=\PPP(+1,n)
$
and $\PPP(\pm i,n)$ as in equation~(\ref{Pn}). 
For example 
$\PP(1)=\PPP(\pm 1,1) = (\lambda_1 -\lambda_2)(\lambda_1-\lambda_3)$,
$\PP(2) = q[2](q^2\lambda_1 -\lambda_2)(q^2\lambda_1-\lambda_3)$.


As before
$
\beta_j = \frac{q \PP(j-1)}{\PP(j)} 
$,
and here
\eql(e05-)
\alpha^0_j = \frac{\prod_{i \neq 1}(-\lambda_i)}{\PP(j)} 
\hspace{.3in} 
\alpha^1_j = 
\frac{q^{2j-4}(q^{2j-2}-1)\lambda_1-q^{2j-2}(\lambda_2+\lambda_3)}{\PP(j)} 
\hspace{.3in} 
\alpha^2_j = \frac{q^{2j-2}}{\PP(j)}  
\eq
Compute $\EEE(\pm,l,j)$ by $\EEE(+,1,0)=1$ and then 
\eql(e04-)
\EEE(+,1,j+1)=
\EEE(+,1,j) 
\left( \beta_{j+1} g_j 
+ \alpha^0_{j+1} + \alpha^1_{j+1} X_{j+1} + \alpha^2_{j+1} X_{j+1}^2  \right) 
\EEE(+,1,j)
\eq 
and $\EEE(-,1,j)=t \EEE(+,1,j)$, and so on.  
For example 
\[   \EEE(+,1,1) 
= \frac{ X_1^2 - (\lambda_2+\lambda_3)X_1+\lambda_2\lambda_3}{\PP(1)}  
\]
and
\[
\EEE(+,1,2) = \EEE(+,1,1) \left( 
\frac{q\PP(1) g_1 
  + \lambda_2 \lambda_3 
  + ((q^2-1)\lambda_1 -q^2(\lambda_2+\lambda_3)) X_2
  + q^2 X_2^2
  }{\PP(2)} 
\right) \EEE(+,1,1)   
\]
\[
= q \frac{( X_1^2 - (\lambda_2+\lambda_3)X_1+\lambda_2\lambda_3 )}{\PP(1)}
 \frac{ ( X_2^2 - (\lambda_2+\lambda_3)X_2+\lambda_2\lambda_3 ) (g_1 +q^{-1})
  }
{\PP(2)}
\]





Given these results, let us consider the generalised Soergel procedure for
$d=3$ corresponding to that in \S\ref{decnos} for $d=2$. 
For illustration 
we consider a field $k$ in which the
equations $\frac{\lambda_i}{\lambda_1}=q^{-n_i}$ ($i=2,3$) 
and $\frac{\lambda_1}{\lambda_3}=q^{-n_1}$ are solved for
positive integer $n_i$ {\em only} in case $n_2=2$, $n_3=4$ 
(in particular, in this instance $q$ is not a root of unity). 
Figure~\ref{A2Soergel02}
illustrates the location of the corresponding reflection hyperplanes
(shown as thick lines)
in the $\HD$, i.e. $A_2$, weight space. 
The Kazhdan--Lusztig polynomials for this geometry are given in 
the figure at \cite[p.1289]{MartinWoodcockLevy2000}. 
The claim, then, is that if $\mu$ is a weight in the fundamental
alcove $A^0$ (such as 0) then $\Delta(\mu)$ has 
simple content (and Loewy structure)
\[
\begin{array}{ccc} & L_{\mu} & \\
L_{\mu.s} && L_{\mu.t} \\
L_{\mu.st} && L_{\mu.ts} \\
& L_{\mu.sts} \end{array} 
\]
(here $s$ and $t$ are the walls of $A^0$ 
--- we are abusing the notation $A.s$ of \S\ref{decnos} to apply to
weights in the obvious way). Of course, localising at small $n$, some
of these simple modules will vanish. Let us consider $n=0,1,2,3$. 

(For any $k$) 
we have $\HD(0,3) \cong k$, with basis $\{ 1 \}$. The weight for the
corresponding simple module is the innermost dot in the figures 
(marked $(,,)$ in figure~\ref{A2Soergel02}(a)).  
For $\HD(1,3)$, the idempotents
$\EEE(+,i,1)$ are all well defined and $1=\sum_{i=1}^3 \EEE(+,i,1)$ is an
idempotent decomposition into primitive idempotents (the corresponding
simples are marked as triangles). 

The set of weights for $\HD(2,3)$ are marked with squares (the weights
corresponding to the multipartitions 
$((2),,)$ and $(,,(2))$ have been explicitly labelled, to fix the
coordinate system). 
Note that $\EEE(+,i,2)$ is divergent in case $i=2,3$. 
The easiest way to see what is happening at $n=2$ is to recall the
$d=3$ version of
the proposition in \S~\ref{hmmm}:

\mprop{
Fixing $k$, there are three ways to quotient to pass
from $\HD(n,3)$ to $\HD(n,2)$ ($n>0$). 
}

The corresponding subsets of the set of
weights for $\HD(2,3)$ lie in three straight
lines --- the dashed lines shown in figure~\ref{A2Soergel02}(a). 
Since we have shown that the $d=2$ algebra is isomorphic to a blob
algebra in each of these cases we can give a complete description. 
The lines marked 12 and 23 correspond to singly critical blob
algebras, with reflection points the intersection points of these
lines with the thick lines shown. Accordingly there are injective 
standard module morphisms (`reflection' morphisms) as indicated by
arrows. 
The other blob line is a semisimple quotient 
(from the parameters it is nominally a singly critical case, 
but the reflection point lies at the outside edge of
localisation of weight space to $n=2$, so there is no image of any
$n=2$ weight in it). 
It is straightforward to show that there are no other non--trivial morphisms. 
Since the $d=3$ alcove structure is determined by the three $d=2$
quotients, and  
since the morphisms indicated are all also blob module morphisms,
their location is necessarily consistent with the $d=3$ Soergel procedure. 


For $\HD(3,3)$ the three blob quotients correspond to the three
dashed lines shown in figure~\ref{A2Soergel02}(b). 
The reflection morphisms shown within these lines are again simply blob
morphisms (NB, the reflection point on the third line now lies within the
localisation of weight space). 
The only question, then, concerns the morphisms into the module 
with $A_2$--weight 0. 
This weight is marked in the figure by the rank of the module, 
$\mbox{Rank}(\Delta(0)|_{n=3})=6$ (this weight coincides with 
the weight for the unique simple of $\HD(0,3)$ in our
scheme, since that module is a localisation of $\Delta(0)|_{n=3}$). 
A straightforward Frobenius reciprocity argument (using the morphisms
at level $n=2$ and the generic restriction rules) shows that there are
at least two homomorphisms into this module; but does not uniquely
determine the domain in either case. 
To confirm the indicated morphisms, 
 consider the walk basis of  $\Delta(0)|_{n=3}$ and the walk
orbits of these walks  --- see figure~\ref{A2sixer}. 
Note from the last of these pictures 
that $321 \sim^{\rG(k)} 333$. It follows that there
  is a homomorphism between the corresponding standard modules by
  ((\ref{standard injection}) proposition) (note that walk 321 is not
  sorted, but that the weaker condition of all steps in a given
  direction being taken consecutively is sufficient for this
  construction). 
The other claimed morphisms follow similarly. 

It is worth noting that $\Delta(0)|_{n=3}$
may generically be defined as 
$\Delta(0)|_{n=3} = \HD(3,3) e_{(1^3)}$, 
where $e_{(1^3)} = e^{-}_{3}$ is the usual $q$--antisymmetriser
\cite{MartinWoodcockLevy2000} 
(i.e., $\Delta(0)|_{n=3}$ is the globalisation of the unique
simple of $\HD(0,3)$). 
This is an illuminating construct to consider in any case. 
The module $\AH(3,3) e_{(1^3)}$ is, of course, much larger. 
By proposition~\ref{AKbasis} it has basis 
$\{ X^a e_{(1^3)} \, | \, a \in \{0,1,2 \}^3 \}$. 
To determine a basis in our case one must use the vanishing of
$\EEE(-,i,2)$ (which generates a significant part of $\AH e_{(1^3)}$). 
This problem is dealt with elsewhere \cite{CoxMartinRyom01}. 

\begin{figure}
\includegraphics{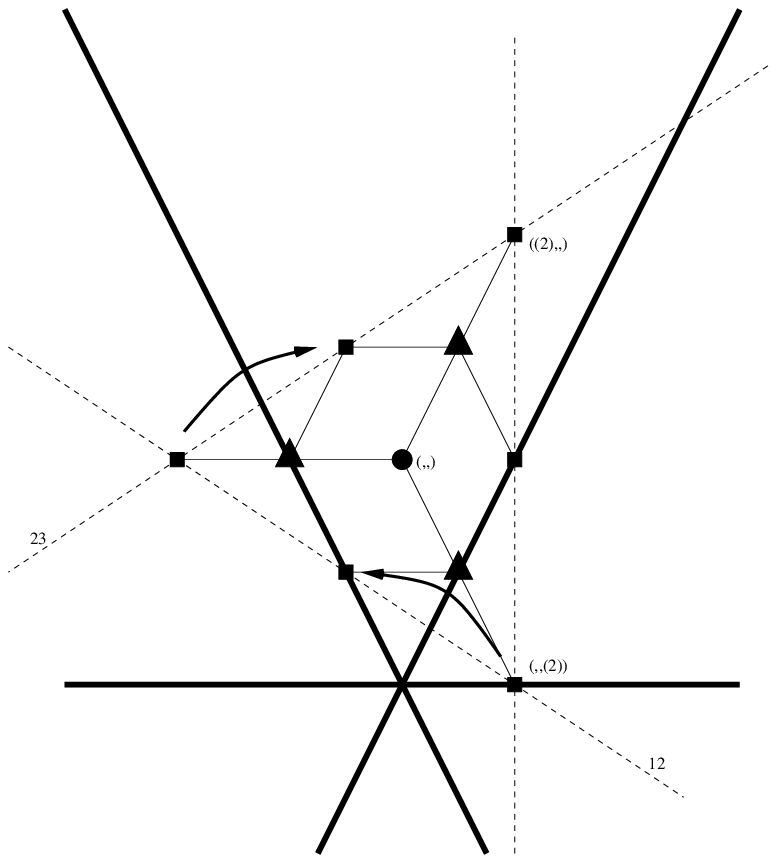} 

\includegraphics{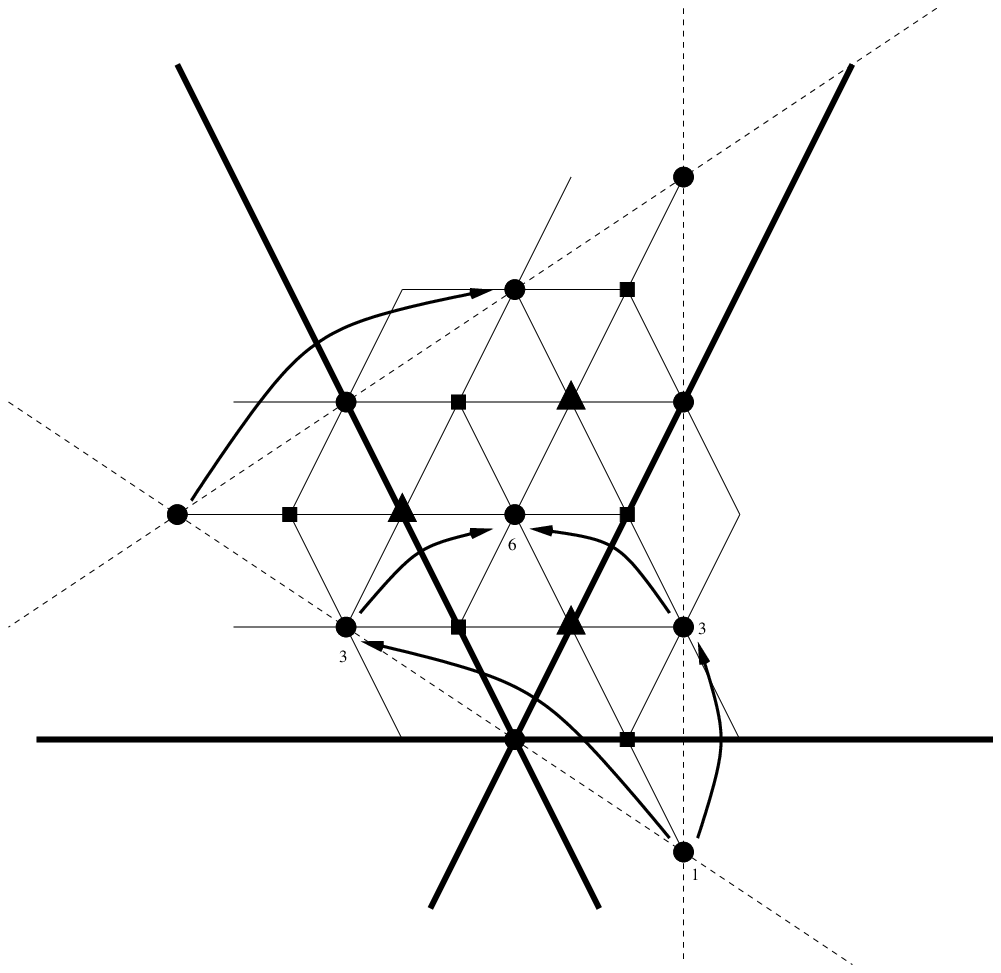} 
\caption{\label{A2Soergel02} Weights of $\HD(n,3)$ in weight space
  formalism 
  (a) for $n=0,1,2$; and  
  (b) for $n=0,1,2,3$.
  Weight set $\weights(0)$ consists of the empty weight (shown
  as a black circle marked (,,) in (a); $\weights(1)$ consists of the three
  weights marked with triangles; and so on (see main text).}
\end{figure}

\begin{figure}
\includegraphics{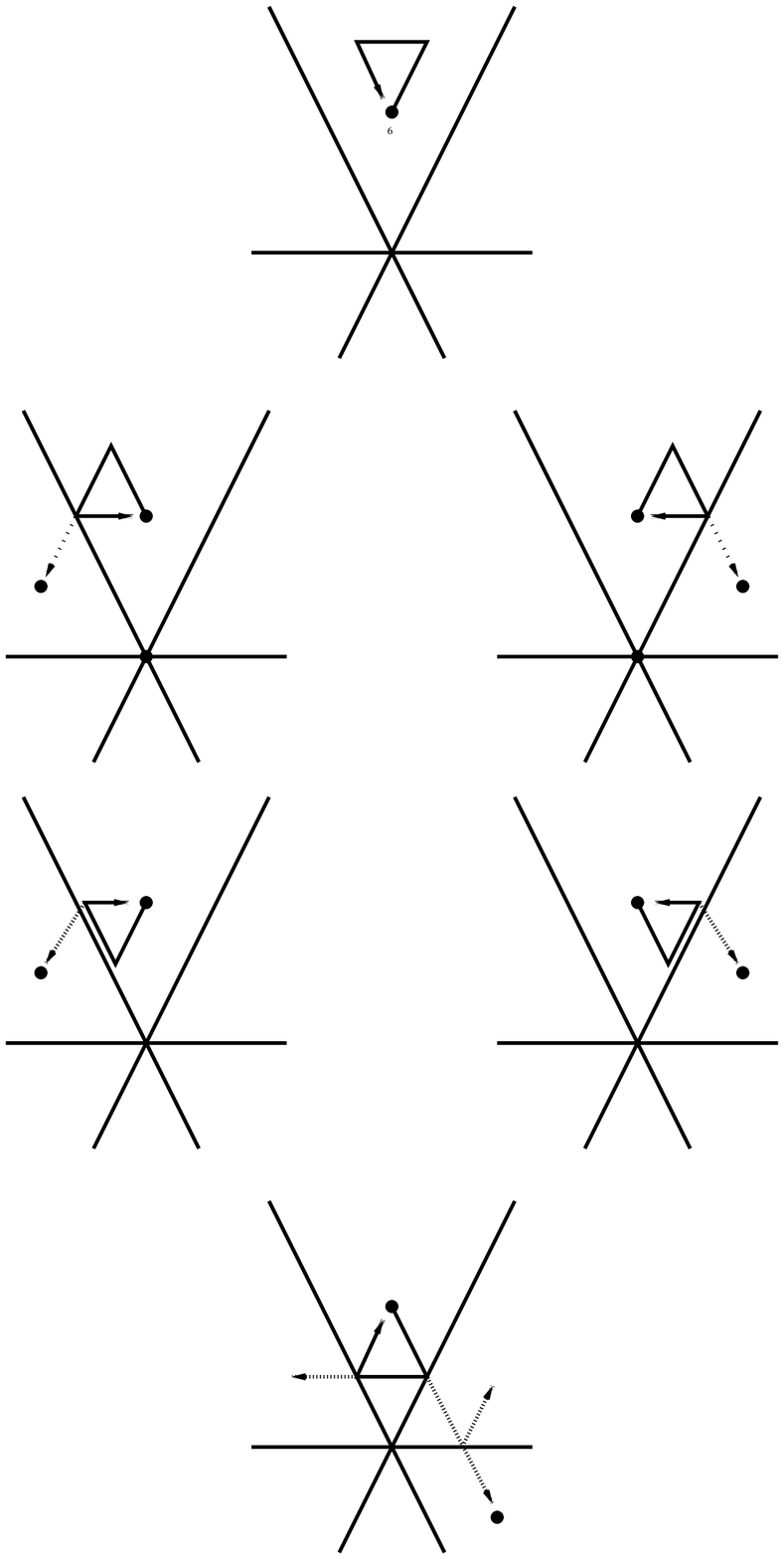} 
\caption{\label{A2sixer} The six walks in $T_0(3)$, and their walk orbits.}
\end{figure}



\subsection{Conjectured basis for
            $\HD(n,d)$}\label{conjectures}

Consider $w \in S_n$, a  permutation, and let $\goth S(w)$ be a 
corresponding reduced word in the Coxeter generators $\{ \sigma_i
\}$. 
The set map from $\goth S(S_n)$ to $\Bbas_n$ given by $\sigma_i
\rightarrow g_i$ is an isomorphism \cite{Humphreys90}. 
Thus from $w \in \Bbas_n$ 
we may read off a permutation (which we will also call $w$). 
Given such a permutation 
define symmetric relation $(w) \subset \underline{n} \times \underline{n}  $ 
by $i(w)j$ if $i>j$ implies $w(i)<w(j)$ (and symmetry). 
That is, $i(w)j$ if the lines $i$ and $j$ cross in the diagram of $w$. 

\mdef \label{claim1} {\bf Claim:}  
\[
\{ X^{(a_1,a_2,\cdots)}w \in \AKbas_n \; 
| \; i(w)j \mbox{ implies } a_i \neq a_j \}
\]
is a basis for $\HD(n,d)$.  

{\em Idea of proof:}
Note that the dimensions are right 
(consider the Robinson--Schensted correspondence 
in the form in \cite{StantonWhite86}, or, for example, \cite{Stanley97}).
Thus it is enough to show spanning or linear independence.
For the latter it is convenient to have a representation ${\cal R}$ of
$\HD$ (if linear independent in ${\cal R}$ then claim proven {\em and}
${\cal R}$ faithful). 

The remainder of the paper is concerned with representations of
$\HD$. 



\vspace{.1in}
\def\mm#1{{\cal M}^{#1}}
\newcommand{\mU}{{\cal  U}}
\newcommand{\TSR}{tensor space representation}%
\newcommand{\FTSR}{faithful \TSR}%
\newcommand{\qp}{{\bf q}}  
\newcommand{\qa}{q}        
\newcommand{\qb}{-q^{-1}}  
\newcommand{\rank}{N}      

\section{On representations of $b_n$ and $\HD(n,d)$}\label{s4}\label{Sreps}

The ordinary Temperley--Lieb algebra has a 
powerful diagram calculus (see \cite{Martin91} for a review). 
One motivation for the introduction of the 
blob algebra \cite{MartinSaleur94a} was to 
bring the utility of such a calculus to the representation theory of 
the periodic/affine algebras studied in \cite{MartinSaleur93}
(see also \cite[Ch.8]{Martin91}, 
\cite{PasquierSaleur90,Jones94b,GrahamLehrer98,GreenRM01} and references
therein). 
In fact all the finite irreducible representations of these 
infinite dimensional algebras may be
constructed using \cite{MartinSaleur94a} (although completeness is not
shown there, see \cite{GrahamLehrer01,MartinWoodcock2000}). 
From the point of view of lattice Statistical Mechanics, $b_n$ 
also renders the `seam' boundary conditions 
(as in \cite{BaxterKellandWu76}, and cf. \cite{Martin86a} for example) 
of the ice--type model \cite{Baxter82} 
into the algebraic formalism of Temperley and Lieb \cite{TemperleyLieb71}.
Indeed there are a number of mathematical and physical reasons 
(in addition to the pursuit of our ingredient I2 for $\HD(-,d)$) 
why a faithful `tensor space' representation of the blob would be useful 
(cf. \cite[\S4]{MartinSaleur94a}, 
\cite{PasquierSaleur90,BehrendPearce97,Donkin93}). 


By a {\em tensor space representation} of $b_n$ (or $\HD(n,d)$) 
we mean a representation {\em for each}
$n$ with underlying module of form $V(n)=V_{aux}\otimes V^{\otimes n}$
as a vector space, 
with $V_{aux}$ and $V$ finite dimensional vector spaces (cf.
\cite{KashiwaraMiwaStern95}); 
on which the generators act `locally' (in
particular $b_{n-1} \subset b_{n}$ acts trivially on the last factor
$V$, so the restriction of $V(n)$ is a manifest direct sum of dim~$V$
copies of $V(n-1)$); and which is well defined in arbitrary
specialisations.

Let $e \in \HD(n+d,d)$ be idempotent. A module $V(n+d)$ for a 
tensor space representation over $K$ 
is {\em globalisable by $e$} if $e$ projects the last $d$ tensor factors 
$V^{\otimes d} \rightarrow K$, 
and acts trivially on other tensor space factors. 


To establish  ingredient I2 {\em generically} for $\HD(-,d)$ 
(i.e., that there exists idempotent $e_d \in \HD(n+d,d)$ such that 
\eql(fembed) e_d \HD(n+d,d) e_d \cong  \HD(n,d)\eq 
cf. \cite{CPS88B,DlabRingel89}) 
certain generalised {\em braid diagrams} may be used 
\cite{\MartinWoodcockLevy2000\pre}. 
However, 

\mpr{
Suppose $V(n+d)$ is the module for a \FTSR\ of $\HD(n+d,d)$ over $\C$, 
and is globalisable by $e_d$. 
Then  equation~(\ref{fembed}) holds. 
}
\newline
To see this note that under these assumptions 
the actions of $e_d$ and $\HD(n,d)$ on $V(n+d)$
commute, and hence they commute in $\HD(n+d,d)$. Thus 
$\HD(n,d) e_d \subseteq e_d \HD(n+d,d) e_d \subseteq \HD(n+d,d)$
is a sequence of inclusions of algebras, and 
$\HD(n,d) e_d$ is isomorphic to $\HD(n,d)$. 
Thus the action of 
$e_d \HD(n+d,d) e_d$ on $V(n+d)$
would be isomorphic to the $\HD(n,d)$ action on $V(n)$. 
Since the latter has a
trivial kernel, this would establish equation~(\ref{fembed}) 
in general.

 
Any \FTSR\ in which $X$ acts 
non--trivially only in the first normal tensor factor $V$ 
(and an otherwise redundant factor $V^{aux}$) 
would be a likely candidate, because of the way ingredient~I2 
works in the $A_n$ case \cite{Jimbo85}. 
Let us briefly review this. 


\mdef \label{ordinaryTS}
The ordinary $H_n$ action (dual to that of 
$U_qsl_N$) on $V_{\rank}^{n}$ is
$$M^q_N:\braid{n} \rightarrow \End(V_N^{n})$$
$$g_i \mapsto 1\otimes 1\otimes\ldots\otimes 
                 \mm{q} \otimes 1\otimes 1\otimes\ldots\otimes 1$$
where $1$ denotes the unit matrix; and $\mm{q} $ (the $i^{th}$ factor) is given by 
$$\mm{q}|_{N=2}
=\mat{cccc} 
\qa^{}&0&0&0\\ 
0&\qa\qb&-1&0 \\
0&-1&0&0 \\
0&0&0&\qa^{} \tam$$
$$\mm{q}|_{N=3}
=\mat{cccccccccc} 
\qa&0        &0        &0 \\ 
0     &\qb+\qa&0        &-1&0 \\
0     &0        &\qb+\qa&0 &0      &0        &-1&0 \\
0     &-1       &0        &0 &0    &0 &0  &0 &0 \\
0     &0        &0        &0 &\qa &0   &0 &0 & 0\\
0     &0        &0        &0 &0      &\qb+\qa&0 &-1 &0\\
0     &0        &-1       &0 &0      &0        &0 &0  &0 \\
0     &0        &0        &0 &0      &-1       &0 &0  &0 \\
0     &0        &0        &0 &0      &0        &0 &0 &\qa 
\tam$$ 
and so on. Each such  
obeys a quadratic local relation with coefficients in $R=\Z[q,q^{-1}]$ 
(specifically $M_N((g_1-\qa)(g_1 \qb))=0$).  


Putting $U_i = g_i-\qa$ we have 
\eql(TL3)
M^q_2(U_{i}U_{i\pm 1}U_{i}-U_{i})=0
\eq
so $M^q_2$ factors through $\TL_n(q)$ 
(in fact it is faithful on $\TL_n(q)$, i.e. $\TL_n(q)=H^2_n(q)$). 

Recall from \S\ref{role of 123} that  $e^-_N$ denotes the
$H_N$ $q$--symmetriser (normalisable as an idempotent
provided that $[N]!$ is invertible in $K$ \cite{Martin91}). 
For $n\geq N$,
$V^n_N$ is
the tensor space module for
$H^N_n(q)$. 
It is easy to check that $V^n_N$ is globalisable by $e^-_N$.

\newcommand{\permm}{{\cal P}}

\mdef 
\label{charges}
Recall that the $\TL_n$ action on $V_2^n$ breaks up directly, over any field,
into summands $\permm_{\lambda}$ of fixed `charge' or weight $\lambda\in\N_0$, 
and then 
\eql(TLtilt)
\permm_{\lambda} = {\huge +}_{\mu \geq \lambda} \Delta_{\mu}
\eq
(generically a direct sum).


It is desirable
not only to have representations of $b_n$ that act on tensor space, but also that
they preserve some version of this charge conservation --- i.e. they are a
direct sum of analogues of permutation representations. 
(Note for example that
the tensor space representation in \cite[\S4]{MartinSaleur94a} is 
neither full tilting \cite{Donkin93}
nor charge conserving.) 

One way to 
proceed is to search for maps from $b_n$ to $\TL_{n'}$ (some $n'$)
(resp. $\HD(n,d)$ to $H^d_{n'}$),
and hence obtain $b_n$--modules by restriction of $M^q_2$. 
Another possibility is
to enrich suitable $T_{n'}$--modules with the property of
$b_n$--module by determining an action of the blob operator on them.  
We begin by investigating the latter.

\subsection{Generalised bialgebra construction}\label{pseudo cabling map}

Let $(M,\circ,e)$ be a finite monoid, and $A$ a $K$--algebra 
with basis $M$ and multiplication defined on this basis by 
$m_1 m_2 = k_{m_1,m_2}(\qp) m_1 \circ m_2$ where $\qp$ is some set of
parameters and $k_{12}(\qp) \in K$. 
(The possibilities for the coefficients will in  general by 
constrained by $M$, 
$$m_1 (m_2 m_3) = m_1 k_{m_2,m_3}(\qp) (m_2 \circ m_3) = 
k_{m_2,m_3}(\qp) k_{m_1,m_2 \circ m_3}(\qp) m_1 \circ (m_2 \circ m_3)
$$ $$= k_{m_1,m_2}(\qp) (m_1 \circ m_2) m_3 = 
k_{m_1,m_2}(\qp) k_{m_1\circ m_2,m_3}(\qp) (m_1 \circ m_2) \circ m_3$$
so $k_{m_2,m_3}(\qp) k_{m_1,m_2 \circ m_3}(\qp)
= k_{m_1,m_2}(\qp) k_{m_1\circ m_2,m_3}(\qp)$,
but there are plentiful solutions --- for example, 
any finite group algebra.)
Suppose there is a triple of points in parameter space for which 
\eql(kk=k)
k_{m_1,m_2}(\qp')k_{m_1,m_2}(\qp'')=k_{m_1,m_2}(\qp)
\eq
for all
$m_1,m_2$, then $A$ has a kind of generalised coproduct:
$$A \hookrightarrow A' \times A''$$
$$ m \mapsto (m,m) $$
making it a kind of generalised bialgebra, since 
\[
\xymatrix{  m_1m_2 \ar@{=}[d] \ar@{|->}[dr] \\
k_{m_1,m_2}(\qp) m_1 \circ m_2 \ar@{|->}[d] 
                             & (m_1m_2,m_1m_2) \ar@{=}[d] \\
k_{m_1,m_2}(\qp) (m_1 \circ m_2,m_1 \circ m_2) \ar@{=}[r] 
               & k_{m_1,m_2}(\qp')k_{m_1,m_2}(\qp'') (m_1 \circ m_2,m_1 \circ m_2)}
\]
(cf. \cite[\S1.1.3(iv)]{Joseph95}, for example). In particular the
coproduct is an algebra morphism and 
if a submanifold $S$ of parameter space can be found for which each pair
$\qp', \qp'' \in S$ has a $\qp \in S$ satisfying equation(\ref{kk=k})
then the sum over all $\qp \in S$ 
of categories of (left) modules is closed under tensor product. 
(The example of group algebras is the usual bialgebra and tensor product.)

It is easy to show, using the diagram calculus (or via a mild
generalisation of 
the \ThL\ diagram variant of cabling 
\cite[\S A(iii)]{Martin89b}), that $T_n(q)$ is an
algebra of this type, with $\qp=\{ q \}$. 
The diagram in figure~\ref{comults}
illustrates the coproduct on $U_1 U_2 \in T_4(q)$, using lines of different
thickness for different $q$. 
The cabling--like visualisation of the two factors, in which they are
embedded in a single pseudodiagram, is possible because the
thin and thick lines are arranged into subdiagrams which never meet in any
composition. 
\begin{figure}
\[
\raisebox{-50pt}{\input{./xfig/comult102n.eepic}{}}
\; \hspace{.1in} \mapsto \hspace{.1in} 
\left( \raisebox{-50pt}{\xfigeps{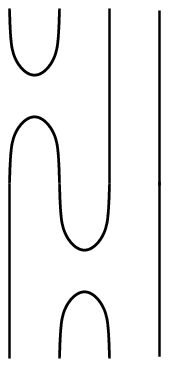}{}} \right. ,
\left. \raisebox{-50pt}{\xfigeps{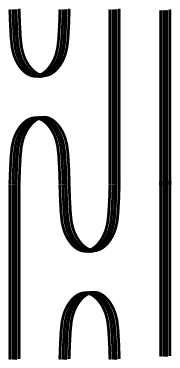}{}} \right)
\hspace{.1in} \sim \hspace{.1in} 
\raisebox{-60pt}{\input{./xfig/comult104xn.eepic}{}}
\]
\caption{\label{comults}}
\end{figure}
It will be evident that  
the set of conditions (\ref{kk=k}) include 
$-[2]_q = [2]_{q'} [2]_{q''}$ in this case 
(consider a composition in which a closed loop arises, such as
$U_1U_1$); and that this is the only  non--trivial condition. 
Through the cabling picture we may pass to another visualisation, in which
the thinner lines have been reflected in a vertical line at the
left edge of the diagram (cf. \cite{tomDieck94}), as illustrated in
figure~\ref{comults15}. 
\begin{figure}
\[
\raisebox{-50pt}{\xfigeps{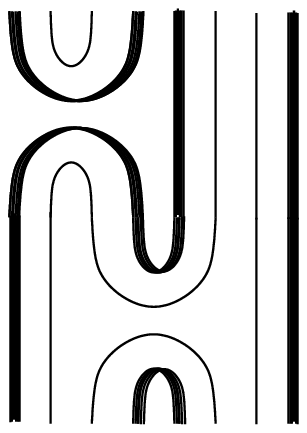}{}}
\; \hspace{.1in} \mapsto \hspace{.1in} 
\raisebox{-50pt}{\xfigeps{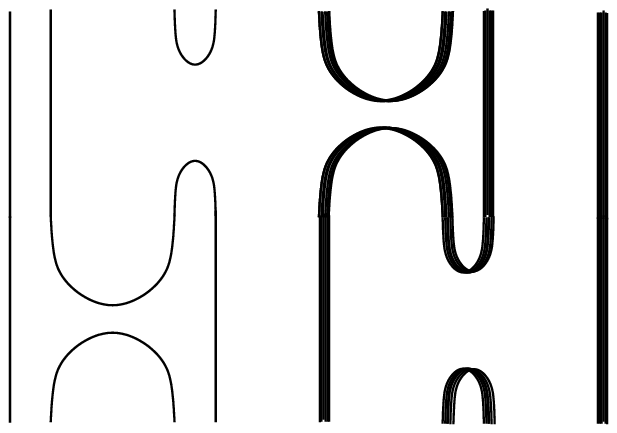}{}}
\]
\caption{\label{comults15}}
\end{figure}
There is no significant difference between these two visualisations,
except that it is perhaps slightly easier to 
describe the construction of certain tensor product
representations {\em explicitly} using the reflected form, as we will see.
Again in the reflected form 
we may view the picture as a single pseudodiagram (but
again 
there is no sense in which the right and
left hand sides ever touch). 
 
For example we may tensor together two tensor space representations in the form:
$$\rho_{} : \TL_n(q) \longrightarrow End(V_2^{2n})$$ 
$$\rho_{}: U_i \mapsto M_2^s(U_{n-i}) M_2^t(U_{n+i})$$
Here the set of conditions (\ref{kk=k}) reduce 
to $-[2]_q = [2]_s [2]_t$ 
via, e.g.,
\[
\xymatrix{ \rho_{}(U_i U_i )  \ar@{=}[d] \ar@{=}[r] 
                             & \rho_{}(U_i)\rho_{}(U_i) \\
 -[2]_q \rho_{}( U_i) \ar@{=}[r] 
                             & [2]_s [2]_t \rho_{}(U_i) . \ar@{=}[u] }  
\]

Set
$$\mU^q(x)=\mat{cccc} 
0&0&0&0\\ 
0&q&1&0 \\
0&1&q^{-1}&0 \\
0&0&0&x \tam$$
and $\mU^q =\mU^q(0)$ (cf. \cite{Deguchi89}). 
Just as $M_2^r(U_i)$ is a matrix acting trivially on every tensor
factor except the $i^{th}$ and $(i+1)^{th}$, where it acts as $-\mU^r$, so
let $M_2^{r,x}(U_i)$ denote a matrix differing from this only in acting
like $-\mU^r(x)$ in that position. 


Note
that 
\eql(rst) 
(\mU^{s}\otimes \mU^{t}) 
   (1 \otimes \mU^{r}(x) \otimes 1) 
      (\mU^{s}\otimes\mU^{t})
= 
\left( \frac{r}{st} +\frac{st}{r} +x \frac{t}{s} \right)  
  (\mU^{s}\otimes \mU^{t})
\eq
for any $r,s,t,x$ (an explicit calculation). 


\mpr{ \label{rho rep}
Fix $q,\ym$, 
put $q=e^{i\mu_q}$, $u_0=[\ym]_q e_-$, 
and choose $r,s,t$ such that
\eql(var1) - \cos(\mu_q)=2\cos(\mu_s)\cos(\mu_t) \eq
\eql(var2)
- \frac{\sin((\ym-1)\mu_q)}{\sin(\ym \mu_q)}
    = \frac{\cos(\mu_s + \mu_t -\mu_r)}{\cos(\mu_r)} . 
\eq
(NB, exclude $m=0$ and caveat $q=1$. A convenient realisation is 
$r=i(-q)^{\ym}$, 
$s=-i\sqrt{iq}$, $t=-\sqrt{iq}$, i.e. 
$\mu_r=m(\mu_q + \pi) +\frac{\pi}{2}$,
$\mu_s=\frac{\mu_q + \pi}{2}-\frac{3\pi}{4}$, 
$\mu_t=\frac{\mu_q + \pi}{2}+\frac{3\pi}{4}$ 
--- so rational $\frac{\mu_q}{\pi}$ and $m$ gives rational  
$\frac{\mu_r}{\pi}$, $\frac{\mu_s}{\pi}$,  $\frac{\mu_t}{\pi}$.) 
Then there is an algebra homomorphism 
$$\rho_{0} : b_n(q,\ym) \longrightarrow End(V_2^{2n})$$ 
given by 
\eql(mape)
\rho_{0}: e_- \mapsto \frac{1}{[2]_r} M_2^r(U_n)
\eq
$$\rho_{0}: U_i \mapsto M_2^s(U_{n-i}) M_2^t(U_{n+i}) . $$
}
\newline
{\em Proof:} We may readily verify $\rho_0(e_- e_- = e_-)$. 
%
%
The relations for $T_n(q) \subset b_n $ may be checked in $\rho_{}$, i.e. 
as above (NB $[m]_q=\frac{\sin(m\mu_q)}{\sin(\mu_q)}$). 
There remains $U_1 e_- U_1 \propto U_1$ (relation(\ref{blob3})). 
This
is validated by the
  explicit calculation in equation(\ref{rst}):
$$U_1 e_- U_1 = \frac{[m-1]_q}{[m]_q} U_1 
  \mapsto - \frac{\left( \frac{r}{st} +\frac{st}{r}
                             \right)}{[2]_r}\rho_{0}(U_1) . $$ 
\Qed

Another way to see this is to note that, in tensor space, 
equation~(\ref{rst}) allows us to make sense of an extension of the
pseudodiagrams in figure~\ref{comults15} 
in which the left and right hand sides meet as in
figure~\ref{comult16x2n2} (specifically, this figure may be replaced 
by a scalar multiple of 
one in which the loop composed of mixed line segments is omitted). 
\begin{figure}
\input{./xfig/comult16x2n2.eepic}
\caption{\label{comult16x2n2}}
\end{figure}

\mdef There is a similar homomorphism $\rho_s$ which simply replaces
equation(\ref{mape}) with 
\eql(mapes)
\rho_{s}: e_- \mapsto \frac{1}{[2]_r} M_2^{r,[2]_r}(U_n)
\eq
and equation(\ref{var2}) with 
\eql(var2s)
- \frac{\sin((m-1)\mu_q)}{\sin(m \mu_q)}
 = 
\frac{\cos(\mu_s + \mu_t -\mu_r)}{\cos(\mu_r)} +e^{i(\mu_t-\mu_s)} . 
\eq
This mild complication has the merit that the blob/box symmetry of the
algebra \cite{MartinSaleur94a} maps this representation to 
one of the same type. 

The isomorphism of $b_n$  to its opposite (note the symmetry of the
relations under writing back to front) provides an automorphism which
is fixed by $\rho_0$ (the representations of the generators are
symmetric matrices), thus $\rho_0$ is contravariant self--dual. 
The same is true of $\rho_s$. 

These very exciting representations merit further study. 
Cox, Martin and Ryom--Hansen have recently
shown
\cite{CoxMartinRyom01} that 
they are generically faithful. 
Other intriguing questions are: 
Are they (full) tilting? 
How do they generalize to higher $d$? 
The latter question is not trivial. We have made considerable use of
the \ThL\ diagram calculus here, and there is no such powerful tool
in evidence for higher $d$, short of the braid group itself. 
The remainder of the paper is essentially concerned with addressing
this question. 
\subsection{Other constructions} \label{Other constructions}

The blob algebra $b_n(q,\ym)$ is 
(at least) {\em singly critical} when $\ym \in \Z^+$, 
cf. section~\ref{decnos} --- if $q$ is not a root of unity then the
procedure there still works, but with $l$ set unboundedly large. 
This means that the usual Pascal triangle of standard modules 
\cite{MartinSaleur94a}
is complicated by at least one wall of `reflection' homomorphisms 
 from outside to inside across the alcove wall at $\ym$,
as exemplified in figure~\ref{pascal2}. 
There is a suggestive combinatorial coincidence with $\TL_n$ 
manifest in the dimension of heads of the blob standard modules in
certain singly critical (i.e. $\ym \in \Z$) cases. 
This starts with the $\ym=\pm 1$ cases, 
where there is a $b_n$ quotient given by $\blobe \mapsto 0$ (resp. 1). 
It follows immediately from the relations
that this quotient is isomorphic to 
$\TL_n$ (and so the coincidence is explained). 
Figure~\ref{pascal2} exhibits a similar phenomenon at $\ym=-2$. 
These may be regarded as special cases of a more general `braid construction'. 

In the next sections we describe this braid construction, and review aspects
of the connection between the blob algebra and periodic (and affine)
systems which lead to other useful maps from $b_n$ 
(and other $\AH(n,d)$ quotients) into ordinary Hecke quotients, and
hence into tensor space.

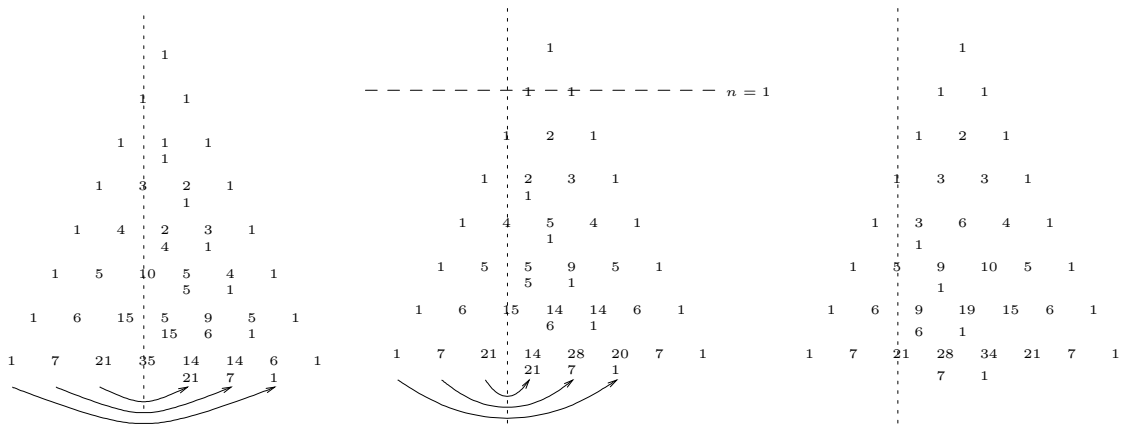
\begin{figure}
\input{./xfig/pascal11.eepic}
\input{./xfig/pascal2.eepic}
\input{./xfig/pascal3.eepic}
\caption{\label{pascal2} The Pascal triangle of $b_n$ standard modules
  (layers $n=0,1,...,7$ are shown) 
  complicated at $\ym \in \Z$ by reflection homomorphisms. 
  Here we show, by their dimensions, the composition factors of each
  standard in cases $\ym=-1,-2,-3$.}
\end{figure}

\section{Braids and the blob approach to periodic systems}\label{other reps}


In this paper the `$A$-type braid group' 
$\braid{}$ is the group of braidings of a row 
of initially vertical strings numbered from the top left: $1,2,...$, 
which braidings are trivial on all but finitely many strings. 
The subgroup $\braid{n}$ acts trivially on all but the first $n$
strings (so $\braid{0}=\braid{1} \subset \braid{2} \subset ...$). 
Thus $H_n$ is a quotient of $\C\braid{n}$, and 
$g_i$
is 
the braid in which string $i$ crosses over string
$i+1$. Evidently such elements, and inverses, generate the group 
\cite{MagnusKarrasSolitar66,Birman75}. 
\newcommand{\rotor}{(\!-\!)}
Let $\rotor : \braid{n} \rightarrow \braid{n}$ be the automorphism
given by $g_{i} \mapsto g_{n-i}$ 
\cite[\S5.7.2]{Martin91} (cf. \cite{tomDieck94}). 

Occasionally we shall need to refer to the group of braidings of precisely $n$
strings (as for example in the Young subgroup construction 
--- see \cite[\S13.1 p323]{Martin91}). 
This group is obviously isomorphic to $\braid{n}$, and we will
distinguish them only by context. 
Let $1^m$ denote the identity
element on precisely $m$ strings, and 
$$1^m_C: \braid{n} \rightarrow \braid{mn}$$ 
the corresponding {\em cabling} morphism (replace each string with $m$
parallel strings). 

Let 
$ \Delta: \braid{n} \rightarrow \braid{n}\times \braid{n}$ 
be the group comultiplication. 
Let 
$Y: \braid{n}\times \braid{n} \hookrightarrow \braid{2n}$ 
be the natural `Young' embedding, extended 
to (a version of) the full braid group by extending the numbering of
strings to $\Z\setminus\{ 0 \}$ 
--- i.e. essentially the full line not the half line --- 
and placing the second copy of $\braid{n}$ on the `minus' side. 
Let $F$ be the map back to the full braid group proper got by 
{\em folding} the 
left hand side of the plane over onto the right hand side at a point
slightly shifted from the origin (so that each negatively numbered
string starting point lies just to the left of its positively numbered
version and the system is bounded on the left again), 
then renumbering --- see figure~\ref{Millef1}. 
Let $S$ be the map back 
to $\braid{}$ 
got by renumbering $i \mapsto i+n+1$ ($i<0$)
and $i \mapsto i+n$ ($i>0$) and discarding all strings numbered less
than 1. 
\newcommand{\Ftilde}{1^2_F}
\newcommand{\Stilde}{1^2_S}
Now define a map $\Ftilde$ from  
$ \braid{n} \rightarrow \braid{2n}$ by commutativity of: 
\[ \xymatrix{
\braid{n}   \ar@{->}[d]^{\Delta} \ar@{->}[rr]^{\Ftilde} & &  \braid{2n}  \\
{\braid{n}\times \braid{n}} \ar@{->}[r]^{1\times\rotor}
                       &   {\braid{n}\times \braid{n}} \ar@{->}[r]^-{Y} 
                                         &  \braid{2n} \ar@{->}[u]^F
}
\]
and similarly for $\Stilde$.
The map $\Ftilde$  is similar also to the  $m=2$ cabling morphism 
in that each string now has a partner running  parallel to it, but the
over/under information is not the same. 
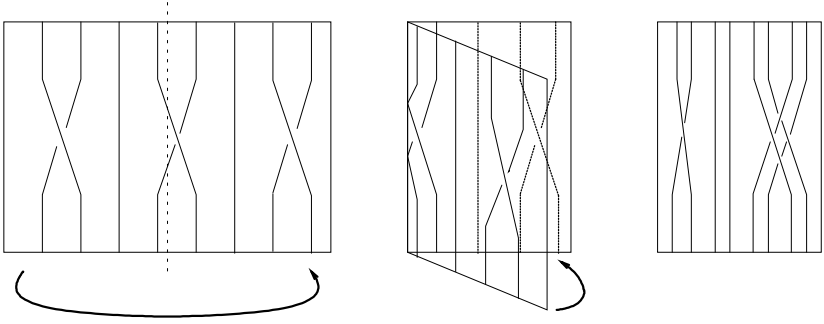
\begin{figure}
\input{./xfig/millef2.eepic}
\caption{\label{Millef1} Illustration of the folding injection of
  $\braid{n}$ into itself. }
\end{figure}

 

\mdef 
We now recall certain constructions
from \cite[\S5.7]{Martin91} (some
changes of notation will be necessary).
Let $\C\abraid{n}$ denote the algebra associated to affine graph 
$\hat{A}_n$  there
(strictly speaking only every other such graph has a pregraph, but
this need not concern us); 
let $\abraid{n}$ denote the underlying group 
 (i.e. the affinization
of the ordinary braid group); and let $g_.$ denote
the `extra' generator associated to the affinizing vertex in
$\hat{A}_n$ cf. $A_{n-1}$. 
Define 
$$G=g_1g_2\ldots g_{n-1}$$
in $\braid{n}$ (the element in which string 1 crosses over strings 2
to $n$); and note that $G g_i G^{-1}=g_{i+1}$ ($i<n-1$). 
\newline{\bf Proposition.}{\em 
\mbox{\rm \,\cite{Martin91}}
There is a homomorphism $\phi_0 : \abraid{n} \mapsto \braid{n}$ 
given by identification on the $\braid{n}$ subalgebra and 
$g_. \mapsto G g_{n-1} G^{-1}$. }


As noted for example 
in \cite[\S3]{Martin89b} there are actually a number of closely 
related ways of building
representations of the affine (or periodic) case, corresponding to the
choice of periodic boundary conditions (the `cohomological seam') 
in a physical system. 
This was systematized, in \cite[\S3]{MartinSaleur94a}, by the
introduction of the idempotent blob generator $\blobe$. 
\newcommand{\gb}{g_-}
\newcommand{\Gb}{G_-}
Using this one
builds an invertible generator $\gb=1+a\blobe$ ($a$ a suitable constant)
obeying 
\eql(blobg)
\gb g_1 \gb g_1 = g_1 \gb g_1 \gb
\eq
(see also proposition (\ref{->blob})) and defines $\Gb = \gb G$. 
Then 
\eql(bhom)
g_. \mapsto \Gb g_{n-1} (\Gb)^{-1}
\eq
 defines a generalisation $\phi_a$ of
$\phi_0$ for each suitable choice of $a,y_e$ 
(see \cite[eq.(25)]{MartinSaleur94a}). 

Note that neither equation(\ref{local}) nor the blob construction for $\gb$ are
needed to verify the map in eqn.(\ref{bhom}); only
equation(\ref{blobg}) is necessary. 
Thus the map generalises to $\aH{n}$, and even to the level of braids:


\mdef
The connection between the $B$--type and periodic
systems now follows, in as much as
$\bbraid{n}$ may be realized as the group of braids on the cylinder, 
whereupon
$\pi = g_{n-1} g_{n-2} \ldots g_1 g_0$ 
is the braid got from the identity braid by turning
the bottom edge of the cylinder through one vertex clockwise (i.e. so
as to take vertex 1 into vertex 2, and so on); 
and $\pi'$, the corresponding generalisation of $\Gb$ (i.e. with $\gb$
replaced by $g_0$), is the anticlockwise turn.
Thus $\bbraid{n}$ may be thought 
of as having affine $\abraid{n}$ as a subgroup --- 
with $g_. = \pi g_1\pi^{-1} = \pi' g_{n-1} (\pi')^{-1}$,
cf. \cite[\S3]{MartinSaleur94a}. 

\subsection{Homomorphisms of $\bbraid{n}$ to $\braid{n'}$}

\mdef
 Let $A$ = $\braid{n}$ and $\Nset{n}= \{1,2,...,n\}$. 
 For $b \in A$ and $i \in \Nset{n}$ 
 define $b(i)$ to be the final position of string $i$ in $b$. 

 Let $p$ be any partition of 
(equivalently, equivalence relation on) $\Nset{n}$. 
 Then for each such $p$ there is a subset of $A$ = $\braid{n}$ such that 
 $b \in A$ implies $b(i)\sim i$. 
 This subset is a subgroup --- call it $p-\braid{n}$, or just $p$--braid. 
 E.g. if $p$ is the `trivial' relation ($\{\Nset{n}\}$) then $p$--braid = $A$;
      if $p$ is the identity relation then $p$-braid is the pure braid
 group (the normal subgroup whose quotient is the permutation action
 on $\Nset{n}$, $b(i)$, described above).

 For convenience when dealing with general $n$ we will describe a partition 
 for which each $i > m$, some $m$, is in the same part by only giving the 
 other parts. Thus $\{\}$-braid = $\{\Nset{n}\}$-braid = $A$; and 
 $\{\!\{1\}\!\}$-braid$=:A'$ is the group which is pure on string 1.  


\mdef
Let $L_i \in \braid{n}$ be the pure braid which takes 
string $i$ behind all earlier strings 
then back in front of them (i.e. $L_1=1$ and $L_{i+1}$=$g_iL_ig_i$). 
See figure~\ref{Li}.
\begin{figure}
\input{./xfig/Li2.eepic}
\caption{\label{Li} $L_i$ --- here the shaded area represents the
  identity braid on the first $i-1$ strings.}
\end{figure}
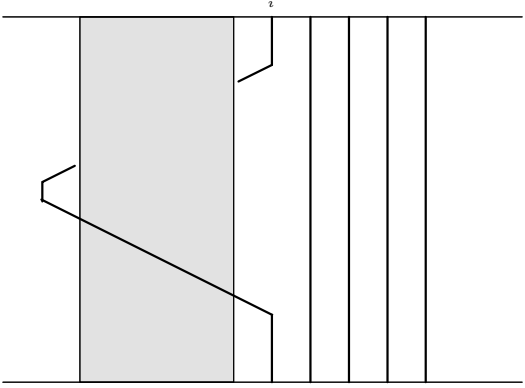
\begin{figure}
\xfigeps{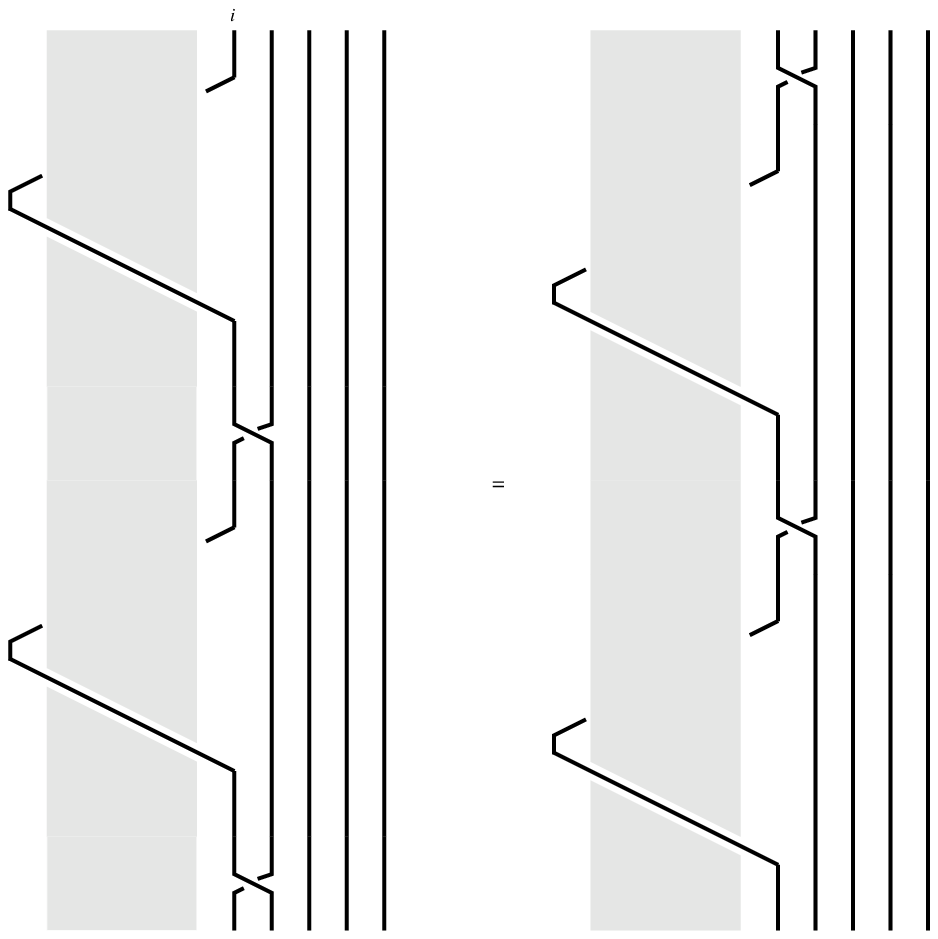}
\caption{\label{g0g1} 
Schematic for generalisations of 
$L_ig_iL_ig_i=L_iL_{i+1} = L_{i+1}L_i = g_iL_ig_iL_i$ as in $f_{b}$ 
($b$ lives in the shaded region).}
\end{figure}
Note that 
in the Hecke algebra quotient $L_i$ is the image of 
$X_i \in \AH(n,1) |_{\lambda_1=1}$; that 
$L_iL_j=L_jL_i$ and that $C_n=\prod_{i=1}^{n} L_i$ is the 
(`clockwise' or $g_i$--built) pure twist element of $\braid{n}$,
denoted $M^2$ in \cite[\S5.7.2]{Martin91}. 


\mdef
\mprop{\label{prop f_b}
Let $b$ be any element of $\braid{m-1}$. 
 Then there is an injective group homomorphism
 $$ f^m_b: \bbraid{n}
                  \; \longrightarrow \; \{\Nset{m-1}\}\!-\!\braid{n+m-1} $$
$$              g_0 \mapsto  C_m b $$
$$              g_i \mapsto  g_{i+m-1}   .   $$
}
\newline
{\em Proof:} consider figure~\ref{g0g1}. This checks the key relation
explicitly. 

For example, $f^2_1$ is $g_0 \mapsto g_1g_1$, $g_1 \mapsto g_2  $ 
and so on. $f^1_1$ is $g_0 \mapsto 1$, $g_i \mapsto g_i$ ($i>0$).



\mdef There is a group homomorphism extending the $1^2_C$ cabling
morphism (cf. \cite[Ch.13]{Martin91}; note also \cite{DateJimboMiwaOkado88}): 
$$ 1^2_C: \bbraid{n} 
                  \;\longrightarrow\; \braid{2n} $$
$$              g_0 \mapsto  g_1 $$
$$              g_i \mapsto  1^2_C(g_{i})   .   $$
There is a similar homomorphism extending the  $\Stilde$
morphism:
$$ \Stilde: \bbraid{n} 
                  \;\longrightarrow\; \braid{2n} $$
$$              g_0 \mapsto  g_n $$
$$              g_i \mapsto  g_{n-i} g_{n+i}     $$
(cf. figure~\ref{Millef1}); and another extending $\Ftilde$.   



\section{On tensor space representations of $\AH({n,d})$}\label{last s}

The constructions above allow us to build representations of
$\bbraid{n}$ from type--$A$ representations (and if these are charge
conserving tensor space representations then these properties will be
preserved, in some sense, as will globalisability). 
Our strategy now 
in searching for maps from $b_n$ to $\TL_n$ (and generalisations to
$d>2$) 
may be summarized by the following picture. 
\newcommand{\XX}{\mbox{${\cal X}$}}
\newcommand{\chippy}{\End(V^{\otimes n'})}
\[
\setlength{\unitlength}{0.00083333in}
\begingroup\makeatletter\ifx\SetFigFont\undefined%
\gdef\SetFigFont#1#2#3#4#5{%
  \reset@font\fontsize{#1}{#2pt}%
  \fontfamily{#3}\fontseries{#4}\fontshape{#5}%
  \selectfont}%
\fi\endgroup%
{\renewcommand{\dashlinestretch}{30}
\begin{picture}(3738,2820)(0,-10)
\dashline{60.000}(2700,2475)(3600,450)
\path(3523.849,547.473)(3600.000,450.000)(3578.678,571.842)
\path(450,600)(1200,1350)
\path(1136.360,1243.934)(1200.000,1350.000)(1093.934,1286.360)
\path(1650,1800)(2400,2550)
\path(2336.360,2443.934)(2400.000,2550.000)(2293.934,2486.360)
\path(600,225)(1500,225)
\path(1380.000,195.000)(1500.000,225.000)(1380.000,255.000)
\path(2400,225)(3300,225)
\path(3180.000,195.000)(3300.000,225.000)(3180.000,255.000)
\put(0,150){\makebox(0,0)[lb]{\smash{{{\SetFigFont{12}{14.4}{\rmdefault}{\mddefault}{\updefault}$\aH{n}$}}}}}
\put(1800,150){\makebox(0,0)[lb]{\smash{{{\SetFigFont{12}{14.4}{\rmdefault}{\mddefault}{\updefault}$H_{n'}$}}}}}
\put(3600,150){\makebox(0,0)[lb]{\smash{{{\SetFigFont{12}{14.4}{\rmdefault}{\mddefault}{\updefault}$\chippy$}}}}}
\put(1275,1425){\makebox(0,0)[lb]{\smash{{{\SetFigFont{12}{14.4}{\rmdefault}{\mddefault}{\updefault}$\AH(n,d)$}}}}}
\put(2475,2625){\makebox(0,0)[lb]{\smash{{{\SetFigFont{12}{14.4}{\rmdefault}{\mddefault}{\updefault}$b_n$}}}}}
\put(975,0){\makebox(0,0)[lb]{\smash{{{\SetFigFont{12}{14.4}{\rmdefault}{\mddefault}{\updefault}$\XX$}}}}}
\put(2575,0){\makebox(0,0)[lb]{\smash{{{\SetFigFont{12}{14.4}{\rmdefault}{\mddefault}{\updefault}${\cal  R}$}}}}}
\put(3225,1575){\makebox(0,0)[lb]{\smash{{{\SetFigFont{12}{14.4}{\rmdefault}{\mddefault}{\updefault}?}}}}}
\put(505,900){\makebox(0,0)[lb]{\smash{{{\SetFigFont{12}{14.4}{\rmdefault}{\mddefault}{\updefault}${\Psi_d}$}}}}}
\put(1775,2100){\makebox(0,0)[lb]{\smash{{{\SetFigFont{12}{14.4}{\rmdefault}{\mddefault}{\updefault}${\phi}$}}}}}
\end{picture}
}
\]
The Northeast pointing maps are the canonical quotients from sections
\ref{s1} and \ref{s3} (we will call the combined map $\eta$); 
the map ${\cal R}$ is the ordinary representation on
tensor space; and the dotted line is
the desired map occurring if ${\cal R} \circ \XX$ factors through $b_n$.   
The exercise is to find among the maps $\XX$ those for which 
${\cal R} \circ \XX$ so factors. 
(As we will see, the candidates we noted in 
\S\ref{Other constructions}
for maps from $b_n$ to $\TL_n$ are just
special cases of the simple representations $f^m_b$ of $\bbraid{n}$ above.) 

By (\ref{ordinaryTS}) each $\bbraid{}$ quotient of form 
\eql(aux) \bbraid{n} \stackrel{f^m_b}{\longrightarrow} \braid{m+n}
                 \stackrel{M_N}{\longrightarrow} V_N^{m+n} \eq
(`auxiliary space' construction) or
$$ \bbraid{n} \stackrel{1^2}{\longrightarrow} \braid{2n}
                 \stackrel{M_N}{\longrightarrow} V_N^{2n} $$
(`cabling related' construction)
factors through a partial specialisation of
$\AH(n,d)$ (some $d$) in which the $g_0$ (i.e. $X$) 
eigenvalues are determined in terms of $q$, but $q$ is indeterminate. 
Obviously the $M_2 \circ f^2_1$ quotient obeys a quadratic $f(g_0)=0$, 
so $d=2$, so it factors through some specialisation of Hecke$_B(n)$; 
but $f(g_0)$ has coefficients in $\Z[q,q^{-1}]$, so still this is not generic 
(Hecke$_B$ has two parameters, $q$ and $Q$, say).


Since all our maps $M_2 \circ \XX$ map into
the $\TL_n$ action on $V_2^n$ 
(for some $n$), their image breaks up at least as far as in eqn.(\ref{TLtilt}).


Possible next steps here are: 
(i) to investigate the generic irreducible content of 
 the $f_b$ representations (of whatever specialisation of $\AH(n,d)$ they 
 might provide); 
 and (ii) to investigate what portion of parameter space 
 is actually accessible by this construction 
 (i.e. what eigenvalues of $g_0$ we can realise by varying $b$). 
Fixing $q$,  
this portion is discrete, i.e. of measure zero, but then so is the
 (at least singly critical) portion we are most interested in, so  it is not
 necessarily too restricted. 

NB, if we want to access a dense subset of parameter space this {\em
  cannot} be via $M_N \circ f^m_b$, since this depends continuously
  only on $q$. Instead we could look at maps ending on, say,
  $\otimes_{i=1}^d V^n_N |_{q=q_i}$. 

\subsection{On cabling related maps}\label{cabling map}
Underlying the map $\Ftilde$ is the full range of direct product
representations of $\braid{n}$. In particular we can regard 
$\otimes_i M^{q_i}_N$ as a representation of $\braid{n}$. 
In general we may not assume that these representations will factor
through any particular Hecke quotient, but if one does then it could
provide a generalisation of the extension of  $\Ftilde$ to
$\bbraid{n}$. 
We form 
$R^{\otimes}(g_i) = M^{q}_N(g_{n-i})\otimes M^{r}_N(g_{n+i})$
and try, say, 
$R^{\otimes}(g_0) = 1\otimes M^{s}_N(g_{n-i})\otimes 1$
and compute $R^{\otimes}(g_0g_1g_0g_1-g_1g_0g_1g_0)$. 
In this particular case, with $N=3$, the image vanishes only when
$q=r=s$ (a brute force calculation). 
An analogous deformation of the extended $1^2_C$ map at $N=3$
fails in the same way. 

The $1^2_C$ cabling map with $N=2$ does not factor through $b_n$ in
general either. 
Further investigations are hindered by the magnitude of the
computations required, but these negative results serve well to
illustrate the extraordinary nature of the $\rho$--representations in
proposition(\ref{rho rep}). 

\subsection{On $N=2$ auxiliary spaces, $b_n$,  and `The coincidence'}

\newcommand{\makeR}{{\goth R}}

In (\ref{aux})
we require that $M_N \circ f^m_b (x) $ is also a representation ($\makeR$
say) of $\eta(x)$ for all $x$. 
Firstly, 
 $$\eta: g_i \mapsto \qa +   U_i  $$
so equation(\ref{TL3}) verifies equation(\ref{blob4}).  
We also require 
$\makeR(U_1 \blobe U_1 -  k_-U_1)=0$ for some $k_-$ 
(the relation (\ref{e02}) is 
not 
sufficient to ensure this). 
Whenever we {\em find} a map, the other question is: Is it faithful?


\mdef
Since 
$$ \eta: g_0 \mapsto \alpha 1 + \beta \blobe $$
the spectrum of $M_2(f^m_b(g_0))$ 
must be quadratic if the map is to factor through the blob as it stands. 
By (\ref{charges})
this spectrum may 
be determined from the action on the zero charge sector $\permm_0$ of $V_2^n$. 
The following lists
are the eigenvalues with multiplicities in this sector, arranged by
standard $\TL_n$--module factor with the `spine' module on the left
and so on. \newline 
$M_2(L_1)$: $\{ 1 \}$ \newline 
$M_2(L_2  = C_2)$: $\{ q^2, q^{-2} \}$ \newline 
$M_2(L_3)$: $\{ q^2, q^{-2}, q^{-4} \}$ $\;$
$M_2(C_3)$: $\{ 1, 1, q^{-6} \}$ \newline 
$M_2(L_4)$: $\{ 1, 1, q^2, q^{-4}, q^{-4}, q^{-6} \}$ $\;$
$M_2(C_4)$: $\{ 1, 1, q^{-4}, q^{-4}, q^{-4}, q^{-12} \}$ \newline 
$M_2(C_5)$: $\{ q^{-4}, q^{-4}, q^{-4}, q^{-4}, q^{-4}, 
                q^{-10}, q^{-10}, q^{-10}, q^{-10}, q^{-20} \}$ \newline 
$M_2(C_6)$: $\{ q^{-8},.., q^{-8}, q^{-12},.., q^{-12},  
                q^{-20}, .., q^{-20}, q^{-32} \}$ \newline 
The pattern for $C_i$ will be obvious. 

Note, therefore, that $C_1,C_2,C_3$ are the only possibilities
here (unless we further specialize to $q$ a root of unity). 
`Null' twist $C_1$ corresponds to the $m=\pm 1$ case already discussed. 
For the other cases it remains to check that the candidates for images of the
generators obey $U_1\blobe U_1 = k_-U_1$ for some scalar $k_-$. 


\mdef
Example: $M_2 \circ f^2_1$. 
(Let $u_1,u_0$ denote the canonical preimages along $\eta$ of $U_1$ and
$\blobe$.) 
An elementary calculation finds a value of $k_-$ for which 
$$ M_2(f^2_1( u_1 u_0 u_1 - k_- u_1 ) )=0 $$
This value then determines that the blob parameter $m=-2$ here. 

Since  span(1,$g_1$)=span(1,$g_1g_1$) here 
(and $V_2^n$ is a faithful $\TL_n$--module)
the image of $\bbraid{n}$ here is the whole of $\TL_{n+1}$. Thus
(\ref{charges}) determines the structure of $V^{n+1}_2$ as a
$\bbraid{n}$--module. 
For $n=2$ it is  
\eql(v23)
V_2^3 = 1+(2+1)+(2+1)+1 , 
\eq
(representing summands by their dimensions) 
which is the structure as a $\TL_3$--module. 
Since we hit the whole of $\TL_3$ 
equation(\ref{v23}) is the irreducible decomposition with $q$ generic.
Let us call the two inequivalent modules here $M^1,M^2$. 
The generic simples of $\AH(2,2)$, 
as indexed by their 2--partitions (see \S1), are 
$$ 
\begin{array}{|r|c|c|c|c|c|}  \hline
\mbox{2--partition} & ((2),) & ((1),(1)) & (,(2)) & ((1^2),) & (,(1^2)) \\ \hline
\mbox{dimension }   &  1     &    2        &1  &     1&           1 \\ \hline
\end{array}
$$
all but the last two of which survive the quotient to $b_2$. 
Note, therefore, that if the $\TL_3$ standard $M^2$ breaks up no further
(i.e. $q$--generically) it IS a
blob representation for
$m=-2$, 
but that $M_2 \circ f^2_1$ cannot be a faithful $b_n$--module. 
 
For $\bbraid{3}$ we have
$ V_2^4 = 1+(3+1)+(2+3+1)+(3+1)+1 $ (as a $\TL_4$--module) 
 cf. Hecke$_B$:
{\small $$
\begin{array}{|r|c|c|c|c|c|c|c|c|c|c|}  \hline
\mbox{2--partition}&((3),)&((2,1),)&((1^3),)&((2),(1))&((1^2),(1))&((1),(2))&((1),(1^2))&(,(3))&..
 \\ \hline
\mbox{dimension }   & 1* & 2 & 1 & 3* & 3 & 3* & 3 & 1* & .. \\ \hline
\end{array}
$$}
(blob representations indicated with a *).
Recall that at $m=-2$ blob standards break up as shown in figure~\ref{pascal2}.
In particular at $n=3$ 
$$\begin{array}{llccc} 1 &2 \; 3 \; 1   \\  \longrightarrow &1     \end{array}.$$
Note that here,
and for all $n$, the heads of the $b_n$ standards to  the right of the $m=-2$
line may indeed be identified with 
the (generic) irreducible $\TL_{n+1}$--modules.
This is neat, but it follows that none of these
representations are faithful. 

To summarize the last 2 sections,  we have not been successful in
generalising the $\rho$--representations. 
The search for full tilting modules for general $d$ continues, and 
we report these negative
results partly to avoid unnecessary duplication later. 
More positively, 
the representations
we have found are interesting from the point of view of Yang--Baxter
equations in Physics \cite{Baxter82}, but we will discuss these
applications elsewhere. 



\[ \]
\noindent
{\bf Acknowledgements} Thanks are due to S Dasmahapatra and the
members of the Donkin seminar for useful conversations during 1997/98, 
and EPSRC for funding part of this work under GRJ29923, GRJ55069
and GRM22536. 
DW has now left the subject, and PPM would like to thank Anton Cox and
Steen Ryom--Hansen for encouraging him to finish the paper after an
extended hiatus, and also for several useful conversations. 

I would like to commend the following recent papers on related topics
to the reader: \cite{SakamotoShoji99,Grojnowski99}.

\[ \]

\appendix 
\begin{center} {\LARGE Appendix} \end{center}
\section{On the Bernstein centre $Z(\aH{n})$ and $Z(\AH(n,d))$}
Some of the manipulations of ideals in \S\ref{s3} and thereafter 
are not trivial. 
The following
mechanical exposition of the Bernstein centre 
and its image in  $Z(\AH(n,d))$ may help the reader to
see where they come  from. 

\mdef Following on from equation(\ref{e012'}) 
define 
$\gh_i = [X_i , g_i ]$ 
(thus $\gh_i = (X_i - X_{i+1}) g_i + (q-q^{-1}) X_{i+1} $). Then 
\eql(ghi) \gh_i X_j = X_{\sigma_i(j)} \gh_i  \eq 
{\footnotesize 
\[ \gh_i \gh_{i+1} \gh_i 
= X_{i} g_{i} ( X_{i+1} g_{i+1} 
     - g_{i+1} X_{i+1}  ) X_{i} g_{i} + \cdots
= X_{i} g_{i}  X_{i} ( X_{i+1} g_{i+1}  ) g_{i}
     - X_{i} g_{i} ( g_{i+1}  g_{i} X_{i}  g_{i} ) X_{i} g_{i} + \cdots
\]
\[
=\cdots = X_{i+1} g_{i+1}    X_{i}    g_{i}^{-1}    X_{i+1}    g_{i+1}
     - g_{i+1}  X_{i+1}   X_{i} g_{i}^{-1}  g_{i+1}   X_{i+1}  +
     \cdots
\] 
\[
=  X_{i+1} g_{i+1}    X_{i}    (g_{i}-(q-q^{-1}))    X_{i+1}    g_{i+1}
     - g_{i+1}  X_{i+1}   X_{i}  (g_{i}-(q-q^{-1}))  g_{i+1}   X_{i+1}  + \cdots
\]
}
\eql(ghgh)
= X_{i+1} g_{i+1}    X_{i}    g_{i}    X_{i+1}    g_{i+1}
     - g_{i+1}  X_{i+1}   X_{i} g_{i}  g_{i+1}   X_{i+1}  + \cdots
= \gh_{i+1} \gh_i \gh_{i+1} . 
\eq
Let $K$ be our ground ring (an integral domain), 
$K[X_-]$ the ring of polynomials in the
$X_i$s, and $K(X_-)$ the field of fractions. 
Note that we can write 
$$ \alpha_i \gh_i =  g_i + \beta_i $$
where both $\alpha_i$ and $\beta_i$ lie in $K(X_-)$. 
(NB, in our quotient $\Psi_d$, the
image of $K[X_-]$ itself generically, but not always, contains 
$\alpha_i$ and $\beta_i$.
For example, when $d=1$ and $q^2=1$ then $X_i - X_{i+1}$ is not invertible.) 
It then follows from equation(\ref{ghi}) that any element of the extension
of $\aH{n}$ by $K(X_-)$ can be expressed in the form 
$$ h = \sum_{w \in \Bbas_n} c_w \gh_{w} $$ 
where $c_w \in K(X_-)$ and $\gh_w$ is obtained by
putting hats on the generators in $w$. 
Indeed the extension may be decomposed as 
\eql(decomposition)
\bigoplus_{w \in \Bbas_n} K(X_-) \gh_w
\eq
(an induction on the usual {\em  length} function on $\Bbas_n$). 
Suppose $h$ has at least one $w\neq 1$ with $c_w \neq 0$ (i.e. $h
\not\in K(X_-)$ subalgebra). 
Then there
is at least one $i$ such that $w(i)\neq i $ (under the obvious
generalisation of the $\sigma_i$ action in equation(\ref{ghi})) and 
$$ h X_i = \sum_w c_w \gh_w X_i = \sum_w c_w X_{w(i)} \gh_w
$$ 
so the $\gh_w$ component of $hX_i -X_i h$ is 
$c_w (X_{w(i)}-X_i) \neq 0$. 
Thus $Z_{K[X_-]}(\aH{n})=K[X_-] \supseteq Z(\aH{n})$. 
But with $c \in K[X_-]$, 
then $\gh_w c - c \gh_w = (c^w - c) \gh_w$ so $Z(\aH{n})=\asp$ as
Bernstein says. 

\mdef
Naturally $Z(\AH(n,d)) \supseteq \Psi_d(\asp)$, depending in principle
on the ground ring. 
The argument above mostly works in this case (to show equality), 
although the possible
specialisations of the ground ring become more restricted ($\lambda_i
\neq 0$, plus the restrictions already mentioned, for example). 
Note also that $\Psi_d(K[X_-])$ is finite dimensional. 
This makes it interesting to study $\Psi_d(\asp)$
--- an algebra which is, in a sense, more complicated that $\asp$
itself. 
For example, $\AH(2,1) \cong H_n$ and dim$(Z(H_2)=H_2)=2$,
so dim$(\Psi_1(\asp))=2$. Here a basis for $\Psi_1(\asp)$ is $\{
1,X+X_2 = \lam_1(1+ g_1^2) \}$ (note that this example illustrates the
problem with $q^2=1$). 

For another example, recall that dim$(Z(\AH(2,2)))=5$ generically,
so dim$(\Psi_1(\asp))=5$. Here a basis for $\Psi_2(\asp)$ is 
$\{ 1,X+X_2,XX_2,X^2+X_2^2,(X+X_2)XX_2  \}$ (this is not supposed to be
obvious!).
Another basis, convenient for comparison with the basis $ $ of $\AH(2,2)$, is 
\eql(basisII)
\{ 1,X+X_2,XX_2,(X+X_2-(\lam_1 + \lam_2))g_1, (XX_2 -\lam_1 \lam_2)g_1
\} .
\eq


\mdef
The set of monomial symmetric polynomials in two variables $X_1,X_2$ 
(a basis for $Z(\aH{2})$) may be
indexed by the set $\Lambda^2$ of Young diagrams of not more than two rows 
(write $m^a = (X^a)^{\Sigma}$, then 
$m^{(0)}=1$ and $m^{(1)}=X_1+X_2$, and so  on). 
NB, The set of such polynomials in which the degree of no individual
variable exceeds $d-1$ 
is {\em not} a basis for $\Psi_d(\asp)$ in general, 
as we see from the examples above.
\bibliographystyle{amsplain}
\bibliography{new31,main}
\end{document}

%% file: xfig/pKL.eepic
\setlength{\unitlength}{0.00083333in}
\begingroup\makeatletter\ifx\SetFigFont\undefined%
\gdef\SetFigFont#1#2#3#4#5{%
  \reset@font\fontsize{#1}{#2pt}%
  \fontfamily{#3}\fontseries{#4}\fontshape{#5}%
  \selectfont}%
\fi\endgroup%
{\renewcommand{\dashlinestretch}{30}
\begin{picture}(5357,4587)(0,-10)
\texture{8101010 10000000 444444 44000000 11101 11000000 444444 44000000 
	101010 10000000 444444 44000000 10101 1000000 444444 44000000 
	101010 10000000 444444 44000000 11101 11000000 444444 44000000 
	101010 10000000 444444 44000000 10101 1000000 444444 44000000 }
\shade\path(2175,2412)(2775,2412)(2775,1812)
	(2175,1812)(2175,2412)
\path(2175,2412)(2775,2412)(2775,1812)
	(2175,1812)(2175,2412)
\path(2250,3012)(2141,2957)(2250,2903)
	(2141,2848)(2250,2794)(2141,2739)
	(2250,2685)(2141,2630)(2250,2576)
	(2141,2521)(2250,2467)(2141,2412)
\path(2850,3612)(2741,3557)(2850,3503)
	(2741,3448)(2850,3394)(2741,3339)
	(2850,3285)(2741,3230)(2850,3176)
	(2741,3121)(2850,3067)(2741,3012)
\path(1650,3612)(1541,3557)(1650,3503)
	(1541,3448)(1650,3394)(1541,3339)
	(1650,3285)(1541,3230)(1650,3176)
	(1541,3121)(1650,3067)(1541,3012)
\path(3450,4212)(3341,4157)(3450,4103)
	(3341,4048)(3450,3994)(3341,3939)
	(3450,3885)(3341,3830)(3450,3776)
	(3341,3721)(3450,3667)(3341,3612)
\path(2250,4212)(2141,4157)(2250,4103)
	(2141,4048)(2250,3994)(2141,3939)
	(2250,3885)(2141,3830)(2250,3776)
	(2141,3721)(2250,3667)(2141,3612)
\path(1050,4212)(941,4157)(1050,4103)
	(941,4048)(1050,3994)(941,3939)
	(1050,3885)(941,3830)(1050,3776)
	(941,3721)(1050,3667)(941,3612)
\path(2850,1812)(2741,1757)(2850,1703)
	(2741,1648)(2850,1594)(2741,1539)
	(2850,1485)(2741,1430)(2850,1376)
	(2741,1321)(2850,1267)(2741,1212)
\path(3450,1212)(3341,1157)(3450,1103)
	(3341,1048)(3450,994)(3341,939)
	(3450,885)(3341,830)(3450,776)
	(3341,721)(3450,667)(3341,612)
\path(2250,1212)(2141,1157)(2250,1103)
	(2141,1048)(2250,994)(2141,939)
	(2250,885)(2141,830)(2250,776)
	(2141,721)(2250,667)(2141,612)
\path(2850,612)(2741,557)(2850,503)
	(2741,448)(2850,394)(2741,339)
	(2850,285)(2741,230)(2850,176)
	(2741,121)(2850,67)(2741,12)
\path(4050,612)(3941,557)(4050,503)
	(3941,448)(4050,394)(3941,339)
	(4050,285)(3941,230)(4050,176)
	(3941,121)(4050,67)(3941,12)
\path(1650,612)(1541,557)(1650,503)
	(1541,448)(1650,394)(1541,339)
	(1650,285)(1541,230)(1650,176)
	(1541,121)(1650,67)(1541,12)
\path(375,4212)(4575,4212)(4575,12)
	(375,12)(375,4212)
\path(375,3612)(4575,3612)
\path(4575,3012)(375,3012)
\path(375,2412)(4575,2412)
\path(375,1812)(4575,1812)
\path(375,1212)(4575,1212)
\path(375,612)(4575,612)
\path(975,4212)(975,12)
\path(1575,4212)(1575,12)
\path(2175,4212)(2175,12)
\path(2775,4212)(2775,12)
\path(3375,4212)(3375,12)
\path(3975,4212)(3975,12)
\path(4875,1512)(4876,1512)(4879,1511)
	(4884,1509)(4891,1506)(4901,1502)
	(4913,1497)(4927,1491)(4943,1484)
	(4960,1475)(4978,1465)(4997,1454)
	(5016,1441)(5035,1427)(5053,1411)
	(5072,1392)(5091,1371)(5109,1347)
	(5127,1319)(5144,1287)(5160,1251)
	(5175,1212)(5186,1174)(5196,1136)
	(5204,1100)(5210,1067)(5214,1037)
	(5218,1011)(5220,989)(5222,969)
	(5224,953)(5224,938)(5225,925)
	(5225,912)(5225,899)(5224,886)
	(5224,871)(5222,855)(5220,835)
	(5218,813)(5214,787)(5210,757)
	(5204,724)(5196,688)(5186,650)
	(5175,612)(5160,573)(5144,537)
	(5127,505)(5109,477)(5091,453)
	(5072,432)(5053,413)(5035,397)
	(5016,383)(4997,370)(4978,359)
	(4960,349)(4943,340)(4927,333)
	(4913,327)(4901,322)(4891,318)(4875,312)
\path(4976.826,382.225)(4875.000,312.000)(4997.893,326.045)
\path(2400,3462)(2401,3465)(2402,3471)
	(2404,3481)(2407,3496)(2411,3513)
	(2414,3532)(2417,3553)(2419,3574)
	(2419,3596)(2418,3618)(2415,3641)
	(2409,3665)(2400,3687)(2390,3704)
	(2380,3718)(2371,3729)(2365,3736)
	(2361,3740)(2358,3743)(2357,3743)
	(2356,3743)(2355,3743)(2352,3743)
	(2347,3743)(2339,3745)(2325,3748)
	(2306,3752)(2281,3757)(2250,3762)
	(2220,3766)(2189,3768)(2158,3770)
	(2128,3771)(2099,3772)(2071,3772)
	(2043,3771)(2016,3770)(1990,3769)
	(1965,3768)(1942,3766)(1921,3765)
	(1904,3764)(1875,3762)
\path(1948.274,3789.607)(1875.000,3762.000)(1951.370,3744.713)
\path(2550,3462)(2551,3465)(2553,3471)
	(2556,3481)(2561,3495)(2567,3514)
	(2572,3535)(2578,3558)(2583,3583)
	(2587,3608)(2589,3633)(2588,3659)
	(2585,3685)(2578,3711)(2566,3737)
	(2550,3762)(2529,3783)(2507,3800)
	(2485,3811)(2464,3817)(2444,3820)
	(2426,3819)(2410,3816)(2394,3812)
	(2378,3808)(2361,3804)(2344,3803)
	(2326,3803)(2306,3807)(2286,3813)
	(2266,3824)(2250,3837)(2240,3850)
	(2234,3862)(2231,3874)(2231,3884)
	(2232,3893)(2236,3900)(2240,3907)
	(2245,3913)(2250,3918)(2256,3924)
	(2262,3929)(2269,3935)(2276,3942)
	(2284,3950)(2292,3959)(2302,3968)
	(2313,3978)(2325,3987)(2344,3997)
	(2363,4003)(2381,4006)(2399,4006)
	(2415,4004)(2430,4001)(2445,3997)(2475,3987)
\path(2396.734,3989.372)(2475.000,3987.000)(2410.964,4032.062)
\put(1800,2637){\makebox(0,0)[lb]{\smash{{{\SetFigFont{12}{14.4}{\rmdefault}{\mddefault}{\updefault}1}}}}}
\put(3000,837){\makebox(0,0)[lb]{\smash{{{\SetFigFont{12}{14.4}{\rmdefault}{\mddefault}{\updefault}$v$}}}}}
\put(1200,237){\makebox(0,0)[lb]{\smash{{{\SetFigFont{12}{14.4}{\rmdefault}{\mddefault}{\updefault}$v$}}}}}
\put(1800,237){\makebox(0,0)[lb]{\smash{{{\SetFigFont{12}{14.4}{\rmdefault}{\mddefault}{\updefault}$v^2$}}}}}
\put(4200,237){\makebox(0,0)[lb]{\smash{{{\SetFigFont{12}{14.4}{\rmdefault}{\mddefault}{\updefault}1}}}}}
\put(3600,837){\makebox(0,0)[lb]{\smash{{{\SetFigFont{12}{14.4}{\rmdefault}{\mddefault}{\updefault}1}}}}}
\put(2400,837){\makebox(0,0)[lb]{\smash{{{\SetFigFont{12}{14.4}{\rmdefault}{\mddefault}{\updefault}$v^2$}}}}}
\put(1800,837){\makebox(0,0)[lb]{\smash{{{\SetFigFont{12}{14.4}{\rmdefault}{\mddefault}{\updefault}$v$}}}}}
\put(3000,1437){\makebox(0,0)[lb]{\smash{{{\SetFigFont{12}{14.4}{\rmdefault}{\mddefault}{\updefault}1}}}}}
\put(2400,2037){\makebox(0,0)[lb]{\smash{{{\SetFigFont{12}{14.4}{\rmdefault}{\mddefault}{\updefault}1}}}}}
\put(2400,2637){\makebox(0,0)[lb]{\smash{{{\SetFigFont{12}{14.4}{\rmdefault}{\mddefault}{\updefault}$v$}}}}}
\put(1800,3237){\makebox(0,0)[lb]{\smash{{{\SetFigFont{12}{14.4}{\rmdefault}{\mddefault}{\updefault}$v$}}}}}
\put(1200,3237){\makebox(0,0)[lb]{\smash{{{\SetFigFont{12}{14.4}{\rmdefault}{\mddefault}{\updefault}1}}}}}
\put(2400,3237){\makebox(0,0)[lb]{\smash{{{\SetFigFont{12}{14.4}{\rmdefault}{\mddefault}{\updefault}$v^2$}}}}}
\put(3000,3237){\makebox(0,0)[lb]{\smash{{{\SetFigFont{12}{14.4}{\rmdefault}{\mddefault}{\updefault}$v$}}}}}
\put(2250,237){\makebox(0,0)[lb]{\smash{{{\SetFigFont{12}{14.4}{\rmdefault}{\mddefault}{\updefault}$v^3\! +\! v$}}}}}
\put(2850,237){\makebox(0,0)[lb]{\smash{{{\SetFigFont{12}{14.4}{\rmdefault}{\mddefault}{\updefault}$1\! +\! v^2$}}}}}
\put(5250,912){\makebox(0,0)[lb]{\smash{{{\SetFigFont{12}{14.4}{\rmdefault}{\mddefault}{\updefault}-}}}}}
\put(3600,237){\makebox(0,0)[lb]{\smash{{{\SetFigFont{12}{14.4}{\rmdefault}{\mddefault}{\updefault}$v$}}}}}
\put(2400,1437){\makebox(0,0)[lb]{\smash{{{\SetFigFont{12}{14.4}{\rmdefault}{\mddefault}{\updefault}$v$}}}}}
\put(0,2037){\makebox(0,0)[lb]{\smash{{{\SetFigFont{12}{14.4}{\rmdefault}{\mddefault}{\updefault}0}}}}}
\put(0,2637){\makebox(0,0)[lb]{\smash{{{\SetFigFont{12}{14.4}{\rmdefault}{\mddefault}{\updefault}1}}}}}
\put(0,3237){\makebox(0,0)[lb]{\smash{{{\SetFigFont{12}{14.4}{\rmdefault}{\mddefault}{\updefault}2}}}}}
\put(0,3837){\makebox(0,0)[lb]{\smash{{{\SetFigFont{12}{14.4}{\rmdefault}{\mddefault}{\updefault}3}}}}}
\put(0,1437){\makebox(0,0)[lb]{\smash{{{\SetFigFont{12}{14.4}{\rmdefault}{\mddefault}{\updefault}-1}}}}}
\put(0,837){\makebox(0,0)[lb]{\smash{{{\SetFigFont{12}{14.4}{\rmdefault}{\mddefault}{\updefault}-2}}}}}
\put(0,237){\makebox(0,0)[lb]{\smash{{{\SetFigFont{12}{14.4}{\rmdefault}{\mddefault}{\updefault}-3}}}}}
\put(2400,4437){\makebox(0,0)[lb]{\smash{{{\SetFigFont{12}{14.4}{\rmdefault}{\mddefault}{\updefault}0}}}}}
\put(1800,4437){\makebox(0,0)[lb]{\smash{{{\SetFigFont{12}{14.4}{\rmdefault}{\mddefault}{\updefault}1}}}}}
\put(1200,4437){\makebox(0,0)[lb]{\smash{{{\SetFigFont{12}{14.4}{\rmdefault}{\mddefault}{\updefault}2}}}}}
\put(600,4437){\makebox(0,0)[lb]{\smash{{{\SetFigFont{12}{14.4}{\rmdefault}{\mddefault}{\updefault}3}}}}}
\put(3000,4437){\makebox(0,0)[lb]{\smash{{{\SetFigFont{12}{14.4}{\rmdefault}{\mddefault}{\updefault}-1}}}}}
\put(3600,4437){\makebox(0,0)[lb]{\smash{{{\SetFigFont{12}{14.4}{\rmdefault}{\mddefault}{\updefault}-2}}}}}
\put(4200,4437){\makebox(0,0)[lb]{\smash{{{\SetFigFont{12}{14.4}{\rmdefault}{\mddefault}{\updefault}-3}}}}}
\put(1800,3837){\makebox(0,0)[lb]{\smash{{{\SetFigFont{12}{14.4}{\rmdefault}{\mddefault}{\updefault}$v^2$}}}}}
\put(2475,3837){\makebox(0,0)[lb]{\smash{{{\SetFigFont{12}{14.4}{\rmdefault}{\mddefault}{\updefault}$v^3$}}}}}
\end{picture}
}

%% file: xfig/young2.latex
\setlength{\unitlength}{1579sp}%
\begingroup\makeatletter\ifx\SetFigFont\undefined%
\gdef\SetFigFont#1#2#3#4#5{%
  \reset@font\fontsize{#1}{#2pt}%
  \fontfamily{#3}\fontseries{#4}\fontshape{#5}%
  \selectfont}%
\fi\endgroup%
\begin{picture}(624,324)(1189,-73)
\thinlines
\put(1201,-61){\framebox(600,300){}}
\put(1501,239){\line( 0,-1){300}}
\end{picture}

%% file: xfig/young4.latex
\setlength{\unitlength}{1579sp}%
\begingroup\makeatletter\ifx\SetFigFont\undefined%
\gdef\SetFigFont#1#2#3#4#5{%
  \reset@font\fontsize{#1}{#2pt}%
  \fontfamily{#3}\fontseries{#4}\fontshape{#5}%
  \selectfont}%
\fi\endgroup%
\begin{picture}(1224,324)(1189,-73)
\thinlines
\special{ps: gsave 0 0 0 setrgbcolor}\put(1201,-61){\framebox(600,300){}}
\special{ps: gsave 0 0 0 setrgbcolor}\put(1501,239){\line( 0,-1){300}}
\special{ps: grestore}\special{ps: gsave 0 0 0 setrgbcolor}\put(1801,239){\line( 1, 0){300}}
\put(2101,239){\line( 0,-1){300}}
\put(2101,-61){\line(-1, 0){300}}
\special{ps: grestore}\special{ps: gsave 0 0 0 setrgbcolor}\put(2101,-61){\line( 1, 0){300}}
\put(2401,-61){\line( 0, 1){300}}
\put(2401,239){\line(-1, 0){300}}
\special{ps: grestore}\end{picture}

%% file: xfig/young3.latex
\setlength{\unitlength}{1579sp}%
\begingroup\makeatletter\ifx\SetFigFont\undefined%
\gdef\SetFigFont#1#2#3#4#5{%
  \reset@font\fontsize{#1}{#2pt}%
  \fontfamily{#3}\fontseries{#4}\fontshape{#5}%
  \selectfont}%
\fi\endgroup%
\begin{picture}(924,324)(1189,-73)
\thinlines
\special{ps: gsave 0 0 0 setrgbcolor}\put(1201,-61){\framebox(600,300){}}
\special{ps: gsave 0 0 0 setrgbcolor}\put(1501,239){\line( 0,-1){300}}
\special{ps: grestore}\special{ps: gsave 0 0 0 setrgbcolor}\put(1801,239){\line( 1, 0){300}}
\put(2101,239){\line( 0,-1){300}}
\put(2101,-61){\line(-1, 0){300}}
\special{ps: grestore}\end{picture}

%% file: xfig/young243.latex
\setlength{\unitlength}{1579sp}%
\begingroup\makeatletter\ifx\SetFigFont\undefined%
\gdef\SetFigFont#1#2#3#4#5{%
  \reset@font\fontsize{#1}{#2pt}%
  \fontfamily{#3}\fontseries{#4}\fontshape{#5}%
  \selectfont}%
\fi\endgroup%
\begin{picture}(1224,1524)(1189,-973)
\thinlines
\special{ps: gsave 0 0 0 setrgbcolor}\put(1201,-61){\framebox(600,300){}}
\special{ps: gsave 0 0 0 setrgbcolor}\put(1501,239){\line( 0,-1){300}}
\special{ps: grestore}\special{ps: gsave 0 0 0 setrgbcolor}\put(1201,-361){\framebox(1200,300){}}
\special{ps: gsave 0 0 0 setrgbcolor}\put(1501,-61){\line( 0,-1){300}}
\special{ps: grestore}\special{ps: gsave 0 0 0 setrgbcolor}\put(1801,-61){\line( 0,-1){300}}
\special{ps: grestore}\special{ps: gsave 0 0 0 setrgbcolor}\put(2101,-61){\line( 0,-1){300}}
\special{ps: grestore}\special{ps: gsave 0 0 0 setrgbcolor}\put(1201,-661){\framebox(900,300){}}
\special{ps: gsave 0 0 0 setrgbcolor}\put(1501,-361){\line( 0,-1){300}}
\special{ps: grestore}\special{ps: gsave 0 0 0 setrgbcolor}\put(1801,-361){\line( 0,-1){300}}
\special{ps: grestore}\special{ps: gsave 0 0 0 setrgbcolor}\multiput(2101,539)(0.00000,-272.72727){6}{\line( 0,-1){136.364}}
\special{ps: grestore}\end{picture}

%% file: xfig/A2eg_alt.eepic
\setlength{\unitlength}{0.00033333in}
\begingroup\makeatletter\ifx\SetFigFont\undefined%
\gdef\SetFigFont#1#2#3#4#5{%
  \reset@font\fontsize{#1}{#2pt}%
  \fontfamily{#3}\fontseries{#4}\fontshape{#5}%
  \selectfont}%
\fi\endgroup%
{\renewcommand{\dashlinestretch}{30}
\begin{picture}(2124,1839)(0,-10)
\put(420,522){\blacken\ellipse{150}{150}}
\put(420,522){\ellipse{150}{150}}
\thicklines
\path(1069,919)(841,530)(385,530)
\thinlines
\path(12,912)(2112,912)
\path(623.378,1723.463)(537.000,1812.000)(571.551,1693.230)
\path(537,1812)(1587,12)
\path(537,12)(1587,1812)
\path(1552.449,1693.230)(1587.000,1812.000)(1500.622,1723.463)
\end{picture}
}

%% file: xfig/young24.latex
\setlength{\unitlength}{1579sp}%
\begingroup\makeatletter\ifx\SetFigFont\undefined%
\gdef\SetFigFont#1#2#3#4#5{%
  \reset@font\fontsize{#1}{#2pt}%
  \fontfamily{#3}\fontseries{#4}\fontshape{#5}%
  \selectfont}%
\fi\endgroup%
\begin{picture}(1224,1224)(1189,-673)
\thinlines
\special{ps: gsave 0 0 0 setrgbcolor}\put(1201,-61){\framebox(600,300){}}
\special{ps: gsave 0 0 0 setrgbcolor}\put(1501,239){\line( 0,-1){300}}
\special{ps: grestore}\special{ps: gsave 0 0 0 setrgbcolor}\put(1201,-361){\framebox(1200,300){}}
\special{ps: gsave 0 0 0 setrgbcolor}\put(1501,-61){\line( 0,-1){300}}
\special{ps: grestore}\special{ps: gsave 0 0 0 setrgbcolor}\put(1801,-61){\line( 0,-1){300}}
\special{ps: grestore}\special{ps: gsave 0 0 0 setrgbcolor}\put(2101,-61){\line( 0,-1){300}}
\special{ps: grestore}\special{ps: gsave 0 0 0 setrgbcolor}\multiput(2401,539)(0.00000,-266.66667){5}{\line( 0,-1){133.333}}
\special{ps: grestore}\end{picture}

%% file: xfig/reflect1.eepic
\setlength{\unitlength}{0.00083333in}
\begingroup\makeatletter\ifx\SetFigFont\undefined%
\gdef\SetFigFont#1#2#3#4#5{%
  \reset@font\fontsize{#1}{#2pt}%
  \fontfamily{#3}\fontseries{#4}\fontshape{#5}%
  \selectfont}%
\fi\endgroup%
{\renewcommand{\dashlinestretch}{30}
\begin{picture}(3624,842)(0,-10)
\put(611,375){\blacken\ellipse{302}{302}}
\put(611,375){\ellipse{302}{302}}
\put(1811,375){\blacken\ellipse{302}{302}}
\put(1811,375){\ellipse{302}{302}}
\put(3011,375){\blacken\ellipse{302}{302}}
\put(3011,375){\ellipse{302}{302}}
\put(1211,375){\blacken\ellipse{152}{152}}
\put(1211,375){\ellipse{152}{152}}
\put(2411,375){\blacken\ellipse{152}{152}}
\put(2411,375){\ellipse{152}{152}}
\path(12,375)(3612,375)
\path(3012,600)(3009,602)(3002,607)
	(2990,615)(2973,627)(2951,641)
	(2927,657)(2900,675)(2873,692)
	(2846,708)(2820,723)(2796,736)
	(2772,748)(2750,758)(2728,767)
	(2706,775)(2685,781)(2662,787)
	(2643,792)(2623,796)(2603,800)
	(2581,803)(2559,806)(2536,809)
	(2512,811)(2488,812)(2463,814)
	(2438,814)(2412,815)(2386,814)
	(2361,814)(2336,812)(2312,811)
	(2288,809)(2265,806)(2243,803)
	(2221,800)(2201,796)(2181,792)
	(2162,787)(2139,781)(2118,775)
	(2096,767)(2074,758)(2052,748)
	(2028,736)(2004,723)(1978,708)
	(1951,692)(1924,675)(1897,657)
	(1873,641)(1851,627)(1812,600)
\blacken\path(1893.587,692.971)(1812.000,600.000)(1927.739,643.639)(1893.587,692.971)
\put(537,0){\makebox(0,0)[lb]{\smash{{{\SetFigFont{12}{14.4}{\rmdefault}{\mddefault}{\updefault}-2}}}}}
\put(1812,0){\makebox(0,0)[lb]{\smash{{{\SetFigFont{12}{14.4}{\rmdefault}{\mddefault}{\updefault}0}}}}}
\put(3012,0){\makebox(0,0)[lb]{\smash{{{\SetFigFont{12}{14.4}{\rmdefault}{\mddefault}{\updefault}2}}}}}
\end{picture}
}

%% file: xfig/comult102n.eepic
\setlength{\unitlength}{0.00033333in}
\begingroup\makeatletter\ifx\SetFigFont\undefined%
\gdef\SetFigFont#1#2#3#4#5{%
  \reset@font\fontsize{#1}{#2pt}%
  \fontfamily{#3}\fontseries{#4}\fontshape{#5}%
  \selectfont}%
\fi\endgroup%
{\renewcommand{\dashlinestretch}{30}
\begin{picture}(1918,4288)(0,-10)
\thicklines
\path(1270,4240)(1270,2140)
\path(70,4240)(70,4237)(70,4229)
	(70,4216)(71,4196)(72,4170)
	(73,4139)(74,4103)(75,4065)
	(77,4026)(79,3987)(81,3950)
	(83,3914)(86,3881)(89,3850)
	(92,3822)(95,3796)(99,3772)
	(104,3750)(109,3729)(114,3709)
	(120,3690)(127,3670)(135,3650)
	(144,3631)(153,3612)(164,3593)
	(175,3575)(187,3558)(200,3541)
	(213,3525)(227,3511)(241,3497)
	(256,3485)(270,3474)(285,3465)
	(300,3457)(314,3451)(328,3446)
	(342,3443)(356,3441)(370,3440)
	(384,3441)(398,3443)(412,3446)
	(426,3451)(440,3457)(455,3465)
	(470,3474)(484,3485)(499,3497)
	(513,3511)(527,3525)(540,3541)
	(553,3558)(565,3575)(576,3593)
	(587,3612)(596,3631)(605,3650)
	(613,3670)(620,3690)(626,3709)
	(631,3729)(636,3750)(641,3772)
	(645,3796)(648,3822)(651,3850)
	(654,3881)(657,3914)(659,3950)
	(661,3987)(663,4026)(665,4065)
	(666,4103)(667,4139)(668,4170)
	(669,4196)(670,4216)(670,4229)
	(670,4237)(670,4240)
\path(70,2140)(70,2143)(70,2151)
	(70,2164)(71,2184)(72,2210)
	(73,2241)(74,2277)(75,2315)
	(77,2354)(79,2393)(81,2430)
	(83,2466)(86,2499)(89,2530)
	(92,2558)(95,2584)(99,2608)
	(104,2630)(109,2651)(114,2671)
	(120,2690)(127,2710)(135,2730)
	(144,2749)(153,2768)(164,2787)
	(175,2805)(187,2822)(200,2839)
	(213,2855)(227,2869)(241,2883)
	(256,2895)(270,2906)(285,2915)
	(300,2923)(314,2929)(328,2934)
	(342,2937)(356,2939)(370,2940)
	(384,2939)(398,2937)(412,2934)
	(426,2929)(440,2923)(455,2915)
	(470,2906)(484,2895)(499,2883)
	(513,2869)(527,2855)(540,2839)
	(553,2822)(565,2805)(576,2787)
	(587,2768)(596,2749)(605,2730)
	(613,2710)(620,2690)(626,2671)
	(631,2651)(636,2630)(641,2608)
	(645,2584)(648,2558)(651,2530)
	(654,2499)(657,2466)(659,2430)
	(661,2393)(663,2354)(665,2315)
	(666,2277)(667,2241)(668,2210)
	(669,2184)(670,2164)(670,2151)
	(670,2143)(670,2140)
\path(1233,4233)(1233,2133)
\path(33,4233)(33,4230)(33,4222)
	(33,4209)(34,4189)(35,4163)
	(36,4132)(37,4096)(38,4058)
	(40,4019)(42,3980)(44,3943)
	(46,3907)(49,3874)(52,3843)
	(55,3815)(58,3789)(62,3765)
	(67,3743)(72,3722)(77,3702)
	(83,3683)(90,3663)(98,3643)
	(107,3624)(116,3605)(127,3586)
	(138,3568)(150,3551)(163,3534)
	(176,3518)(190,3504)(204,3490)
	(219,3478)(233,3467)(248,3458)
	(263,3450)(277,3444)(291,3439)
	(305,3436)(319,3434)(333,3433)
	(347,3434)(361,3436)(375,3439)
	(389,3444)(403,3450)(418,3458)
	(433,3467)(447,3478)(462,3490)
	(476,3504)(490,3518)(503,3534)
	(516,3551)(528,3568)(539,3586)
	(550,3605)(559,3624)(568,3643)
	(576,3663)(583,3683)(589,3702)
	(594,3722)(599,3743)(604,3765)
	(608,3789)(611,3815)(614,3843)
	(617,3874)(620,3907)(622,3943)
	(624,3980)(626,4019)(628,4058)
	(629,4096)(630,4132)(631,4163)
	(632,4189)(633,4209)(633,4222)
	(633,4230)(633,4233)
\path(33,2133)(33,2136)(33,2144)
	(33,2157)(34,2177)(35,2203)
	(36,2234)(37,2270)(38,2308)
	(40,2347)(42,2386)(44,2423)
	(46,2459)(49,2492)(52,2523)
	(55,2551)(58,2577)(62,2601)
	(67,2623)(72,2644)(77,2664)
	(83,2683)(90,2703)(98,2723)
	(107,2742)(116,2761)(127,2780)
	(138,2798)(150,2815)(163,2832)
	(176,2848)(190,2862)(204,2876)
	(219,2888)(233,2899)(248,2908)
	(263,2916)(277,2922)(291,2927)
	(305,2930)(319,2932)(333,2933)
	(347,2932)(361,2930)(375,2927)
	(389,2922)(403,2916)(418,2908)
	(433,2899)(447,2888)(462,2876)
	(476,2862)(490,2848)(503,2832)
	(516,2815)(528,2798)(539,2780)
	(550,2761)(559,2742)(568,2723)
	(576,2703)(583,2683)(589,2664)
	(594,2644)(599,2623)(604,2601)
	(608,2577)(611,2551)(614,2523)
	(617,2492)(620,2459)(622,2423)
	(624,2386)(626,2347)(628,2308)
	(629,2270)(630,2234)(631,2203)
	(632,2177)(633,2157)(633,2144)
	(633,2136)(633,2133)
\path(70,2140)(70,40)
\path(33,2133)(33,33)
\path(1885,2176)(1885,76)
\path(1848,2169)(1848,69)
\path(1885,4216)(1885,2116)
\path(1848,4209)(1848,2109)
\path(1270,2140)(1270,2137)(1270,2129)
	(1270,2116)(1269,2096)(1268,2070)
	(1267,2039)(1266,2003)(1265,1965)
	(1263,1926)(1261,1887)(1259,1850)
	(1257,1814)(1254,1781)(1251,1750)
	(1248,1722)(1245,1696)(1241,1672)
	(1236,1650)(1231,1629)(1226,1609)
	(1220,1590)(1213,1570)(1205,1550)
	(1196,1531)(1187,1512)(1176,1493)
	(1165,1475)(1153,1458)(1140,1441)
	(1127,1425)(1113,1411)(1099,1397)
	(1084,1385)(1070,1374)(1055,1365)
	(1040,1357)(1026,1351)(1012,1346)
	(998,1343)(984,1341)(970,1340)
	(956,1341)(942,1343)(928,1346)
	(914,1351)(900,1357)(885,1365)
	(870,1374)(856,1385)(841,1397)
	(827,1411)(813,1425)(800,1441)
	(787,1458)(775,1475)(764,1493)
	(753,1512)(744,1531)(735,1550)
	(727,1570)(720,1590)(714,1609)
	(709,1629)(704,1650)(699,1672)
	(695,1696)(692,1722)(689,1750)
	(686,1781)(683,1814)(681,1850)
	(679,1887)(677,1926)(675,1965)
	(674,2003)(673,2039)(672,2070)
	(671,2096)(670,2116)(670,2129)
	(670,2137)(670,2140)
\path(1270,40)(1270,43)(1270,51)
	(1270,64)(1269,84)(1268,110)
	(1267,141)(1266,177)(1265,215)
	(1263,254)(1261,293)(1259,330)
	(1257,366)(1254,399)(1251,430)
	(1248,458)(1245,484)(1241,508)
	(1236,530)(1231,551)(1226,571)
	(1220,590)(1213,610)(1205,630)
	(1196,649)(1187,668)(1176,687)
	(1165,705)(1153,722)(1140,739)
	(1127,755)(1113,769)(1099,783)
	(1084,795)(1070,806)(1055,815)
	(1040,823)(1026,829)(1012,834)
	(998,837)(984,839)(970,840)
	(956,839)(942,837)(928,834)
	(914,829)(900,823)(885,815)
	(870,806)(856,795)(841,783)
	(827,769)(813,755)(800,739)
	(787,722)(775,705)(764,687)
	(753,668)(744,649)(735,630)
	(727,610)(720,590)(714,571)
	(709,551)(704,530)(699,508)
	(695,484)(692,458)(689,430)
	(686,399)(683,366)(681,330)
	(679,293)(677,254)(675,215)
	(674,177)(673,141)(672,110)
	(671,84)(670,64)(670,51)
	(670,43)(670,40)
\path(1233,2133)(1233,2130)(1233,2122)
	(1233,2109)(1232,2089)(1231,2063)
	(1230,2032)(1229,1996)(1228,1958)
	(1226,1919)(1224,1880)(1222,1843)
	(1220,1807)(1217,1774)(1214,1743)
	(1211,1715)(1208,1689)(1204,1665)
	(1199,1643)(1194,1622)(1189,1602)
	(1183,1583)(1176,1563)(1168,1543)
	(1159,1524)(1150,1505)(1139,1486)
	(1128,1468)(1116,1451)(1103,1434)
	(1090,1418)(1076,1404)(1062,1390)
	(1047,1378)(1033,1367)(1018,1358)
	(1003,1350)(989,1344)(975,1339)
	(961,1336)(947,1334)(933,1333)
	(919,1334)(905,1336)(891,1339)
	(877,1344)(863,1350)(848,1358)
	(833,1367)(819,1378)(804,1390)
	(790,1404)(776,1418)(763,1434)
	(750,1451)(738,1468)(727,1486)
	(716,1505)(707,1524)(698,1543)
	(690,1563)(683,1583)(677,1602)
	(672,1622)(667,1643)(662,1665)
	(658,1689)(655,1715)(652,1743)
	(649,1774)(646,1807)(644,1843)
	(642,1880)(640,1919)(638,1958)
	(637,1996)(636,2032)(635,2063)
	(634,2089)(633,2109)(633,2122)
	(633,2130)(633,2133)
\path(1233,33)(1233,36)(1233,44)
	(1233,57)(1232,77)(1231,103)
	(1230,134)(1229,170)(1228,208)
	(1226,247)(1224,286)(1222,323)
	(1220,359)(1217,392)(1214,423)
	(1211,451)(1208,477)(1204,501)
	(1199,523)(1194,544)(1189,564)
	(1183,583)(1176,603)(1168,623)
	(1159,642)(1150,661)(1139,680)
	(1128,698)(1116,715)(1103,732)
	(1090,748)(1076,762)(1062,776)
	(1047,788)(1033,799)(1018,808)
	(1003,816)(989,822)(975,827)
	(961,830)(947,832)(933,833)
	(919,832)(905,830)(891,827)
	(877,822)(863,816)(848,808)
	(833,799)(819,788)(804,776)
	(790,762)(776,748)(763,732)
	(750,715)(738,698)(727,680)
	(716,661)(707,642)(698,623)
	(690,603)(683,583)(677,564)
	(672,544)(667,523)(662,501)
	(658,477)(655,451)(652,423)
	(649,392)(646,359)(644,323)
	(642,286)(640,247)(638,208)
	(637,170)(636,134)(635,103)
	(634,77)(633,57)(633,44)
	(633,36)(633,33)
\end{picture}
}

%% file: xfig/comult104xn.eepic
\setlength{\unitlength}{0.00027500in}
\begingroup\makeatletter\ifx\SetFigFont\undefined%
\gdef\SetFigFont#1#2#3#4#5{%
  \reset@font\fontsize{#1}{#2pt}%
  \fontfamily{#3}\fontseries{#4}\fontshape{#5}%
  \selectfont}%
\fi\endgroup%
{\renewcommand{\dashlinestretch}{30}
\begin{picture}(4266,6125)(0,-10)
\thicklines
\path(633,3115)(633,40)
\path(48,3063)(48,63)
\path(70,3070)(70,70)
\path(130,3070)(130,70)
\path(100,3100)(100,100)
\path(182,3070)(182,70)
\path(160,3070)(160,70)
\path(3648,3136)(3648,61)
\path(4233,3084)(4233,84)
\path(4211,3091)(4211,91)
\path(4151,3091)(4151,91)
\path(4181,3121)(4181,121)
\path(4099,3091)(4099,91)
\path(4121,3091)(4121,91)
\path(3648,6031)(3648,2956)
\path(4233,5979)(4233,2979)
\path(4211,5986)(4211,2986)
\path(4151,5986)(4151,2986)
\path(4181,6016)(4181,3016)
\path(4099,5986)(4099,2986)
\path(4121,5986)(4121,2986)
\path(3033,6040)(3033,3040)
\path(2433,6040)(2433,3040)
\path(2455,6047)(2455,3047)
\path(2515,6047)(2515,3047)
\path(2485,6077)(2485,3077)
\path(2567,6048)(2567,3048)
\path(2537,6048)(2537,3048)
\path(633,6040)(633,6037)(633,6029)
	(633,6016)(634,5996)(635,5970)
	(636,5939)(637,5903)(638,5865)
	(640,5826)(642,5787)(644,5750)
	(646,5714)(649,5681)(652,5650)
	(655,5622)(658,5596)(662,5572)
	(667,5550)(672,5529)(677,5509)
	(683,5490)(690,5470)(698,5450)
	(707,5431)(716,5412)(727,5393)
	(738,5375)(750,5358)(763,5341)
	(776,5325)(790,5311)(804,5297)
	(819,5285)(833,5274)(848,5265)
	(863,5257)(877,5251)(891,5246)
	(905,5243)(919,5241)(933,5240)
	(947,5241)(961,5243)(975,5246)
	(989,5251)(1003,5257)(1018,5265)
	(1033,5274)(1047,5285)(1062,5297)
	(1076,5311)(1090,5325)(1103,5341)
	(1116,5358)(1128,5375)(1139,5393)
	(1150,5412)(1159,5431)(1168,5450)
	(1176,5470)(1183,5490)(1189,5509)
	(1194,5529)(1199,5550)(1204,5572)
	(1208,5596)(1211,5622)(1214,5650)
	(1217,5681)(1220,5714)(1222,5750)
	(1224,5787)(1226,5826)(1228,5865)
	(1229,5903)(1230,5939)(1231,5970)
	(1232,5996)(1233,6016)(1233,6029)
	(1233,6037)(1233,6040)
\path(33,6040)(33,6037)(33,6032)
	(34,6021)(34,6006)(35,5984)
	(37,5956)(39,5922)(41,5884)
	(43,5841)(47,5796)(50,5750)
	(54,5703)(59,5656)(64,5611)
	(69,5567)(75,5526)(81,5486)
	(88,5449)(96,5414)(104,5381)
	(112,5350)(122,5321)(132,5293)
	(143,5266)(156,5240)(169,5215)
	(183,5190)(197,5167)(212,5145)
	(229,5122)(246,5100)(264,5078)
	(283,5056)(304,5034)(326,5013)
	(348,4992)(372,4971)(397,4950)
	(422,4930)(449,4911)(476,4892)
	(504,4875)(532,4858)(561,4842)
	(590,4827)(619,4813)(649,4800)
	(678,4789)(707,4778)(736,4769)
	(765,4762)(794,4755)(822,4750)
	(850,4745)(878,4742)(905,4741)
	(933,4740)(961,4741)(988,4742)
	(1016,4745)(1044,4750)(1072,4755)
	(1101,4762)(1130,4769)(1159,4778)
	(1188,4789)(1217,4800)(1247,4813)
	(1276,4827)(1305,4842)(1334,4858)
	(1362,4875)(1390,4892)(1417,4911)
	(1444,4930)(1469,4950)(1494,4971)
	(1518,4992)(1540,5013)(1562,5034)
	(1583,5056)(1602,5078)(1620,5100)
	(1637,5122)(1654,5145)(1669,5167)
	(1683,5190)(1697,5215)(1710,5240)
	(1723,5266)(1734,5293)(1744,5321)
	(1754,5350)(1762,5381)(1770,5414)
	(1778,5449)(1785,5486)(1791,5526)
	(1797,5567)(1802,5611)(1807,5656)
	(1812,5703)(1816,5750)(1819,5796)
	(1823,5841)(1825,5884)(1827,5922)
	(1829,5956)(1831,5984)(1832,6006)
	(1832,6021)(1833,6032)(1833,6037)(1833,6040)
\path(33,3040)(33,3043)(33,3048)
	(34,3059)(34,3074)(35,3096)
	(37,3124)(39,3158)(41,3196)
	(43,3239)(47,3284)(50,3330)
	(54,3377)(59,3424)(64,3469)
	(69,3513)(75,3554)(81,3594)
	(88,3631)(96,3666)(104,3699)
	(112,3730)(122,3759)(132,3787)
	(143,3814)(156,3840)(169,3865)
	(183,3890)(197,3913)(212,3935)
	(229,3958)(246,3980)(264,4002)
	(283,4024)(304,4046)(326,4067)
	(348,4088)(372,4109)(397,4130)
	(422,4150)(449,4169)(476,4188)
	(504,4205)(532,4222)(561,4238)
	(590,4253)(619,4267)(649,4280)
	(678,4291)(707,4302)(736,4311)
	(765,4318)(794,4325)(822,4330)
	(850,4335)(878,4338)(905,4339)
	(933,4340)(961,4339)(988,4338)
	(1016,4335)(1044,4330)(1072,4325)
	(1101,4318)(1130,4311)(1159,4302)
	(1188,4291)(1217,4280)(1247,4267)
	(1276,4253)(1305,4238)(1334,4222)
	(1362,4205)(1390,4188)(1417,4169)
	(1444,4150)(1469,4130)(1494,4109)
	(1518,4088)(1540,4067)(1562,4046)
	(1583,4024)(1602,4002)(1620,3980)
	(1637,3958)(1654,3935)(1669,3913)
	(1683,3890)(1697,3865)(1710,3840)
	(1723,3814)(1734,3787)(1744,3759)
	(1754,3730)(1762,3699)(1770,3666)
	(1778,3631)(1785,3594)(1791,3554)
	(1797,3513)(1802,3469)(1807,3424)
	(1812,3377)(1816,3330)(1819,3284)
	(1823,3239)(1825,3196)(1827,3158)
	(1829,3124)(1831,3096)(1832,3074)
	(1832,3059)(1833,3048)(1833,3043)(1833,3040)
\path(1233,3040)(1233,3037)(1233,3032)
	(1234,3021)(1234,3006)(1235,2984)
	(1237,2956)(1239,2922)(1241,2884)
	(1243,2841)(1247,2796)(1250,2750)
	(1254,2703)(1259,2656)(1264,2611)
	(1269,2567)(1275,2526)(1281,2486)
	(1288,2449)(1296,2414)(1304,2381)
	(1312,2350)(1322,2321)(1332,2293)
	(1343,2266)(1356,2240)(1369,2215)
	(1383,2190)(1397,2167)(1412,2145)
	(1429,2122)(1446,2100)(1464,2078)
	(1483,2056)(1504,2034)(1526,2013)
	(1548,1992)(1572,1971)(1597,1950)
	(1622,1930)(1649,1911)(1676,1892)
	(1704,1875)(1732,1858)(1761,1842)
	(1790,1827)(1819,1813)(1849,1800)
	(1878,1789)(1907,1778)(1936,1769)
	(1965,1762)(1994,1755)(2022,1750)
	(2050,1745)(2078,1742)(2105,1741)
	(2133,1740)(2161,1741)(2188,1742)
	(2216,1745)(2244,1750)(2272,1755)
	(2301,1762)(2330,1769)(2359,1778)
	(2388,1789)(2417,1800)(2447,1813)
	(2476,1827)(2505,1842)(2534,1858)
	(2562,1875)(2590,1892)(2617,1911)
	(2644,1930)(2669,1950)(2694,1971)
	(2718,1992)(2740,2013)(2762,2034)
	(2783,2056)(2802,2078)(2820,2100)
	(2837,2122)(2854,2145)(2869,2167)
	(2883,2190)(2897,2215)(2910,2240)
	(2923,2266)(2934,2293)(2944,2321)
	(2954,2350)(2962,2381)(2970,2414)
	(2978,2449)(2985,2486)(2991,2526)
	(2997,2567)(3002,2611)(3007,2656)
	(3012,2703)(3016,2750)(3019,2796)
	(3023,2841)(3025,2884)(3027,2922)
	(3029,2956)(3031,2984)(3032,3006)
	(3032,3021)(3033,3032)(3033,3037)(3033,3040)
\path(1233,40)(1233,43)(1233,48)
	(1234,59)(1234,74)(1235,96)
	(1237,124)(1239,158)(1241,196)
	(1243,239)(1247,284)(1250,330)
	(1254,377)(1259,424)(1264,469)
	(1269,513)(1275,554)(1281,594)
	(1288,631)(1296,666)(1304,699)
	(1312,730)(1322,759)(1332,787)
	(1343,814)(1356,840)(1369,865)
	(1383,890)(1397,913)(1412,935)
	(1429,958)(1446,980)(1464,1002)
	(1483,1024)(1504,1046)(1526,1067)
	(1548,1088)(1572,1109)(1597,1130)
	(1622,1150)(1649,1169)(1676,1188)
	(1704,1205)(1732,1222)(1761,1238)
	(1790,1253)(1819,1267)(1849,1280)
	(1878,1291)(1907,1302)(1936,1311)
	(1965,1318)(1994,1325)(2022,1330)
	(2050,1335)(2078,1338)(2105,1339)
	(2133,1340)(2161,1339)(2188,1338)
	(2216,1335)(2244,1330)(2272,1325)
	(2301,1318)(2330,1311)(2359,1302)
	(2388,1291)(2417,1280)(2447,1267)
	(2476,1253)(2505,1238)(2534,1222)
	(2562,1205)(2590,1188)(2617,1169)
	(2644,1150)(2669,1130)(2694,1109)
	(2718,1088)(2740,1067)(2762,1046)
	(2783,1024)(2802,1002)(2820,980)
	(2837,958)(2854,935)(2869,913)
	(2883,890)(2897,865)(2910,840)
	(2923,814)(2934,787)(2944,759)
	(2954,730)(2962,699)(2970,666)
	(2978,631)(2985,594)(2991,554)
	(2997,513)(3002,469)(3007,424)
	(3012,377)(3016,330)(3019,284)
	(3023,239)(3025,196)(3027,158)
	(3029,124)(3031,96)(3032,74)
	(3032,59)(3033,48)(3033,43)(3033,40)
\path(633,3040)(633,3043)(633,3051)
	(633,3064)(634,3084)(635,3110)
	(636,3141)(637,3177)(638,3215)
	(640,3254)(642,3293)(644,3330)
	(646,3366)(649,3399)(652,3430)
	(655,3458)(658,3484)(662,3508)
	(667,3530)(672,3551)(677,3571)
	(683,3590)(690,3610)(698,3630)
	(707,3649)(716,3668)(727,3687)
	(738,3705)(750,3722)(763,3739)
	(776,3755)(790,3769)(804,3783)
	(819,3795)(833,3806)(848,3815)
	(863,3823)(877,3829)(891,3834)
	(905,3837)(919,3839)(933,3840)
	(947,3839)(961,3837)(975,3834)
	(989,3829)(1003,3823)(1018,3815)
	(1033,3806)(1047,3795)(1062,3783)
	(1076,3769)(1090,3755)(1103,3739)
	(1116,3722)(1128,3705)(1139,3687)
	(1150,3668)(1159,3649)(1168,3630)
	(1176,3610)(1183,3590)(1189,3571)
	(1194,3551)(1199,3530)(1204,3508)
	(1208,3484)(1211,3458)(1214,3430)
	(1217,3399)(1220,3366)(1222,3330)
	(1224,3293)(1226,3254)(1228,3215)
	(1229,3177)(1230,3141)(1231,3110)
	(1232,3084)(1233,3064)(1233,3051)
	(1233,3043)(1233,3040)
\path(34,6032)(34,6029)(34,6024)
	(35,6013)(35,5998)(36,5976)
	(38,5948)(40,5914)(42,5876)
	(44,5833)(48,5788)(51,5742)
	(55,5695)(60,5648)(65,5603)
	(70,5559)(76,5518)(82,5478)
	(89,5441)(97,5406)(105,5373)
	(113,5342)(123,5313)(133,5285)
	(144,5258)(157,5232)(170,5207)
	(184,5182)(198,5159)(213,5137)
	(230,5114)(247,5092)(265,5070)
	(284,5048)(305,5026)(327,5005)
	(349,4984)(373,4963)(398,4942)
	(423,4922)(450,4903)(477,4884)
	(505,4867)(533,4850)(562,4834)
	(591,4819)(620,4805)(650,4792)
	(679,4781)(708,4770)(737,4761)
	(766,4754)(795,4747)(823,4742)
	(851,4737)(879,4734)(906,4733)
	(934,4732)(962,4733)(989,4734)
	(1017,4737)(1045,4742)(1073,4747)
	(1102,4754)(1131,4761)(1160,4770)
	(1189,4781)(1218,4792)(1248,4805)
	(1277,4819)(1306,4834)(1335,4850)
	(1363,4867)(1391,4884)(1418,4903)
	(1445,4922)(1470,4942)(1495,4963)
	(1519,4984)(1541,5005)(1563,5026)
	(1584,5048)(1603,5070)(1621,5092)
	(1638,5114)(1655,5137)(1670,5159)
	(1684,5182)(1698,5207)(1711,5232)
	(1724,5258)(1735,5285)(1745,5313)
	(1755,5342)(1763,5373)(1771,5406)
	(1779,5441)(1786,5478)(1792,5518)
	(1798,5559)(1803,5603)(1808,5648)
	(1813,5695)(1817,5742)(1820,5788)
	(1824,5833)(1826,5876)(1828,5914)
	(1830,5948)(1832,5976)(1833,5998)
	(1833,6013)(1834,6024)(1834,6029)(1834,6032)
\path(71,6039)(71,6036)(71,6031)
	(72,6020)(72,6005)(73,5983)
	(75,5955)(77,5921)(79,5883)
	(81,5840)(85,5795)(88,5749)
	(92,5702)(97,5655)(102,5610)
	(107,5566)(113,5525)(119,5485)
	(126,5448)(134,5413)(142,5380)
	(150,5349)(160,5320)(170,5292)
	(181,5265)(194,5239)(207,5214)
	(221,5189)(235,5166)(250,5144)
	(267,5121)(284,5099)(302,5077)
	(321,5055)(342,5033)(364,5012)
	(386,4991)(410,4970)(435,4949)
	(460,4929)(487,4910)(514,4891)
	(542,4874)(570,4857)(599,4841)
	(628,4826)(657,4812)(687,4799)
	(716,4788)(745,4777)(774,4768)
	(803,4761)(832,4754)(860,4749)
	(888,4744)(916,4741)(943,4740)
	(971,4739)(999,4740)(1026,4741)
	(1054,4744)(1082,4749)(1110,4754)
	(1139,4761)(1168,4768)(1197,4777)
	(1226,4788)(1255,4799)(1285,4812)
	(1314,4826)(1343,4841)(1372,4857)
	(1400,4874)(1428,4891)(1455,4910)
	(1482,4929)(1507,4949)(1532,4970)
	(1556,4991)(1578,5012)(1600,5033)
	(1621,5055)(1640,5077)(1658,5099)
	(1675,5121)(1692,5144)(1707,5166)
	(1721,5189)(1735,5214)(1748,5239)
	(1761,5265)(1772,5292)(1782,5320)
	(1792,5349)(1800,5380)(1808,5413)
	(1816,5448)(1823,5485)(1829,5525)
	(1835,5566)(1840,5610)(1845,5655)
	(1850,5702)(1854,5749)(1857,5795)
	(1861,5840)(1863,5883)(1865,5921)
	(1867,5955)(1869,5983)(1870,6005)
	(1870,6020)(1871,6031)(1871,6036)(1871,6039)
\path(87,6031)(87,6028)(87,6023)
	(88,6012)(88,5997)(89,5975)
	(91,5947)(93,5913)(95,5875)
	(97,5832)(101,5787)(104,5741)
	(108,5694)(113,5647)(118,5602)
	(123,5558)(129,5517)(135,5477)
	(142,5440)(150,5405)(158,5372)
	(166,5341)(176,5312)(186,5284)
	(197,5257)(210,5231)(223,5206)
	(237,5181)(251,5158)(266,5136)
	(283,5113)(300,5091)(318,5069)
	(337,5047)(358,5025)(380,5004)
	(402,4983)(426,4962)(451,4941)
	(476,4921)(503,4902)(530,4883)
	(558,4866)(586,4849)(615,4833)
	(644,4818)(673,4804)(703,4791)
	(732,4780)(761,4769)(790,4760)
	(819,4753)(848,4746)(876,4741)
	(904,4736)(932,4733)(959,4732)
	(987,4731)(1015,4732)(1042,4733)
	(1070,4736)(1098,4741)(1126,4746)
	(1155,4753)(1184,4760)(1213,4769)
	(1242,4780)(1271,4791)(1301,4804)
	(1330,4818)(1359,4833)(1388,4849)
	(1416,4866)(1444,4883)(1471,4902)
	(1498,4921)(1523,4941)(1548,4962)
	(1572,4983)(1594,5004)(1616,5025)
	(1637,5047)(1656,5069)(1674,5091)
	(1691,5113)(1708,5136)(1723,5158)
	(1737,5181)(1751,5206)(1764,5231)
	(1777,5257)(1788,5284)(1798,5312)
	(1808,5341)(1816,5372)(1824,5405)
	(1832,5440)(1839,5477)(1845,5517)
	(1851,5558)(1856,5602)(1861,5647)
	(1866,5694)(1870,5741)(1873,5787)
	(1877,5832)(1879,5875)(1881,5913)
	(1883,5947)(1885,5975)(1886,5997)
	(1886,6012)(1887,6023)(1887,6028)(1887,6031)
\path(124,6038)(124,6035)(124,6030)
	(125,6019)(125,6004)(126,5982)
	(128,5954)(130,5920)(132,5882)
	(134,5839)(138,5794)(141,5748)
	(145,5701)(150,5654)(155,5609)
	(160,5565)(166,5524)(172,5484)
	(179,5447)(187,5412)(195,5379)
	(203,5348)(213,5319)(223,5291)
	(234,5264)(247,5238)(260,5213)
	(274,5188)(288,5165)(303,5143)
	(320,5120)(337,5098)(355,5076)
	(374,5054)(395,5032)(417,5011)
	(439,4990)(463,4969)(488,4948)
	(513,4928)(540,4909)(567,4890)
	(595,4873)(623,4856)(652,4840)
	(681,4825)(710,4811)(740,4798)
	(769,4787)(798,4776)(827,4767)
	(856,4760)(885,4753)(913,4748)
	(941,4743)(969,4740)(996,4739)
	(1024,4738)(1052,4739)(1079,4740)
	(1107,4743)(1135,4748)(1163,4753)
	(1192,4760)(1221,4767)(1250,4776)
	(1279,4787)(1308,4798)(1338,4811)
	(1367,4825)(1396,4840)(1425,4856)
	(1453,4873)(1481,4890)(1508,4909)
	(1535,4928)(1560,4948)(1585,4969)
	(1609,4990)(1631,5011)(1653,5032)
	(1674,5054)(1693,5076)(1711,5098)
	(1728,5120)(1745,5143)(1760,5165)
	(1774,5188)(1788,5213)(1801,5238)
	(1814,5264)(1825,5291)(1835,5319)
	(1845,5348)(1853,5379)(1861,5412)
	(1869,5447)(1876,5484)(1882,5524)
	(1888,5565)(1893,5609)(1898,5654)
	(1903,5701)(1907,5748)(1910,5794)
	(1914,5839)(1916,5882)(1918,5920)
	(1920,5954)(1922,5982)(1923,6004)
	(1923,6019)(1924,6030)(1924,6035)(1924,6038)
\path(87,3049)(87,3052)(87,3057)
	(88,3068)(88,3083)(89,3105)
	(91,3133)(93,3167)(95,3205)
	(97,3248)(101,3293)(104,3339)
	(108,3386)(113,3433)(118,3478)
	(123,3522)(129,3563)(135,3603)
	(142,3640)(150,3675)(158,3708)
	(166,3739)(176,3768)(186,3796)
	(197,3823)(210,3849)(223,3874)
	(237,3899)(251,3922)(266,3944)
	(283,3967)(300,3989)(318,4011)
	(337,4033)(358,4055)(380,4076)
	(402,4097)(426,4118)(451,4139)
	(476,4159)(503,4178)(530,4197)
	(558,4214)(586,4231)(615,4247)
	(644,4262)(673,4276)(703,4289)
	(732,4300)(761,4311)(790,4320)
	(819,4327)(848,4334)(876,4339)
	(904,4344)(932,4347)(959,4348)
	(987,4349)(1015,4348)(1042,4347)
	(1070,4344)(1098,4339)(1126,4334)
	(1155,4327)(1184,4320)(1213,4311)
	(1242,4300)(1271,4289)(1301,4276)
	(1330,4262)(1359,4247)(1388,4231)
	(1416,4214)(1444,4197)(1471,4178)
	(1498,4159)(1523,4139)(1548,4118)
	(1572,4097)(1594,4076)(1616,4055)
	(1637,4033)(1656,4011)(1674,3989)
	(1691,3967)(1708,3944)(1723,3922)
	(1737,3899)(1751,3874)(1764,3849)
	(1777,3823)(1788,3796)(1798,3768)
	(1808,3739)(1816,3708)(1824,3675)
	(1832,3640)(1839,3603)(1845,3563)
	(1851,3522)(1856,3478)(1861,3433)
	(1866,3386)(1870,3339)(1873,3293)
	(1877,3248)(1879,3205)(1881,3167)
	(1883,3133)(1885,3105)(1886,3083)
	(1886,3068)(1887,3057)(1887,3052)(1887,3049)
\path(124,3042)(124,3045)(124,3050)
	(125,3061)(125,3076)(126,3098)
	(128,3126)(130,3160)(132,3198)
	(134,3241)(138,3286)(141,3332)
	(145,3379)(150,3426)(155,3471)
	(160,3515)(166,3556)(172,3596)
	(179,3633)(187,3668)(195,3701)
	(203,3732)(213,3761)(223,3789)
	(234,3816)(247,3842)(260,3867)
	(274,3892)(288,3915)(303,3937)
	(320,3960)(337,3982)(355,4004)
	(374,4026)(395,4048)(417,4069)
	(439,4090)(463,4111)(488,4132)
	(513,4152)(540,4171)(567,4190)
	(595,4207)(623,4224)(652,4240)
	(681,4255)(710,4269)(740,4282)
	(769,4293)(798,4304)(827,4313)
	(856,4320)(885,4327)(913,4332)
	(941,4337)(969,4340)(996,4341)
	(1024,4342)(1052,4341)(1079,4340)
	(1107,4337)(1135,4332)(1163,4327)
	(1192,4320)(1221,4313)(1250,4304)
	(1279,4293)(1308,4282)(1338,4269)
	(1367,4255)(1396,4240)(1425,4224)
	(1453,4207)(1481,4190)(1508,4171)
	(1535,4152)(1560,4132)(1585,4111)
	(1609,4090)(1631,4069)(1653,4048)
	(1674,4026)(1693,4004)(1711,3982)
	(1728,3960)(1745,3937)(1760,3915)
	(1774,3892)(1788,3867)(1801,3842)
	(1814,3816)(1825,3789)(1835,3761)
	(1845,3732)(1853,3701)(1861,3668)
	(1869,3633)(1876,3596)(1882,3556)
	(1888,3515)(1893,3471)(1898,3426)
	(1903,3379)(1907,3332)(1910,3286)
	(1914,3241)(1916,3198)(1918,3160)
	(1920,3126)(1922,3098)(1923,3076)
	(1923,3061)(1924,3050)(1924,3045)(1924,3042)
\path(34,3048)(34,3051)(34,3056)
	(35,3067)(35,3082)(36,3104)
	(38,3132)(40,3166)(42,3204)
	(44,3247)(48,3292)(51,3338)
	(55,3385)(60,3432)(65,3477)
	(70,3521)(76,3562)(82,3602)
	(89,3639)(97,3674)(105,3707)
	(113,3738)(123,3767)(133,3795)
	(144,3822)(157,3848)(170,3873)
	(184,3898)(198,3921)(213,3943)
	(230,3966)(247,3988)(265,4010)
	(284,4032)(305,4054)(327,4075)
	(349,4096)(373,4117)(398,4138)
	(423,4158)(450,4177)(477,4196)
	(505,4213)(533,4230)(562,4246)
	(591,4261)(620,4275)(650,4288)
	(679,4299)(708,4310)(737,4319)
	(766,4326)(795,4333)(823,4338)
	(851,4343)(879,4346)(906,4347)
	(934,4348)(962,4347)(989,4346)
	(1017,4343)(1045,4338)(1073,4333)
	(1102,4326)(1131,4319)(1160,4310)
	(1189,4299)(1218,4288)(1248,4275)
	(1277,4261)(1306,4246)(1335,4230)
	(1363,4213)(1391,4196)(1418,4177)
	(1445,4158)(1470,4138)(1495,4117)
	(1519,4096)(1541,4075)(1563,4054)
	(1584,4032)(1603,4010)(1621,3988)
	(1638,3966)(1655,3943)(1670,3921)
	(1684,3898)(1698,3873)(1711,3848)
	(1724,3822)(1735,3795)(1745,3767)
	(1755,3738)(1763,3707)(1771,3674)
	(1779,3639)(1786,3602)(1792,3562)
	(1798,3521)(1803,3477)(1808,3432)
	(1813,3385)(1817,3338)(1820,3292)
	(1824,3247)(1826,3204)(1828,3166)
	(1830,3132)(1832,3104)(1833,3082)
	(1833,3067)(1834,3056)(1834,3051)(1834,3048)
\path(71,3041)(71,3044)(71,3049)
	(72,3060)(72,3075)(73,3097)
	(75,3125)(77,3159)(79,3197)
	(81,3240)(85,3285)(88,3331)
	(92,3378)(97,3425)(102,3470)
	(107,3514)(113,3555)(119,3595)
	(126,3632)(134,3667)(142,3700)
	(150,3731)(160,3760)(170,3788)
	(181,3815)(194,3841)(207,3866)
	(221,3891)(235,3914)(250,3936)
	(267,3959)(284,3981)(302,4003)
	(321,4025)(342,4047)(364,4068)
	(386,4089)(410,4110)(435,4131)
	(460,4151)(487,4170)(514,4189)
	(542,4206)(570,4223)(599,4239)
	(628,4254)(657,4268)(687,4281)
	(716,4292)(745,4303)(774,4312)
	(803,4319)(832,4326)(860,4331)
	(888,4336)(916,4339)(943,4340)
	(971,4341)(999,4340)(1026,4339)
	(1054,4336)(1082,4331)(1110,4326)
	(1139,4319)(1168,4312)(1197,4303)
	(1226,4292)(1255,4281)(1285,4268)
	(1314,4254)(1343,4239)(1372,4223)
	(1400,4206)(1428,4189)(1455,4170)
	(1482,4151)(1507,4131)(1532,4110)
	(1556,4089)(1578,4068)(1600,4047)
	(1621,4025)(1640,4003)(1658,3981)
	(1675,3959)(1692,3936)(1707,3914)
	(1721,3891)(1735,3866)(1748,3841)
	(1761,3815)(1772,3788)(1782,3760)
	(1792,3731)(1800,3700)(1808,3667)
	(1816,3632)(1823,3595)(1829,3555)
	(1835,3514)(1840,3470)(1845,3425)
	(1850,3378)(1854,3331)(1857,3285)
	(1861,3240)(1863,3197)(1865,3159)
	(1867,3125)(1869,3097)(1870,3075)
	(1870,3060)(1871,3049)(1871,3044)(1871,3041)
\path(1945,3047)(1945,3044)(1945,3036)
	(1945,3023)(1946,3003)(1947,2977)
	(1948,2946)(1949,2910)(1950,2872)
	(1952,2833)(1954,2794)(1956,2757)
	(1958,2721)(1961,2688)(1964,2657)
	(1967,2629)(1970,2603)(1974,2579)
	(1979,2557)(1984,2536)(1989,2516)
	(1995,2497)(2002,2477)(2010,2457)
	(2019,2438)(2028,2419)(2039,2400)
	(2050,2382)(2062,2365)(2075,2348)
	(2088,2332)(2102,2318)(2116,2304)
	(2131,2292)(2145,2281)(2160,2272)
	(2175,2264)(2189,2258)(2203,2253)
	(2217,2250)(2231,2248)(2245,2247)
	(2259,2248)(2273,2250)(2287,2253)
	(2301,2258)(2315,2264)(2330,2272)
	(2345,2281)(2359,2292)(2374,2304)
	(2388,2318)(2402,2332)(2415,2348)
	(2428,2365)(2440,2382)(2451,2400)
	(2462,2419)(2471,2438)(2480,2457)
	(2488,2477)(2495,2497)(2501,2516)
	(2506,2536)(2511,2557)(2516,2579)
	(2520,2603)(2523,2629)(2526,2657)
	(2529,2688)(2532,2721)(2534,2757)
	(2536,2794)(2538,2833)(2540,2872)
	(2541,2910)(2542,2946)(2543,2977)
	(2544,3003)(2545,3023)(2545,3036)
	(2545,3044)(2545,3047)
\path(1833,3040)(1833,3037)(1833,3029)
	(1833,3016)(1834,2996)(1835,2970)
	(1836,2939)(1837,2903)(1838,2865)
	(1840,2826)(1842,2787)(1844,2750)
	(1846,2714)(1849,2681)(1852,2650)
	(1855,2622)(1858,2596)(1862,2572)
	(1867,2550)(1872,2529)(1877,2509)
	(1883,2490)(1890,2470)(1898,2450)
	(1907,2431)(1916,2412)(1927,2393)
	(1938,2375)(1950,2358)(1963,2341)
	(1976,2325)(1990,2311)(2004,2297)
	(2019,2285)(2033,2274)(2048,2265)
	(2063,2257)(2077,2251)(2091,2246)
	(2105,2243)(2119,2241)(2133,2240)
	(2147,2241)(2161,2243)(2175,2246)
	(2189,2251)(2203,2257)(2218,2265)
	(2233,2274)(2247,2285)(2262,2297)
	(2276,2311)(2290,2325)(2303,2341)
	(2316,2358)(2328,2375)(2339,2393)
	(2350,2412)(2359,2431)(2368,2450)
	(2376,2470)(2383,2490)(2389,2509)
	(2394,2529)(2399,2550)(2404,2572)
	(2408,2596)(2411,2622)(2414,2650)
	(2417,2681)(2420,2714)(2422,2750)
	(2424,2787)(2426,2826)(2428,2865)
	(2429,2903)(2430,2939)(2431,2970)
	(2432,2996)(2433,3016)(2433,3029)
	(2433,3037)(2433,3040)
\path(1870,3040)(1870,3037)(1870,3029)
	(1870,3016)(1871,2996)(1872,2970)
	(1873,2939)(1874,2903)(1875,2865)
	(1877,2826)(1879,2787)(1881,2750)
	(1883,2714)(1886,2681)(1889,2650)
	(1892,2622)(1895,2596)(1899,2572)
	(1904,2550)(1909,2529)(1914,2509)
	(1920,2490)(1927,2470)(1935,2450)
	(1944,2431)(1953,2412)(1964,2393)
	(1975,2375)(1987,2358)(2000,2341)
	(2013,2325)(2027,2311)(2041,2297)
	(2056,2285)(2070,2274)(2085,2265)
	(2100,2257)(2114,2251)(2128,2246)
	(2142,2243)(2156,2241)(2170,2240)
	(2184,2241)(2198,2243)(2212,2246)
	(2226,2251)(2240,2257)(2255,2265)
	(2270,2274)(2284,2285)(2299,2297)
	(2313,2311)(2327,2325)(2340,2341)
	(2353,2358)(2365,2375)(2376,2393)
	(2387,2412)(2396,2431)(2405,2450)
	(2413,2470)(2420,2490)(2426,2509)
	(2431,2529)(2436,2550)(2441,2572)
	(2445,2596)(2448,2622)(2451,2650)
	(2454,2681)(2457,2714)(2459,2750)
	(2461,2787)(2463,2826)(2465,2865)
	(2466,2903)(2467,2939)(2468,2970)
	(2469,2996)(2470,3016)(2470,3029)
	(2470,3037)(2470,3040)
\path(1908,3040)(1908,3037)(1908,3029)
	(1908,3016)(1909,2996)(1910,2970)
	(1911,2939)(1912,2903)(1913,2865)
	(1915,2826)(1917,2787)(1919,2750)
	(1921,2714)(1924,2681)(1927,2650)
	(1930,2622)(1933,2596)(1937,2572)
	(1942,2550)(1947,2529)(1952,2509)
	(1958,2490)(1965,2470)(1973,2450)
	(1982,2431)(1991,2412)(2002,2393)
	(2013,2375)(2025,2358)(2038,2341)
	(2051,2325)(2065,2311)(2079,2297)
	(2094,2285)(2108,2274)(2123,2265)
	(2138,2257)(2152,2251)(2166,2246)
	(2180,2243)(2194,2241)(2208,2240)
	(2222,2241)(2236,2243)(2250,2246)
	(2264,2251)(2278,2257)(2293,2265)
	(2308,2274)(2322,2285)(2337,2297)
	(2351,2311)(2365,2325)(2378,2341)
	(2391,2358)(2403,2375)(2414,2393)
	(2425,2412)(2434,2431)(2443,2450)
	(2451,2470)(2458,2490)(2464,2509)
	(2469,2529)(2474,2550)(2479,2572)
	(2483,2596)(2486,2622)(2489,2650)
	(2492,2681)(2495,2714)(2497,2750)
	(2499,2787)(2501,2826)(2503,2865)
	(2504,2903)(2505,2939)(2506,2970)
	(2507,2996)(2508,3016)(2508,3029)
	(2508,3037)(2508,3040)
\path(1900,3040)(1900,3037)(1900,3029)
	(1900,3016)(1901,2996)(1902,2970)
	(1903,2939)(1904,2903)(1905,2865)
	(1907,2826)(1909,2787)(1911,2750)
	(1913,2714)(1916,2681)(1919,2650)
	(1922,2622)(1925,2596)(1929,2572)
	(1934,2550)(1939,2529)(1944,2509)
	(1950,2490)(1957,2470)(1965,2450)
	(1974,2431)(1983,2412)(1994,2393)
	(2005,2375)(2017,2358)(2030,2341)
	(2043,2325)(2057,2311)(2071,2297)
	(2086,2285)(2100,2274)(2115,2265)
	(2130,2257)(2144,2251)(2158,2246)
	(2172,2243)(2186,2241)(2200,2240)
	(2214,2241)(2228,2243)(2242,2246)
	(2256,2251)(2270,2257)(2285,2265)
	(2300,2274)(2314,2285)(2329,2297)
	(2343,2311)(2357,2325)(2370,2341)
	(2383,2358)(2395,2375)(2406,2393)
	(2417,2412)(2426,2431)(2435,2450)
	(2443,2470)(2450,2490)(2456,2509)
	(2461,2529)(2466,2550)(2471,2572)
	(2475,2596)(2478,2622)(2481,2650)
	(2484,2681)(2487,2714)(2489,2750)
	(2491,2787)(2493,2826)(2495,2865)
	(2496,2903)(2497,2939)(2498,2970)
	(2499,2996)(2500,3016)(2500,3029)
	(2500,3037)(2500,3040)
\path(1945,33)(1945,36)(1945,44)
	(1945,57)(1946,77)(1947,103)
	(1948,134)(1949,170)(1950,208)
	(1952,247)(1954,286)(1956,323)
	(1958,359)(1961,392)(1964,423)
	(1967,451)(1970,477)(1974,501)
	(1979,523)(1984,544)(1989,564)
	(1995,583)(2002,603)(2010,623)
	(2019,642)(2028,661)(2039,680)
	(2050,698)(2062,715)(2075,732)
	(2088,748)(2102,762)(2116,776)
	(2131,788)(2145,799)(2160,808)
	(2175,816)(2189,822)(2203,827)
	(2217,830)(2231,832)(2245,833)
	(2259,832)(2273,830)(2287,827)
	(2301,822)(2315,816)(2330,808)
	(2345,799)(2359,788)(2374,776)
	(2388,762)(2402,748)(2415,732)
	(2428,715)(2440,698)(2451,680)
	(2462,661)(2471,642)(2480,623)
	(2488,603)(2495,583)(2501,564)
	(2506,544)(2511,523)(2516,501)
	(2520,477)(2523,451)(2526,423)
	(2529,392)(2532,359)(2534,323)
	(2536,286)(2538,247)(2540,208)
	(2541,170)(2542,134)(2543,103)
	(2544,77)(2545,57)(2545,44)
	(2545,36)(2545,33)
\path(1833,40)(1833,43)(1833,51)
	(1833,64)(1834,84)(1835,110)
	(1836,141)(1837,177)(1838,215)
	(1840,254)(1842,293)(1844,330)
	(1846,366)(1849,399)(1852,430)
	(1855,458)(1858,484)(1862,508)
	(1867,530)(1872,551)(1877,571)
	(1883,590)(1890,610)(1898,630)
	(1907,649)(1916,668)(1927,687)
	(1938,705)(1950,722)(1963,739)
	(1976,755)(1990,769)(2004,783)
	(2019,795)(2033,806)(2048,815)
	(2063,823)(2077,829)(2091,834)
	(2105,837)(2119,839)(2133,840)
	(2147,839)(2161,837)(2175,834)
	(2189,829)(2203,823)(2218,815)
	(2233,806)(2247,795)(2262,783)
	(2276,769)(2290,755)(2303,739)
	(2316,722)(2328,705)(2339,687)
	(2350,668)(2359,649)(2368,630)
	(2376,610)(2383,590)(2389,571)
	(2394,551)(2399,530)(2404,508)
	(2408,484)(2411,458)(2414,430)
	(2417,399)(2420,366)(2422,330)
	(2424,293)(2426,254)(2428,215)
	(2429,177)(2430,141)(2431,110)
	(2432,84)(2433,64)(2433,51)
	(2433,43)(2433,40)
\path(1870,40)(1870,43)(1870,51)
	(1870,64)(1871,84)(1872,110)
	(1873,141)(1874,177)(1875,215)
	(1877,254)(1879,293)(1881,330)
	(1883,366)(1886,399)(1889,430)
	(1892,458)(1895,484)(1899,508)
	(1904,530)(1909,551)(1914,571)
	(1920,590)(1927,610)(1935,630)
	(1944,649)(1953,668)(1964,687)
	(1975,705)(1987,722)(2000,739)
	(2013,755)(2027,769)(2041,783)
	(2056,795)(2070,806)(2085,815)
	(2100,823)(2114,829)(2128,834)
	(2142,837)(2156,839)(2170,840)
	(2184,839)(2198,837)(2212,834)
	(2226,829)(2240,823)(2255,815)
	(2270,806)(2284,795)(2299,783)
	(2313,769)(2327,755)(2340,739)
	(2353,722)(2365,705)(2376,687)
	(2387,668)(2396,649)(2405,630)
	(2413,610)(2420,590)(2426,571)
	(2431,551)(2436,530)(2441,508)
	(2445,484)(2448,458)(2451,430)
	(2454,399)(2457,366)(2459,330)
	(2461,293)(2463,254)(2465,215)
	(2466,177)(2467,141)(2468,110)
	(2469,84)(2470,64)(2470,51)
	(2470,43)(2470,40)
\path(1908,40)(1908,43)(1908,51)
	(1908,64)(1909,84)(1910,110)
	(1911,141)(1912,177)(1913,215)
	(1915,254)(1917,293)(1919,330)
	(1921,366)(1924,399)(1927,430)
	(1930,458)(1933,484)(1937,508)
	(1942,530)(1947,551)(1952,571)
	(1958,590)(1965,610)(1973,630)
	(1982,649)(1991,668)(2002,687)
	(2013,705)(2025,722)(2038,739)
	(2051,755)(2065,769)(2079,783)
	(2094,795)(2108,806)(2123,815)
	(2138,823)(2152,829)(2166,834)
	(2180,837)(2194,839)(2208,840)
	(2222,839)(2236,837)(2250,834)
	(2264,829)(2278,823)(2293,815)
	(2308,806)(2322,795)(2337,783)
	(2351,769)(2365,755)(2378,739)
	(2391,722)(2403,705)(2414,687)
	(2425,668)(2434,649)(2443,630)
	(2451,610)(2458,590)(2464,571)
	(2469,551)(2474,530)(2479,508)
	(2483,484)(2486,458)(2489,430)
	(2492,399)(2495,366)(2497,330)
	(2499,293)(2501,254)(2503,215)
	(2504,177)(2505,141)(2506,110)
	(2507,84)(2508,64)(2508,51)
	(2508,43)(2508,40)
\path(1900,40)(1900,43)(1900,51)
	(1900,64)(1901,84)(1902,110)
	(1903,141)(1904,177)(1905,215)
	(1907,254)(1909,293)(1911,330)
	(1913,366)(1916,399)(1919,430)
	(1922,458)(1925,484)(1929,508)
	(1934,530)(1939,551)(1944,571)
	(1950,590)(1957,610)(1965,630)
	(1974,649)(1983,668)(1994,687)
	(2005,705)(2017,722)(2030,739)
	(2043,755)(2057,769)(2071,783)
	(2086,795)(2100,806)(2115,815)
	(2130,823)(2144,829)(2158,834)
	(2172,837)(2186,839)(2200,840)
	(2214,839)(2228,837)(2242,834)
	(2256,829)(2270,823)(2285,815)
	(2300,806)(2314,795)(2329,783)
	(2343,769)(2357,755)(2370,739)
	(2383,722)(2395,705)(2406,687)
	(2417,668)(2426,649)(2435,630)
	(2443,610)(2450,590)(2456,571)
	(2461,551)(2466,530)(2471,508)
	(2475,484)(2478,458)(2481,430)
	(2484,399)(2487,366)(2489,330)
	(2491,293)(2493,254)(2495,215)
	(2496,177)(2497,141)(2498,110)
	(2499,84)(2500,64)(2500,51)
	(2500,43)(2500,40)
\path(183,6040)(183,6037)(183,6032)
	(184,6021)(184,6006)(185,5984)
	(187,5956)(189,5922)(191,5884)
	(193,5841)(197,5796)(200,5750)
	(204,5703)(209,5656)(214,5611)
	(219,5567)(225,5526)(231,5486)
	(238,5449)(246,5414)(254,5381)
	(262,5350)(272,5321)(282,5293)
	(293,5266)(306,5240)(319,5215)
	(333,5190)(347,5167)(362,5145)
	(379,5122)(396,5100)(414,5078)
	(433,5056)(454,5034)(476,5013)
	(498,4992)(522,4971)(547,4950)
	(572,4930)(599,4911)(626,4892)
	(654,4875)(682,4858)(711,4842)
	(740,4827)(769,4813)(799,4800)
	(828,4789)(857,4778)(886,4769)
	(915,4762)(944,4755)(972,4750)
	(1000,4745)(1028,4742)(1055,4741)
	(1083,4740)(1111,4741)(1138,4742)
	(1166,4745)(1194,4750)(1222,4755)
	(1251,4762)(1280,4769)(1309,4778)
	(1338,4789)(1367,4800)(1397,4813)
	(1426,4827)(1455,4842)(1484,4858)
	(1512,4875)(1540,4892)(1567,4911)
	(1594,4930)(1619,4950)(1644,4971)
	(1668,4992)(1690,5013)(1712,5034)
	(1733,5056)(1752,5078)(1770,5100)
	(1787,5122)(1804,5145)(1819,5167)
	(1833,5190)(1847,5215)(1860,5240)
	(1873,5266)(1884,5293)(1894,5321)
	(1904,5350)(1912,5381)(1920,5414)
	(1928,5449)(1935,5486)(1941,5526)
	(1947,5567)(1952,5611)(1957,5656)
	(1962,5703)(1966,5750)(1969,5796)
	(1973,5841)(1975,5884)(1977,5922)
	(1979,5956)(1981,5984)(1982,6006)
	(1982,6021)(1983,6032)(1983,6037)(1983,6040)
\path(175,6033)(175,6030)(175,6025)
	(176,6014)(176,5999)(177,5977)
	(179,5949)(181,5915)(183,5877)
	(185,5834)(189,5789)(192,5743)
	(196,5696)(201,5649)(206,5604)
	(211,5560)(217,5519)(223,5479)
	(230,5442)(238,5407)(246,5374)
	(254,5343)(264,5314)(274,5286)
	(285,5259)(298,5233)(311,5208)
	(325,5183)(339,5160)(354,5138)
	(371,5115)(388,5093)(406,5071)
	(425,5049)(446,5027)(468,5006)
	(490,4985)(514,4964)(539,4943)
	(564,4923)(591,4904)(618,4885)
	(646,4868)(674,4851)(703,4835)
	(732,4820)(761,4806)(791,4793)
	(820,4782)(849,4771)(878,4762)
	(907,4755)(936,4748)(964,4743)
	(992,4738)(1020,4735)(1047,4734)
	(1075,4733)(1103,4734)(1130,4735)
	(1158,4738)(1186,4743)(1214,4748)
	(1243,4755)(1272,4762)(1301,4771)
	(1330,4782)(1359,4793)(1389,4806)
	(1418,4820)(1447,4835)(1476,4851)
	(1504,4868)(1532,4885)(1559,4904)
	(1586,4923)(1611,4943)(1636,4964)
	(1660,4985)(1682,5006)(1704,5027)
	(1725,5049)(1744,5071)(1762,5093)
	(1779,5115)(1796,5138)(1811,5160)
	(1825,5183)(1839,5208)(1852,5233)
	(1865,5259)(1876,5286)(1886,5314)
	(1896,5343)(1904,5374)(1912,5407)
	(1920,5442)(1927,5479)(1933,5519)
	(1939,5560)(1944,5604)(1949,5649)
	(1954,5696)(1958,5743)(1961,5789)
	(1965,5834)(1967,5877)(1969,5915)
	(1971,5949)(1973,5977)(1974,5999)
	(1974,6014)(1975,6025)(1975,6030)(1975,6033)
\path(152,6033)(152,6030)(152,6025)
	(153,6014)(153,5999)(154,5977)
	(156,5949)(158,5915)(160,5877)
	(162,5834)(166,5789)(169,5743)
	(173,5696)(178,5649)(183,5604)
	(188,5560)(194,5519)(200,5479)
	(207,5442)(215,5407)(223,5374)
	(231,5343)(241,5314)(251,5286)
	(262,5259)(275,5233)(288,5208)
	(302,5183)(316,5160)(331,5138)
	(348,5115)(365,5093)(383,5071)
	(402,5049)(423,5027)(445,5006)
	(467,4985)(491,4964)(516,4943)
	(541,4923)(568,4904)(595,4885)
	(623,4868)(651,4851)(680,4835)
	(709,4820)(738,4806)(768,4793)
	(797,4782)(826,4771)(855,4762)
	(884,4755)(913,4748)(941,4743)
	(969,4738)(997,4735)(1024,4734)
	(1052,4733)(1080,4734)(1107,4735)
	(1135,4738)(1163,4743)(1191,4748)
	(1220,4755)(1249,4762)(1278,4771)
	(1307,4782)(1336,4793)(1366,4806)
	(1395,4820)(1424,4835)(1453,4851)
	(1481,4868)(1509,4885)(1536,4904)
	(1563,4923)(1588,4943)(1613,4964)
	(1637,4985)(1659,5006)(1681,5027)
	(1702,5049)(1721,5071)(1739,5093)
	(1756,5115)(1773,5138)(1788,5160)
	(1802,5183)(1816,5208)(1829,5233)
	(1842,5259)(1853,5286)(1863,5314)
	(1873,5343)(1881,5374)(1889,5407)
	(1897,5442)(1904,5479)(1910,5519)
	(1916,5560)(1921,5604)(1926,5649)
	(1931,5696)(1935,5743)(1938,5789)
	(1942,5834)(1944,5877)(1946,5915)
	(1948,5949)(1950,5977)(1951,5999)
	(1951,6014)(1952,6025)(1952,6030)(1952,6033)
\path(175,3041)(175,3044)(175,3049)
	(176,3060)(176,3075)(177,3097)
	(179,3125)(181,3159)(183,3197)
	(185,3240)(189,3285)(192,3331)
	(196,3378)(201,3425)(206,3470)
	(211,3514)(217,3555)(223,3595)
	(230,3632)(238,3667)(246,3700)
	(254,3731)(264,3760)(274,3788)
	(285,3815)(298,3841)(311,3866)
	(325,3891)(339,3914)(354,3936)
	(371,3959)(388,3981)(406,4003)
	(425,4025)(446,4047)(468,4068)
	(490,4089)(514,4110)(539,4131)
	(564,4151)(591,4170)(618,4189)
	(646,4206)(674,4223)(703,4239)
	(732,4254)(761,4268)(791,4281)
	(820,4292)(849,4303)(878,4312)
	(907,4319)(936,4326)(964,4331)
	(992,4336)(1020,4339)(1047,4340)
	(1075,4341)(1103,4340)(1130,4339)
	(1158,4336)(1186,4331)(1214,4326)
	(1243,4319)(1272,4312)(1301,4303)
	(1330,4292)(1359,4281)(1389,4268)
	(1418,4254)(1447,4239)(1476,4223)
	(1504,4206)(1532,4189)(1559,4170)
	(1586,4151)(1611,4131)(1636,4110)
	(1660,4089)(1682,4068)(1704,4047)
	(1725,4025)(1744,4003)(1762,3981)
	(1779,3959)(1796,3936)(1811,3914)
	(1825,3891)(1839,3866)(1852,3841)
	(1865,3815)(1876,3788)(1886,3760)
	(1896,3731)(1904,3700)(1912,3667)
	(1920,3632)(1927,3595)(1933,3555)
	(1939,3514)(1944,3470)(1949,3425)
	(1954,3378)(1958,3331)(1961,3285)
	(1965,3240)(1967,3197)(1969,3159)
	(1971,3125)(1973,3097)(1974,3075)
	(1974,3060)(1975,3049)(1975,3044)(1975,3041)
\path(152,3041)(152,3044)(152,3049)
	(153,3060)(153,3075)(154,3097)
	(156,3125)(158,3159)(160,3197)
	(162,3240)(166,3285)(169,3331)
	(173,3378)(178,3425)(183,3470)
	(188,3514)(194,3555)(200,3595)
	(207,3632)(215,3667)(223,3700)
	(231,3731)(241,3760)(251,3788)
	(262,3815)(275,3841)(288,3866)
	(302,3891)(316,3914)(331,3936)
	(348,3959)(365,3981)(383,4003)
	(402,4025)(423,4047)(445,4068)
	(467,4089)(491,4110)(516,4131)
	(541,4151)(568,4170)(595,4189)
	(623,4206)(651,4223)(680,4239)
	(709,4254)(738,4268)(768,4281)
	(797,4292)(826,4303)(855,4312)
	(884,4319)(913,4326)(941,4331)
	(969,4336)(997,4339)(1024,4340)
	(1052,4341)(1080,4340)(1107,4339)
	(1135,4336)(1163,4331)(1191,4326)
	(1220,4319)(1249,4312)(1278,4303)
	(1307,4292)(1336,4281)(1366,4268)
	(1395,4254)(1424,4239)(1453,4223)
	(1481,4206)(1509,4189)(1536,4170)
	(1563,4151)(1588,4131)(1613,4110)
	(1637,4089)(1659,4068)(1681,4047)
	(1702,4025)(1721,4003)(1739,3981)
	(1756,3959)(1773,3936)(1788,3914)
	(1802,3891)(1816,3866)(1829,3841)
	(1842,3815)(1853,3788)(1863,3760)
	(1873,3731)(1881,3700)(1889,3667)
	(1897,3632)(1904,3595)(1910,3555)
	(1916,3514)(1921,3470)(1926,3425)
	(1931,3378)(1935,3331)(1938,3285)
	(1942,3240)(1944,3197)(1946,3159)
	(1948,3125)(1950,3097)(1951,3075)
	(1951,3060)(1952,3049)(1952,3044)(1952,3041)
\path(1975,3056)(1975,3053)(1975,3045)
	(1975,3032)(1976,3012)(1977,2986)
	(1978,2955)(1979,2919)(1980,2881)
	(1982,2842)(1984,2803)(1986,2766)
	(1988,2730)(1991,2697)(1994,2666)
	(1997,2638)(2000,2612)(2004,2588)
	(2009,2566)(2014,2545)(2019,2525)
	(2025,2506)(2032,2486)(2040,2466)
	(2049,2447)(2058,2428)(2069,2409)
	(2080,2391)(2092,2374)(2105,2357)
	(2118,2341)(2132,2327)(2146,2313)
	(2161,2301)(2175,2290)(2190,2281)
	(2205,2273)(2219,2267)(2233,2262)
	(2247,2259)(2261,2257)(2275,2256)
	(2289,2257)(2303,2259)(2317,2262)
	(2331,2267)(2345,2273)(2360,2281)
	(2375,2290)(2389,2301)(2404,2313)
	(2418,2327)(2432,2341)(2445,2357)
	(2458,2374)(2470,2391)(2481,2409)
	(2492,2428)(2501,2447)(2510,2466)
	(2518,2486)(2525,2506)(2531,2525)
	(2536,2545)(2541,2566)(2546,2588)
	(2550,2612)(2553,2638)(2556,2666)
	(2559,2697)(2562,2730)(2564,2766)
	(2566,2803)(2568,2842)(2570,2881)
	(2571,2919)(2572,2955)(2573,2986)
	(2574,3012)(2575,3032)(2575,3045)
	(2575,3053)(2575,3056)
\path(1990,33)(1990,36)(1990,44)
	(1990,57)(1991,77)(1992,103)
	(1993,134)(1994,170)(1995,208)
	(1997,247)(1999,286)(2001,323)
	(2003,359)(2006,392)(2009,423)
	(2012,451)(2015,477)(2019,501)
	(2024,523)(2029,544)(2034,564)
	(2040,583)(2047,603)(2055,623)
	(2064,642)(2073,661)(2084,680)
	(2095,698)(2107,715)(2120,732)
	(2133,748)(2147,762)(2161,776)
	(2176,788)(2190,799)(2205,808)
	(2220,816)(2234,822)(2248,827)
	(2262,830)(2276,832)(2290,833)
	(2304,832)(2318,830)(2332,827)
	(2346,822)(2360,816)(2375,808)
	(2390,799)(2404,788)(2419,776)
	(2433,762)(2447,748)(2460,732)
	(2473,715)(2485,698)(2496,680)
	(2507,661)(2516,642)(2525,623)
	(2533,603)(2540,583)(2546,564)
	(2551,544)(2556,523)(2561,501)
	(2565,477)(2568,451)(2571,423)
	(2574,392)(2577,359)(2579,323)
	(2581,286)(2583,247)(2585,208)
	(2586,170)(2587,134)(2588,103)
	(2589,77)(2590,57)(2590,44)
	(2590,36)(2590,33)
\end{picture}
}

%% file: xfig/pascal11.eepic
\setlength{\unitlength}{0.00037500in}
\begingroup\makeatletter\ifx\SetFigFont\undefined%
\gdef\SetFigFont#1#2#3#4#5{%
  \reset@font\fontsize{#1}{#2pt}%
  \fontfamily{#3}\fontseries{#4}\fontshape{#5}%
  \selectfont}%
\fi\endgroup%
{\renewcommand{\dashlinestretch}{30}
\begin{picture}(4292,5639)(0,-10)
\dashline{60.000}(1875,5612)(1875,212)
\path(1275,512)(1278,511)(1284,508)
	(1294,503)(1310,496)(1332,486)
	(1359,474)(1390,460)(1425,445)
	(1462,429)(1499,414)(1536,398)
	(1573,384)(1608,371)(1641,359)
	(1673,349)(1702,339)(1730,332)
	(1756,325)(1781,320)(1805,317)
	(1829,314)(1852,312)(1875,312)
	(1898,312)(1921,314)(1945,317)
	(1969,320)(1994,325)(2020,332)
	(2048,339)(2077,349)(2109,359)
	(2142,371)(2177,384)(2214,398)
	(2251,414)(2288,429)(2325,445)
	(2360,460)(2391,474)(2418,486)
	(2440,496)(2475,512)
\path(2378.336,434.825)(2475.000,512.000)(2353.390,489.393)
\path(675,512)(677,511)(681,509)
	(689,506)(702,501)(720,494)
	(743,484)(772,473)(807,459)
	(846,444)(890,427)(937,409)
	(988,390)(1039,371)(1092,351)
	(1145,332)(1198,313)(1250,295)
	(1300,278)(1349,263)(1396,248)
	(1441,234)(1484,222)(1526,211)
	(1565,201)(1604,192)(1640,185)
	(1676,179)(1711,173)(1744,169)
	(1777,166)(1810,164)(1843,162)
	(1875,162)(1907,162)(1940,164)
	(1973,166)(2006,169)(2039,173)
	(2074,179)(2110,185)(2146,192)
	(2185,201)(2224,211)(2266,222)
	(2309,234)(2354,248)(2401,263)
	(2450,278)(2500,295)(2552,313)
	(2605,332)(2658,351)(2711,371)
	(2762,390)(2813,409)(2860,427)
	(2904,444)(2943,459)(2978,473)
	(3007,484)(3030,494)(3048,501)(3075,512)
\path(2975.188,438.942)(3075.000,512.000)(2952.550,494.507)
\path(75,512)(77,511)(80,510)
	(87,507)(97,503)(112,498)
	(133,490)(158,480)(189,468)
	(226,455)(268,439)(314,421)
	(366,402)(421,382)(479,361)
	(540,339)(603,317)(667,294)
	(731,272)(796,250)(859,228)
	(922,207)(984,187)(1044,168)
	(1103,151)(1160,134)(1215,118)
	(1268,103)(1319,90)(1368,78)
	(1416,67)(1463,57)(1508,48)
	(1551,40)(1594,33)(1636,27)
	(1677,22)(1717,19)(1757,16)
	(1796,14)(1836,12)(1875,12)
	(1914,12)(1954,14)(1993,16)
	(2033,19)(2073,22)(2114,27)
	(2156,33)(2199,40)(2242,48)
	(2287,57)(2334,67)(2382,78)
	(2431,90)(2482,103)(2535,118)
	(2590,134)(2647,151)(2706,168)
	(2766,187)(2828,207)(2891,228)
	(2954,250)(3019,272)(3083,294)
	(3147,317)(3210,339)(3271,361)
	(3329,382)(3384,402)(3436,421)
	(3482,439)(3524,455)(3561,468)
	(3592,480)(3617,490)(3638,498)
	(3653,503)(3675,512)
\path(3575.293,438.798)(3675.000,512.000)(3552.575,494.330)
\put(2100,5012){\makebox(0,0)[lb]{\smash{{{\SetFigFont{5}{6.0}{\rmdefault}{\mddefault}{\updefault}1}}}}}
\put(1800,4412){\makebox(0,0)[lb]{\smash{{{\SetFigFont{5}{6.0}{\rmdefault}{\mddefault}{\updefault}1}}}}}
\put(2400,4412){\makebox(0,0)[lb]{\smash{{{\SetFigFont{5}{6.0}{\rmdefault}{\mddefault}{\updefault}1}}}}}
\put(2100,3812){\makebox(0,0)[lb]{\smash{{{\SetFigFont{5}{6.0}{\rmdefault}{\mddefault}{\updefault}1}}}}}
\put(3000,3212){\makebox(0,0)[lb]{\smash{{{\SetFigFont{5}{6.0}{\rmdefault}{\mddefault}{\updefault}1}}}}}
\put(1800,3212){\makebox(0,0)[lb]{\smash{{{\SetFigFont{5}{6.0}{\rmdefault}{\mddefault}{\updefault}3}}}}}
\put(1500,3812){\makebox(0,0)[lb]{\smash{{{\SetFigFont{5}{6.0}{\rmdefault}{\mddefault}{\updefault}1}}}}}
\put(2700,3812){\makebox(0,0)[lb]{\smash{{{\SetFigFont{5}{6.0}{\rmdefault}{\mddefault}{\updefault}1}}}}}
\put(1200,3212){\makebox(0,0)[lb]{\smash{{{\SetFigFont{5}{6.0}{\rmdefault}{\mddefault}{\updefault}1}}}}}
\put(900,2612){\makebox(0,0)[lb]{\smash{{{\SetFigFont{5}{6.0}{\rmdefault}{\mddefault}{\updefault}1}}}}}
\put(1500,2612){\makebox(0,0)[lb]{\smash{{{\SetFigFont{5}{6.0}{\rmdefault}{\mddefault}{\updefault}4}}}}}
\put(2100,2612){\makebox(0,0)[lb]{\smash{{{\SetFigFont{5}{6.0}{\rmdefault}{\mddefault}{\updefault}2}}}}}
\put(2700,2612){\makebox(0,0)[lb]{\smash{{{\SetFigFont{5}{6.0}{\rmdefault}{\mddefault}{\updefault}3}}}}}
\put(3300,2612){\makebox(0,0)[lb]{\smash{{{\SetFigFont{5}{6.0}{\rmdefault}{\mddefault}{\updefault}1}}}}}
\put(600,2012){\makebox(0,0)[lb]{\smash{{{\SetFigFont{5}{6.0}{\rmdefault}{\mddefault}{\updefault}1}}}}}
\put(1200,2012){\makebox(0,0)[lb]{\smash{{{\SetFigFont{5}{6.0}{\rmdefault}{\mddefault}{\updefault}5}}}}}
\put(1800,2012){\makebox(0,0)[lb]{\smash{{{\SetFigFont{5}{6.0}{\rmdefault}{\mddefault}{\updefault}10}}}}}
\put(3000,2012){\makebox(0,0)[lb]{\smash{{{\SetFigFont{5}{6.0}{\rmdefault}{\mddefault}{\updefault}4}}}}}
\put(3600,2012){\makebox(0,0)[lb]{\smash{{{\SetFigFont{5}{6.0}{\rmdefault}{\mddefault}{\updefault}1}}}}}
\put(300,1412){\makebox(0,0)[lb]{\smash{{{\SetFigFont{5}{6.0}{\rmdefault}{\mddefault}{\updefault}1}}}}}
\put(900,1412){\makebox(0,0)[lb]{\smash{{{\SetFigFont{5}{6.0}{\rmdefault}{\mddefault}{\updefault}6}}}}}
\put(1500,1412){\makebox(0,0)[lb]{\smash{{{\SetFigFont{5}{6.0}{\rmdefault}{\mddefault}{\updefault}15}}}}}
\put(2100,1412){\makebox(0,0)[lb]{\smash{{{\SetFigFont{5}{6.0}{\rmdefault}{\mddefault}{\updefault}5}}}}}
\put(2700,1412){\makebox(0,0)[lb]{\smash{{{\SetFigFont{5}{6.0}{\rmdefault}{\mddefault}{\updefault}9}}}}}
\put(3300,1412){\makebox(0,0)[lb]{\smash{{{\SetFigFont{5}{6.0}{\rmdefault}{\mddefault}{\updefault}5}}}}}
\put(3900,1412){\makebox(0,0)[lb]{\smash{{{\SetFigFont{5}{6.0}{\rmdefault}{\mddefault}{\updefault}1}}}}}
\put(0,812){\makebox(0,0)[lb]{\smash{{{\SetFigFont{5}{6.0}{\rmdefault}{\mddefault}{\updefault}1}}}}}
\put(600,812){\makebox(0,0)[lb]{\smash{{{\SetFigFont{5}{6.0}{\rmdefault}{\mddefault}{\updefault}7}}}}}
\put(1200,812){\makebox(0,0)[lb]{\smash{{{\SetFigFont{5}{6.0}{\rmdefault}{\mddefault}{\updefault}21}}}}}
\put(1800,812){\makebox(0,0)[lb]{\smash{{{\SetFigFont{5}{6.0}{\rmdefault}{\mddefault}{\updefault}35}}}}}
\put(2400,812){\makebox(0,0)[lb]{\smash{{{\SetFigFont{5}{6.0}{\rmdefault}{\mddefault}{\updefault}14}}}}}
\put(3000,812){\makebox(0,0)[lb]{\smash{{{\SetFigFont{5}{6.0}{\rmdefault}{\mddefault}{\updefault}14}}}}}
\put(3600,812){\makebox(0,0)[lb]{\smash{{{\SetFigFont{5}{6.0}{\rmdefault}{\mddefault}{\updefault}6}}}}}
\put(4200,812){\makebox(0,0)[lb]{\smash{{{\SetFigFont{5}{6.0}{\rmdefault}{\mddefault}{\updefault}1}}}}}
\put(2100,3587){\makebox(0,0)[lb]{\smash{{{\SetFigFont{5}{6.0}{\rmdefault}{\mddefault}{\updefault}1}}}}}
\put(2400,2987){\makebox(0,0)[lb]{\smash{{{\SetFigFont{5}{6.0}{\rmdefault}{\mddefault}{\updefault}1}}}}}
\put(2100,2387){\makebox(0,0)[lb]{\smash{{{\SetFigFont{5}{6.0}{\rmdefault}{\mddefault}{\updefault}4}}}}}
\put(2700,2387){\makebox(0,0)[lb]{\smash{{{\SetFigFont{5}{6.0}{\rmdefault}{\mddefault}{\updefault}1}}}}}
\put(2400,1787){\makebox(0,0)[lb]{\smash{{{\SetFigFont{5}{6.0}{\rmdefault}{\mddefault}{\updefault}5}}}}}
\put(3000,1787){\makebox(0,0)[lb]{\smash{{{\SetFigFont{5}{6.0}{\rmdefault}{\mddefault}{\updefault}1}}}}}
\put(2100,1187){\makebox(0,0)[lb]{\smash{{{\SetFigFont{5}{6.0}{\rmdefault}{\mddefault}{\updefault}15}}}}}
\put(2700,1187){\makebox(0,0)[lb]{\smash{{{\SetFigFont{5}{6.0}{\rmdefault}{\mddefault}{\updefault}6}}}}}
\put(3300,1187){\makebox(0,0)[lb]{\smash{{{\SetFigFont{5}{6.0}{\rmdefault}{\mddefault}{\updefault}1}}}}}
\put(2400,3212){\makebox(0,0)[lb]{\smash{{{\SetFigFont{5}{6.0}{\rmdefault}{\mddefault}{\updefault}2}}}}}
\put(2400,2012){\makebox(0,0)[lb]{\smash{{{\SetFigFont{5}{6.0}{\rmdefault}{\mddefault}{\updefault}5}}}}}
\put(2400,587){\makebox(0,0)[lb]{\smash{{{\SetFigFont{5}{6.0}{\rmdefault}{\mddefault}{\updefault}21}}}}}
\put(3000,587){\makebox(0,0)[lb]{\smash{{{\SetFigFont{5}{6.0}{\rmdefault}{\mddefault}{\updefault}7}}}}}
\put(3600,587){\makebox(0,0)[lb]{\smash{{{\SetFigFont{5}{6.0}{\rmdefault}{\mddefault}{\updefault}1}}}}}
\end{picture}
}

%% file: xfig/pascal2.eepic
\setlength{\unitlength}{0.00037500in}
\begingroup\makeatletter\ifx\SetFigFont\undefined%
\gdef\SetFigFont#1#2#3#4#5{%
  \reset@font\fontsize{#1}{#2pt}%
  \fontfamily{#3}\fontseries{#4}\fontshape{#5}%
  \selectfont}%
\fi\endgroup%
{\renewcommand{\dashlinestretch}{30}
\begin{picture}(5438,5739)(0,-10)
\dashline{180.000}(12,4587)(4812,4587)
\dashline{60.000}(1962,5712)(1962,12)
\path(1662,612)(1663,610)(1667,605)
	(1672,597)(1680,586)(1691,571)
	(1705,553)(1720,534)(1737,514)
	(1755,493)(1775,473)(1796,454)
	(1818,437)(1842,421)(1869,407)
	(1898,397)(1929,390)(1962,387)
	(1995,390)(2026,397)(2055,407)
	(2082,421)(2106,437)(2128,454)
	(2149,473)(2169,493)(2187,514)
	(2204,534)(2219,553)(2233,571)
	(2244,586)(2262,612)
\path(2218.361,496.261)(2262.000,612.000)(2169.029,530.413)
\path(1062,612)(1063,611)(1065,609)
	(1069,606)(1076,601)(1085,594)
	(1097,584)(1112,573)(1130,559)
	(1150,543)(1173,526)(1199,508)
	(1226,488)(1256,468)(1287,448)
	(1320,427)(1354,407)(1389,387)
	(1425,368)(1464,349)(1503,331)
	(1545,314)(1588,299)(1634,284)
	(1683,271)(1734,260)(1787,250)
	(1844,243)(1902,239)(1962,237)
	(2022,239)(2080,243)(2137,250)
	(2190,260)(2241,271)(2290,284)
	(2336,299)(2379,314)(2421,331)
	(2460,349)(2499,368)(2535,387)
	(2570,407)(2604,427)(2637,448)
	(2668,468)(2698,488)(2725,508)
	(2751,526)(2774,543)(2794,559)
	(2812,573)(2827,584)(2839,594)
	(2848,601)(2862,612)
\path(2786.176,514.272)(2862.000,612.000)(2749.107,561.451)
\path(462,612)(463,611)(465,610)
	(468,608)(474,604)(482,599)
	(492,591)(506,582)(522,571)
	(542,559)(565,544)(590,528)
	(619,510)(650,491)(683,470)
	(719,449)(756,427)(796,405)
	(837,382)(880,360)(924,337)
	(969,315)(1016,293)(1065,271)
	(1114,251)(1166,231)(1219,211)
	(1274,193)(1332,176)(1391,159)
	(1453,144)(1518,130)(1586,118)
	(1657,108)(1730,99)(1806,92)
	(1883,88)(1962,87)(2041,88)
	(2118,92)(2194,99)(2267,108)
	(2338,118)(2406,130)(2471,144)
	(2533,159)(2592,176)(2650,193)
	(2705,211)(2758,231)(2810,251)
	(2859,271)(2908,293)(2955,315)
	(3000,337)(3044,360)(3087,382)
	(3128,405)(3168,427)(3205,449)
	(3241,470)(3274,491)(3305,510)
	(3334,528)(3359,544)(3382,559)
	(3402,571)(3418,582)(3432,591)
	(3442,599)(3450,604)(3462,612)
\path(3378.795,520.474)(3462.000,612.000)(3345.513,570.397)
\put(2487,5112){\makebox(0,0)[lb]{\smash{{{\SetFigFont{5}{6.0}{\rmdefault}{\mddefault}{\updefault}1}}}}}
\put(2187,4512){\makebox(0,0)[lb]{\smash{{{\SetFigFont{5}{6.0}{\rmdefault}{\mddefault}{\updefault}1}}}}}
\put(2787,4512){\makebox(0,0)[lb]{\smash{{{\SetFigFont{5}{6.0}{\rmdefault}{\mddefault}{\updefault}1}}}}}
\put(2487,3912){\makebox(0,0)[lb]{\smash{{{\SetFigFont{5}{6.0}{\rmdefault}{\mddefault}{\updefault}2}}}}}
\put(2787,3312){\makebox(0,0)[lb]{\smash{{{\SetFigFont{5}{6.0}{\rmdefault}{\mddefault}{\updefault}3}}}}}
\put(3387,3312){\makebox(0,0)[lb]{\smash{{{\SetFigFont{5}{6.0}{\rmdefault}{\mddefault}{\updefault}1}}}}}
\put(2187,3312){\makebox(0,0)[lb]{\smash{{{\SetFigFont{5}{6.0}{\rmdefault}{\mddefault}{\updefault}2}}}}}
\put(1887,3912){\makebox(0,0)[lb]{\smash{{{\SetFigFont{5}{6.0}{\rmdefault}{\mddefault}{\updefault}1}}}}}
\put(3087,3912){\makebox(0,0)[lb]{\smash{{{\SetFigFont{5}{6.0}{\rmdefault}{\mddefault}{\updefault}1}}}}}
\put(1587,3312){\makebox(0,0)[lb]{\smash{{{\SetFigFont{5}{6.0}{\rmdefault}{\mddefault}{\updefault}1}}}}}
\put(1287,2712){\makebox(0,0)[lb]{\smash{{{\SetFigFont{5}{6.0}{\rmdefault}{\mddefault}{\updefault}1}}}}}
\put(1887,2712){\makebox(0,0)[lb]{\smash{{{\SetFigFont{5}{6.0}{\rmdefault}{\mddefault}{\updefault}4}}}}}
\put(2487,2712){\makebox(0,0)[lb]{\smash{{{\SetFigFont{5}{6.0}{\rmdefault}{\mddefault}{\updefault}5}}}}}
\put(3087,2712){\makebox(0,0)[lb]{\smash{{{\SetFigFont{5}{6.0}{\rmdefault}{\mddefault}{\updefault}4}}}}}
\put(3687,2712){\makebox(0,0)[lb]{\smash{{{\SetFigFont{5}{6.0}{\rmdefault}{\mddefault}{\updefault}1}}}}}
\put(987,2112){\makebox(0,0)[lb]{\smash{{{\SetFigFont{5}{6.0}{\rmdefault}{\mddefault}{\updefault}1}}}}}
\put(1587,2112){\makebox(0,0)[lb]{\smash{{{\SetFigFont{5}{6.0}{\rmdefault}{\mddefault}{\updefault}5}}}}}
\put(2187,2112){\makebox(0,0)[lb]{\smash{{{\SetFigFont{5}{6.0}{\rmdefault}{\mddefault}{\updefault}5}}}}}
\put(2787,2112){\makebox(0,0)[lb]{\smash{{{\SetFigFont{5}{6.0}{\rmdefault}{\mddefault}{\updefault}9}}}}}
\put(3387,2112){\makebox(0,0)[lb]{\smash{{{\SetFigFont{5}{6.0}{\rmdefault}{\mddefault}{\updefault}5}}}}}
\put(3987,2112){\makebox(0,0)[lb]{\smash{{{\SetFigFont{5}{6.0}{\rmdefault}{\mddefault}{\updefault}1}}}}}
\put(687,1512){\makebox(0,0)[lb]{\smash{{{\SetFigFont{5}{6.0}{\rmdefault}{\mddefault}{\updefault}1}}}}}
\put(1287,1512){\makebox(0,0)[lb]{\smash{{{\SetFigFont{5}{6.0}{\rmdefault}{\mddefault}{\updefault}6}}}}}
\put(1887,1512){\makebox(0,0)[lb]{\smash{{{\SetFigFont{5}{6.0}{\rmdefault}{\mddefault}{\updefault}15}}}}}
\put(2487,1512){\makebox(0,0)[lb]{\smash{{{\SetFigFont{5}{6.0}{\rmdefault}{\mddefault}{\updefault}14}}}}}
\put(3087,1512){\makebox(0,0)[lb]{\smash{{{\SetFigFont{5}{6.0}{\rmdefault}{\mddefault}{\updefault}14}}}}}
\put(3687,1512){\makebox(0,0)[lb]{\smash{{{\SetFigFont{5}{6.0}{\rmdefault}{\mddefault}{\updefault}6}}}}}
\put(4287,1512){\makebox(0,0)[lb]{\smash{{{\SetFigFont{5}{6.0}{\rmdefault}{\mddefault}{\updefault}1}}}}}
\put(387,912){\makebox(0,0)[lb]{\smash{{{\SetFigFont{5}{6.0}{\rmdefault}{\mddefault}{\updefault}1}}}}}
\put(987,912){\makebox(0,0)[lb]{\smash{{{\SetFigFont{5}{6.0}{\rmdefault}{\mddefault}{\updefault}7}}}}}
\put(1587,912){\makebox(0,0)[lb]{\smash{{{\SetFigFont{5}{6.0}{\rmdefault}{\mddefault}{\updefault}21}}}}}
\put(2187,912){\makebox(0,0)[lb]{\smash{{{\SetFigFont{5}{6.0}{\rmdefault}{\mddefault}{\updefault}14}}}}}
\put(3387,912){\makebox(0,0)[lb]{\smash{{{\SetFigFont{5}{6.0}{\rmdefault}{\mddefault}{\updefault}20}}}}}
\put(3987,912){\makebox(0,0)[lb]{\smash{{{\SetFigFont{5}{6.0}{\rmdefault}{\mddefault}{\updefault}7}}}}}
\put(4587,912){\makebox(0,0)[lb]{\smash{{{\SetFigFont{5}{6.0}{\rmdefault}{\mddefault}{\updefault}1}}}}}
\put(2187,3087){\makebox(0,0)[lb]{\smash{{{\SetFigFont{5}{6.0}{\rmdefault}{\mddefault}{\updefault}1}}}}}
\put(2487,2487){\makebox(0,0)[lb]{\smash{{{\SetFigFont{5}{6.0}{\rmdefault}{\mddefault}{\updefault}1}}}}}
\put(2187,1887){\makebox(0,0)[lb]{\smash{{{\SetFigFont{5}{6.0}{\rmdefault}{\mddefault}{\updefault}5}}}}}
\put(2787,1887){\makebox(0,0)[lb]{\smash{{{\SetFigFont{5}{6.0}{\rmdefault}{\mddefault}{\updefault}1}}}}}
\put(2487,1287){\makebox(0,0)[lb]{\smash{{{\SetFigFont{5}{6.0}{\rmdefault}{\mddefault}{\updefault}6}}}}}
\put(3087,1287){\makebox(0,0)[lb]{\smash{{{\SetFigFont{5}{6.0}{\rmdefault}{\mddefault}{\updefault}1}}}}}
\put(2187,687){\makebox(0,0)[lb]{\smash{{{\SetFigFont{5}{6.0}{\rmdefault}{\mddefault}{\updefault}21}}}}}
\put(2787,687){\makebox(0,0)[lb]{\smash{{{\SetFigFont{5}{6.0}{\rmdefault}{\mddefault}{\updefault}7}}}}}
\put(3387,687){\makebox(0,0)[lb]{\smash{{{\SetFigFont{5}{6.0}{\rmdefault}{\mddefault}{\updefault}1}}}}}
\put(2787,912){\makebox(0,0)[lb]{\smash{{{\SetFigFont{5}{6.0}{\rmdefault}{\mddefault}{\updefault}28}}}}}
\put(4962,4512){\makebox(0,0)[lb]{\smash{{{\SetFigFont{5}{6.0}{\rmdefault}{\mddefault}{\updefault}$n=1$}}}}}
\end{picture}
}

%% file: xfig/pascal3.eepic
\setlength{\unitlength}{0.00037500in}
\begingroup\makeatletter\ifx\SetFigFont\undefined%
\gdef\SetFigFont#1#2#3#4#5{%
  \reset@font\fontsize{#1}{#2pt}%
  \fontfamily{#3}\fontseries{#4}\fontshape{#5}%
  \selectfont}%
\fi\endgroup%
{\renewcommand{\dashlinestretch}{30}
\begin{picture}(4292,5739)(0,-10)
\dashline{60.000}(1275,5712)(1275,12)
\put(2100,5112){\makebox(0,0)[lb]{\smash{{{\SetFigFont{5}{6.0}{\rmdefault}{\mddefault}{\updefault}1}}}}}
\put(1800,4512){\makebox(0,0)[lb]{\smash{{{\SetFigFont{5}{6.0}{\rmdefault}{\mddefault}{\updefault}1}}}}}
\put(2400,4512){\makebox(0,0)[lb]{\smash{{{\SetFigFont{5}{6.0}{\rmdefault}{\mddefault}{\updefault}1}}}}}
\put(2100,3912){\makebox(0,0)[lb]{\smash{{{\SetFigFont{5}{6.0}{\rmdefault}{\mddefault}{\updefault}2}}}}}
\put(2400,3312){\makebox(0,0)[lb]{\smash{{{\SetFigFont{5}{6.0}{\rmdefault}{\mddefault}{\updefault}3}}}}}
\put(3000,3312){\makebox(0,0)[lb]{\smash{{{\SetFigFont{5}{6.0}{\rmdefault}{\mddefault}{\updefault}1}}}}}
\put(1800,3312){\makebox(0,0)[lb]{\smash{{{\SetFigFont{5}{6.0}{\rmdefault}{\mddefault}{\updefault}3}}}}}
\put(1500,3912){\makebox(0,0)[lb]{\smash{{{\SetFigFont{5}{6.0}{\rmdefault}{\mddefault}{\updefault}1}}}}}
\put(2700,3912){\makebox(0,0)[lb]{\smash{{{\SetFigFont{5}{6.0}{\rmdefault}{\mddefault}{\updefault}1}}}}}
\put(1200,3312){\makebox(0,0)[lb]{\smash{{{\SetFigFont{5}{6.0}{\rmdefault}{\mddefault}{\updefault}1}}}}}
\put(900,2712){\makebox(0,0)[lb]{\smash{{{\SetFigFont{5}{6.0}{\rmdefault}{\mddefault}{\updefault}1}}}}}
\put(2100,2712){\makebox(0,0)[lb]{\smash{{{\SetFigFont{5}{6.0}{\rmdefault}{\mddefault}{\updefault}6}}}}}
\put(2700,2712){\makebox(0,0)[lb]{\smash{{{\SetFigFont{5}{6.0}{\rmdefault}{\mddefault}{\updefault}4}}}}}
\put(3300,2712){\makebox(0,0)[lb]{\smash{{{\SetFigFont{5}{6.0}{\rmdefault}{\mddefault}{\updefault}1}}}}}
\put(600,2112){\makebox(0,0)[lb]{\smash{{{\SetFigFont{5}{6.0}{\rmdefault}{\mddefault}{\updefault}1}}}}}
\put(1200,2112){\makebox(0,0)[lb]{\smash{{{\SetFigFont{5}{6.0}{\rmdefault}{\mddefault}{\updefault}5}}}}}
\put(1800,2112){\makebox(0,0)[lb]{\smash{{{\SetFigFont{5}{6.0}{\rmdefault}{\mddefault}{\updefault}9}}}}}
\put(2400,2112){\makebox(0,0)[lb]{\smash{{{\SetFigFont{5}{6.0}{\rmdefault}{\mddefault}{\updefault}10}}}}}
\put(3000,2112){\makebox(0,0)[lb]{\smash{{{\SetFigFont{5}{6.0}{\rmdefault}{\mddefault}{\updefault}5}}}}}
\put(3600,2112){\makebox(0,0)[lb]{\smash{{{\SetFigFont{5}{6.0}{\rmdefault}{\mddefault}{\updefault}1}}}}}
\put(300,1512){\makebox(0,0)[lb]{\smash{{{\SetFigFont{5}{6.0}{\rmdefault}{\mddefault}{\updefault}1}}}}}
\put(900,1512){\makebox(0,0)[lb]{\smash{{{\SetFigFont{5}{6.0}{\rmdefault}{\mddefault}{\updefault}6}}}}}
\put(1500,1512){\makebox(0,0)[lb]{\smash{{{\SetFigFont{5}{6.0}{\rmdefault}{\mddefault}{\updefault}9}}}}}
\put(2100,1512){\makebox(0,0)[lb]{\smash{{{\SetFigFont{5}{6.0}{\rmdefault}{\mddefault}{\updefault}19}}}}}
\put(2700,1512){\makebox(0,0)[lb]{\smash{{{\SetFigFont{5}{6.0}{\rmdefault}{\mddefault}{\updefault}15}}}}}
\put(3300,1512){\makebox(0,0)[lb]{\smash{{{\SetFigFont{5}{6.0}{\rmdefault}{\mddefault}{\updefault}6}}}}}
\put(3900,1512){\makebox(0,0)[lb]{\smash{{{\SetFigFont{5}{6.0}{\rmdefault}{\mddefault}{\updefault}1}}}}}
\put(0,912){\makebox(0,0)[lb]{\smash{{{\SetFigFont{5}{6.0}{\rmdefault}{\mddefault}{\updefault}1}}}}}
\put(600,912){\makebox(0,0)[lb]{\smash{{{\SetFigFont{5}{6.0}{\rmdefault}{\mddefault}{\updefault}7}}}}}
\put(1200,912){\makebox(0,0)[lb]{\smash{{{\SetFigFont{5}{6.0}{\rmdefault}{\mddefault}{\updefault}21}}}}}
\put(1800,912){\makebox(0,0)[lb]{\smash{{{\SetFigFont{5}{6.0}{\rmdefault}{\mddefault}{\updefault}28}}}}}
\put(2400,912){\makebox(0,0)[lb]{\smash{{{\SetFigFont{5}{6.0}{\rmdefault}{\mddefault}{\updefault}34}}}}}
\put(3000,912){\makebox(0,0)[lb]{\smash{{{\SetFigFont{5}{6.0}{\rmdefault}{\mddefault}{\updefault}21}}}}}
\put(3600,912){\makebox(0,0)[lb]{\smash{{{\SetFigFont{5}{6.0}{\rmdefault}{\mddefault}{\updefault}7}}}}}
\put(4200,912){\makebox(0,0)[lb]{\smash{{{\SetFigFont{5}{6.0}{\rmdefault}{\mddefault}{\updefault}1}}}}}
\put(1500,2412){\makebox(0,0)[lb]{\smash{{{\SetFigFont{5}{6.0}{\rmdefault}{\mddefault}{\updefault}1}}}}}
\put(1500,1212){\makebox(0,0)[lb]{\smash{{{\SetFigFont{5}{6.0}{\rmdefault}{\mddefault}{\updefault}6}}}}}
\put(2100,1212){\makebox(0,0)[lb]{\smash{{{\SetFigFont{5}{6.0}{\rmdefault}{\mddefault}{\updefault}1}}}}}
\put(1800,612){\makebox(0,0)[lb]{\smash{{{\SetFigFont{5}{6.0}{\rmdefault}{\mddefault}{\updefault}7}}}}}
\put(2400,612){\makebox(0,0)[lb]{\smash{{{\SetFigFont{5}{6.0}{\rmdefault}{\mddefault}{\updefault}1}}}}}
\put(1500,2712){\makebox(0,0)[lb]{\smash{{{\SetFigFont{5}{6.0}{\rmdefault}{\mddefault}{\updefault}3}}}}}
\put(1800,1812){\makebox(0,0)[lb]{\smash{{{\SetFigFont{5}{6.0}{\rmdefault}{\mddefault}{\updefault}1}}}}}
\end{picture}
}

%% file: xfig/millef2.eepic
\setlength{\unitlength}{0.00033333in}
\begingroup\makeatletter\ifx\SetFigFont\undefined%
\gdef\SetFigFont#1#2#3#4#5{%
  \reset@font\fontsize{#1}{#2pt}%
  \fontfamily{#3}\fontseries{#4}\fontshape{#5}%
  \selectfont}%
\fi\endgroup%
{\renewcommand{\dashlinestretch}{30}
\begin{picture}(12774,4981)(0,-10)
\path(10212,4654)(12762,4654)(12762,1054)
	(10212,1054)(10212,4654)
\path(10512,4654)(10512,3754)(10737,1954)(10737,1054)
\path(10737,4654)(10737,3754)(10632,3079)
\path(10587,2854)(10437,1954)(10437,1054)
\path(11112,4654)(11112,1054)
\path(11337,4654)(11337,1054)
\path(11712,4654)(11712,3754)
\path(11712,3754)(12297,1954)(12297,1054)
\path(12312,4654)(12312,3754)(12042,2914)
\path(11967,2779)(11697,1969)(11697,1054)
\path(11937,4654)(11937,3754)(12087,3229)
\path(12147,3109)(12537,1954)(12537,1054)
\path(12537,4654)(12537,3754)(12282,2899)
\path(12222,2779)(12147,2554)
\path(12087,2434)(11937,1954)(11937,1054)
\path(12,4654)(5112,4654)(5112,1054)
	(12,1054)(12,4654)
\dashline{60.000}(2562,4954)(2562,754)
\path(6312,4654)(8862,4654)(8862,1054)
	(6312,1054)(6312,4654)
\path(6312,4654)(6312,1054)(8487,154)
	(8487,3754)(6312,4654)(6312,4654)
\path(612,4654)(612,3754)(1212,1954)(1212,1054)
\path(1212,4654)(1212,3754)(987,3004)
\path(837,2704)(612,1954)
\path(612,1954)(612,1054)
\path(7617,4114)(7617,3154)(8037,1279)(8037,334)
\path(8112,3904)(8112,2974)(7887,2314)(7902,2329)
\path(7782,2104)(7527,1459)(7527,559)
\path(2412,4639)(2412,3754)(3012,1954)(3012,1054)
\path(3012,4654)(3012,3754)(2757,2914)
\path(2682,2734)(2412,1954)(2412,1054)
\path(6462,4579)(6462,3679)(6312,3379)
	(6762,1954)(6762,1054)
\path(6762,4654)(6762,3754)(6537,3004)
\path(6387,2779)(6312,2554)(6462,1879)(6462,979)
\path(1812,4654)(1812,1054)
\path(7062,4354)(7062,754)
\path(4212,4639)(4212,3739)(4812,1939)(4812,1039)
\path(4812,4639)(4812,3739)(4587,2989)
\path(4212,1969)(4212,1069)
\path(4437,2734)(4212,1984)
\path(3612,4639)(3612,1039)
\dottedline{45}(8067,4639)(8067,3739)(8667,1939)(8667,1039)
\dottedline{45}(8622,4654)(8622,3754)(8397,3004)
\dottedline{45}(8067,1954)(8067,1054)
\dottedline{45}(8292,2719)(8067,1969)
\dottedline{45}(7407,4639)(7407,1039)
\thicklines
\path(312,754)(311,753)(310,752)
	(309,750)(307,747)(304,743)
	(300,739)(296,733)(291,726)
	(286,718)(279,709)(273,698)
	(266,687)(258,674)(251,661)
	(243,646)(236,631)(229,616)
	(223,599)(217,582)(213,565)
	(209,548)(206,530)(205,512)
	(206,494)(208,476)(213,458)
	(219,440)(228,422)(240,404)
	(254,386)(272,368)(293,350)
	(317,333)(346,315)(379,297)
	(416,280)(458,263)(506,246)
	(559,229)(618,213)(683,197)
	(753,182)(830,168)(912,154)
	(983,144)(1057,134)(1132,125)
	(1208,117)(1285,110)(1361,103)
	(1436,97)(1509,91)(1581,86)
	(1651,82)(1719,78)(1784,74)
	(1847,71)(1908,69)(1966,66)
	(2022,64)(2076,62)(2127,61)
	(2177,59)(2224,58)(2270,57)
	(2315,56)(2358,55)(2400,55)
	(2441,55)(2482,54)(2522,54)
	(2562,54)(2602,54)(2642,54)
	(2683,55)(2724,55)(2766,55)
	(2809,56)(2854,57)(2900,58)
	(2947,59)(2997,61)(3048,62)
	(3102,64)(3158,66)(3216,69)
	(3277,71)(3340,74)(3405,78)
	(3473,82)(3543,86)(3615,91)
	(3688,97)(3763,103)(3839,110)
	(3916,117)(3992,125)(4067,134)
	(4141,144)(4212,154)(4294,168)
	(4371,182)(4441,197)(4506,213)
	(4565,229)(4618,246)(4666,263)
	(4708,280)(4745,297)(4778,315)
	(4807,333)(4831,350)(4852,368)
	(4870,386)(4884,404)(4896,422)
	(4905,440)(4911,458)(4916,476)
	(4918,494)(4919,512)(4918,530)
	(4915,548)(4911,565)(4907,582)
	(4901,599)(4895,616)(4888,631)
	(4881,646)(4873,661)(4866,674)
	(4858,687)(4851,698)(4845,709)
	(4838,718)(4833,726)(4828,733)
	(4824,739)(4820,743)(4817,747)(4812,754)
\path(4906.161,673.789)(4812.000,754.000)(4857.337,638.915)
\path(8637,154)(8639,154)(8643,155)
	(8651,156)(8662,158)(8678,160)
	(8696,163)(8718,167)(8741,172)
	(8766,177)(8792,184)(8818,191)
	(8844,199)(8869,208)(8895,219)
	(8920,231)(8945,246)(8969,263)
	(8992,282)(9012,304)(9029,327)
	(9042,349)(9052,369)(9059,387)
	(9064,402)(9067,413)(9069,423)
	(9069,430)(9069,436)(9068,442)
	(9067,447)(9065,454)(9063,462)
	(9060,472)(9056,485)(9051,502)
	(9045,522)(9036,547)(9026,574)
	(9012,604)(8996,634)(8978,662)
	(8958,689)(8939,713)(8918,736)
	(8898,757)(8877,777)(8856,795)
	(8836,813)(8815,829)(8795,844)
	(8776,858)(8759,871)(8744,881)
	(8732,890)(8712,904)
\path(8827.512,859.761)(8712.000,904.000)(8793.104,810.608)
\end{picture}
}

%% file: xfig/Li2.eepic
\setlength{\unitlength}{0.00033333in}
\begingroup\makeatletter\ifx\SetFigFont\undefined%
\gdef\SetFigFont#1#2#3#4#5{%
  \reset@font\fontsize{#1}{#2pt}%
  \fontfamily{#3}\fontseries{#4}\fontshape{#5}%
  \selectfont}%
\fi\endgroup%
{\renewcommand{\dashlinestretch}{30}
\begin{picture}(8124,6084)(0,-10)
\texture{8101010 10000000 444444 44000000 11101 11000000 444444 44000000 
	101010 10000000 444444 44000000 10101 1000000 444444 44000000 
	101010 10000000 444444 44000000 11101 11000000 444444 44000000 
	101010 10000000 444444 44000000 10101 1000000 444444 44000000 }
\shade\path(1212,5754)(3612,5754)(3612,3804)
	(3612,54)(1212,54)(1212,5004)(1212,5754)
\path(1212,5754)(3612,5754)(3612,3804)
	(3612,54)(1212,54)(1212,5004)(1212,5754)
\drawline(4212,1104)(4212,1104)
\path(12,5754)(8112,5754)
\path(12,5754)(8112,5754)
\path(12,54)(8112,54)
\path(12,54)(8112,54)
\thicklines
\path(4812,5754)(4812,54)
\path(5412,54)(5412,5754)
\path(6012,5754)(6012,54)
\path(6612,5754)(6612,54)
\path(4212,1104)(4212,54)
\path(612,2904)(4212,1104)
\path(627,2874)(627,3174)(1137,3429)
\path(4212,5004)(3687,4749)
\path(4212,5004)(4212,5754)
\put(4137,5904){\makebox(0,0)[lb]{\smash{{{\SetFigFont{5}{6.0}{\rmdefault}{\mddefault}{\updefault}$i$}}}}}
\end{picture}
}